\documentclass[oneside,leqno,11pt]{article}
\usepackage{amssymb,amsmath,latexsym}

\def\ge{\mathfrak{e}}

\def\gg{\mathfrak{g}}

\def\gk{\mathfrak{k}}
\def\gl{\mathfrak{l}}

\def\gq{\mathfrak{q}}

\def\gs{\mathfrak{s}}
\def\gt{\mathfrak{t}}
\def\gu{\mathfrak{u}}

\def\gso{\mathfrak{s}\mathfrak{o}}
\def\gsu{\mathfrak{s}\mathfrak{u}}
\def\gsp{\mathfrak{s}\mathfrak{p}}

\def\C{\mathbb{C}}

\def\E{\mathbb{E}}

\def\H{\mathbb{H}}

\def\K{\mathbb{K}}
\def\L{\mathbb{L}}

\def\N{\mathbb{N}}

\def\R{\mathbb{R}}

\def\Ind{{\rm Ind}}
\def\Ad{{\rm Ad}}

\def\rank{{\rm rank\,}}

\def\det{{\rm det\,}}
\def\Pf{{\rm Pf\,}}

\def\cA{\mathcal{A}}

\def\cF{\mathcal{F}}

\def\cO{\mathcal{O}}

\def\cS{\mathcal{S}}

\renewcommand{\thesection}{\arabic{section}}
\renewcommand{\theequation}{\thesection.\arabic{equation}}

\newtheorem{theorem}[equation]{Theorem}

\newtheorem{lemma}[equation]{Lemma}

\newtheorem{corollary}[equation]{Corollary}

\newtheorem{proposition}[equation]{Proposition}

\newtheorem{definition}[equation]{Definition}

\newtheorem{remark}[equation]{Remark}

\newtheorem{noname}[equation]{}

\title{Geometry of the Borel -- de Siebenthal Discrete Series}

\author{Bent \O{}rsted\footnote{Research partially supported by the Danish 
Research Council.}\,\,\,\, 
\& \,\,Joseph A. Wolf\footnote{Research partially supported by the NSF,
the Danish Research Council, and IMF, Aarhus University. 
\endgraf
 \endgraf
}}

\date{27 January 2009}

\begin{document}

\maketitle

\abstract{Let $G_0$ be a connected, simply connected real simple Lie group.  
Suppose that $G_0$ has a compact Cartan subgroup $T_0$, so it has discrete
series representations.  Relative to $T_0$ there is a distinguished
positive root system $\Delta^+$ for which there is a unique noncompact
simple root $\nu$, the ``Borel -- de Siebenthal system''.  There is a lot
of fascinating geometry associated to the corresponding ``Borel -- de Siebenthal
discrete series'' representations of $G_0$.  In this paper we explore
some of those geometric aspects and we work out the $K_0$--spectra
of the Borel -- de Siebenthal discrete series representations.  This has 
already been carried out in detail for the case where the associated 
symmetric space
$G_0/K_0$ is of hermitian type, i.e. where $\nu$ has coefficient $1$ in the 
maximal root $\mu$, so we assume that the group $G_0$ is not of hermitian
type, in other words that $\nu$ has coefficient $2$ in $\mu$.  
\medskip

Several authors have studied the case where $G_0/K_0$ is a
quaternionic symmetric space and the inducing holomorphic vector
bundle is a line bundle.  That is the case where $\mu$ is
orthogonal to the compact simple roots and the inducing representation
is $1$--dimensional.}
\bigskip

\setcounter{section}{0}
\section{Introduction}\label{sec1}
\setcounter{equation}{0}

One of Harish--Chandra's great achievements was the
existence theorem for discrete series representations
of a semisimple Lie group.  He characterized the groups
with discrete series representations by the equal rank
condition, he found the explicit formulae on the regular
elliptic set for the characters of the discrete series,
and he showed that those formulae specify the characters.
At the same time (and as a main motivation) he
was able to explicitly construct a particularly
simple series, the holomorphic discrete
series, for those groups where the corresponding
Riemannian symmetric space is a bounded symmetric
domain in a complex Euclidean space. For the other
discrete series representations however, the
actual construction has remained less explicit,
although there are several beautiful realizations.
\medskip

In this paper we initiate the study of a certain family, the
so-called Borel -- de Siebenthal discrete series,
from a point of view as close as possible
to that of Harish--Chandra for the holomorphic
discrete series. This is motivated in part by the
work of Gross and Wallach for the scalar case of the quaternionic
discrete series \cite{GW}. As in that case we obtain in
particular the admissibility of the series
for a small subgroup of the maximal compact
subgroup. At the same time we discover a
rather appealing geometry for the coadjoint
orbit that one wants to attach to the
discrete series in question.
In particular we give a detailed classification of the
possible structures of such orbits in terms of explicit
prehomogeneous vector spaces with relative invariants.
We feel these deserve attention in their own right;
and while we do give the construction for the 
Borel -- de Siebenthal discrete
series here, including the explicit $K_0$-types and the
(important) admissibility for a small subgroup $K_1$ of $K_0$,
we defer further analysis of continuations of the series
to a sequel to this paper.  In particular, we shall then
elucidate the role of the relative invariants in
constructing rather singular representations in the
continuation of the series.  We mention that much of this has
been carried out in the quaternion line bundle case in \cite{GW}.
We also mention the papers \cite{Kn2} and \cite{Kn3} treating
such questions for the indefinite orthogonal and symplectic
groups; here methods from \cite{GW} are used, and the connection
to the continuation of unitary modules in the sense of Vogan
(with criteria for good and fair range of unitarity of
cohomologically induced modules) is made clear. Our approach
seeks to employ analytic methods and the geometry of the orbits,
and in particular to use reproducing kernels, see e.g. \cite{WaW}.  
\medskip

Several questions concerning discrete series representations
may be resolved by our methods, for example the question
of finding admissible branching laws, where one obtains direct
sum decompositions with finite multiplicities. By
applying admissibility of $K_1$ such results may be
obtained in complete analogy with what happens for
holomorphic discrete series representations.
\medskip

In this paper we give a complete description of the geometry of the elliptic
coadjoint orbits corresponding to the Borel -- de Siebenthal discrete series.
They are open $G_0$--orbits in certain complex flag manifolds and we give
precise results on their maximal compact subvarieties (which are compact
hermitian symmetric spaces) and the holomorphic normal bundles to those
subvarieties.  We use this structural information to give a concrete geometric 
construction of the representations in this series, including the
structure of the $K_0$--types. Our construction of the $K_0$--types provides 
an analogue of the $K_0$--type decomposition of holomorphic discrete series
representations.  
\medskip

The quaternionic discrete series, studied by Gross and Wallach \cite{GW}
in the line bundle case, is the special case of the Borel -- de Siebenthal 
discrete series, where the maximal root is compact and is orthogonal to 
all the other compact positive roots, or equivalently
(see \cite{W0}) where $K_0$ has a local direct factor isomorphic to $Sp(1)$.
While each complex simple Lie group $G$ has exactly one noncompact real
form $G_0$ that has quaternionic discrete series representations, every
real simple Lie group $G_0$ with $\rank G_0 = \rank K_0$ has either
holomorphic discrete series representations or 
Borel -- de Siebenthal discrete series representations.  Thus the
Borel -- de Siebenthal discrete series is the natural extension of the
holomorphic discrete series.  
\medskip

Our geometric approach allows us to extend several results from \cite{GW}.
While Gross and Wallach constructed quaternionic discrete series
representations on spaces of holomorphic forms with values in a line
bundle, we also allow vector bundles -- which is natural in the more
general setting considered here.  Our construction also provides good 
concrete examples of minimal cohomology degree realizations of discrete
series representations in the sense of Kostant \cite{Ko3}.
\medskip

Our basic tool is complex differential geometry and the associated 
cohomology groups.  An important component of this is a collection of
basic spectral sequence arguments, already implicit in the paper 
\cite{S1}. See also \cite{S2}, \cite{S3} and \cite{W4}.  Here we make use of
some technical results of M. Eastwood and the second named author from 
\cite{EW} for some crucial identifications of duals of finite dimensional  
representations of reductive Lie groups, in particular for keeping track 
of the action of the center in terms of the highest weights and 
the Dynkin diagrams. 
\medskip

Our results include a careful collection of the data attached to
the orbits in question, and an explicit formula for the
$K_0$--types in the Harish--Chandra module corresponding to
the discrete series in question. As a by--product we find two natural
sets of strongly orthogonal roots, one corresponding to the
hermitian symmetric space $K_0/L_0$ and the other corresponding to
the riemannian symmetric space $G_0/K_0$.  They fit together to realize
the orbit $G_0(z_0) = G_0/L_0$ as a kind of Siegel domain
of Type II. This should provide useful coordinates for explicit
calculations of the elements in the Harish--Chandra module.
\medskip

In Section \ref{sec2} we work out the general structure of the complex
manifold $D = G_0(z_0) \cong G_0/L_0$.  We describe the action of
$\gl_0$ on the tangent and normal spaces to the maximal compact 
subvariety $Y = K_0(z_0) \cong K_0/L_0$, and on their duals.  then
we discuss a negativity condition that is crucial to the realization of 
our discrete series representations.
\medskip

In Section \ref{sec3} we list all instances of simple Lie algebras $\gg_0$
corresponding to Borel -- de Siebenthal root orders.  Setting aside the 
well--understood hermitian symmetric cases, we work out the precise 
structure of the algebras $\gk_0$ and $\gl_0$ and the parts of the complexified
Lie algebra $\gg$ that correspond to the holomorphic tangent space of $Y$ and
the holomorphic normal space of $Y$ in $D$, including the representations
of $\gl_0$ on those two spaces.  In each case this allows explicit 
parameterization of the  Borel -- de Siebenthal discrete series.
\medskip

In Section \ref{sec4}  we consider the prehomogeneous space $(L,\gu_1)$
where $\gu_1 \subset \gg$ represents the holomorphic normal space to
$Y$ in $D$.  There we describe the algebra of relative invariants, using
our knowledge of the representation of $L$ on the symmetric algebra 
$S(\gu_1)$.  In most cases we can be explicit, but in some we must rely
on general results of Sato and Kimura \cite{SK}. These invariants
are (in addition to being interesting in themselves) relevant for
the next step of understanding the analytic continuation of the
discrete series; here the ring of regular functions on the zero
set of an invariant will correspond to a module in this continuation.
We intend to follow this idea in a sequel to this paper.   
\medskip

In Section \ref{sec5} we assemble our preparations and work out the exact
$K$--spectrum of the Borel -- de Siebenthal discrete series representations.
Our main result here, which is the main result of the paper, is
Theorem \ref{bds_line_ds}.
In a final example we look a the "sufficient negativity" condition that
ensures the non-vanishing of exactly the right analytic cohomology group,
and we compare it to the corresponding condition for individual $K$-types
- this will indicate a possibility of continuing the discrete series family.
      
\medskip

\section{Notation and the Basic Fibration}\label{sec2}
\setcounter{equation}{0}
In general we use capital Latin letters for Lie groups with subscript
${}_0$ for real groups and no subscript for complexifications.  
We use the corresponding small Gothic letters for Lie algebras, again
with subscript 
${}_0$ for real Lie algebras and no subscript for complexifications.
Our basic objects are a connected simply connected simple real
Lie group $G_0$, its Lie algebra $\gg_0$, the complexification $G$ of
$G_0$, and the Lie algebra $\gg$ of $G$.  Here $G$ is a connected
simply connected complex Lie group and the inclusion 
$\gg_0 \hookrightarrow \gg$ defines a homomorphism $G_0 \to G$ with
discrete central kernel.
\medskip

When we omit a subscript ${}_0$ where there had been one before, we
mean complexification.
\medskip

Fix a Cartan involution $\theta$ of $G_0$ and $\gg_0$.  The fixed point
set $K_0 = G_0^\theta$ is a maximal compactly embedded subgroup of $G_0$.
As usual, we decompose $\gg_0 = \gk_0 + \gs_0$ and $\gg = \gk + \gs$
into $(\pm 1)$--eigenspaces of $\theta$, where $\gk_0$ (resp. $\gk$)
is the Lie algebra of $K_0$ (resp. $K$).
\medskip

We now make two assumptions:
\begin{equation}\label{basic1}
\rank G_0 = \rank K_0 \text{\,, and the symmetric space }
S_0 := G_0/K_0 \text{ is not of hermitian type.}
\end{equation}
In particular $K_0$ is a maximal compact subgroup of $G_0$.  Both $G_0$ and
$K_0$ are simply connected semisimple groups with finite center.
\medskip

Fix a maximal torus $T_0 \subset K_0$.  Then $T_0$ is a compact Cartan
subgroup of $G_0$, and a celebrated theorem of Harish--Chandra says
that $G_0$ has discrete series representations.
\medskip

The construction of Borel and de Siebenthal \cite{BoS} provides a
positive root system $\Delta^+ = \Delta^+_G$ for $(\gg,\gt)$ such that the
associated simple root system $\Psi = \Psi_G$ contains just one noncompact
root.  We denote
\begin{equation}\label{noncsimple}
\Psi = \{\psi_1, \dots , \psi_\ell\} \text{ (Bourbaki root order) and }
\nu \in \Psi \text{ is the noncompact simple root.}
\end{equation}
Every root $\alpha \in \Delta^+$ has expression 
$\alpha = \sum n_i(\alpha)\psi_i$.  Since we have excluded the 
hermitian case, the coefficient of $\nu$ in the maximal root $\mu$
is $2$.  Further, a root is compact just when the coefficient of
$\nu$ in its expansion is $0$ or $\pm 2$, noncompact just when that
coefficient is $\pm 1$.  Also, $(\Psi \setminus \{\nu\})\cup\{-\mu\}$
is a simple root system for $(\gk,\gt)$.  Grading by the coefficient
$n_\nu$ of $\nu$ we have a parabolic subalgebra of $\gg$ given by
\begin{equation}\label{parabolic}
\gq = \gl + \gu_-, \text{ reductive part } 
	\gl = \gt + \sum_{n_\nu = 0} \gg_\alpha \text{ and nilradical }
	\gu_- = \gu_{-2} + \gu_{-1} 
\end{equation}
where $\gu_i = \sum_{n_\nu = i} \gg_\alpha$.  The opposite parabolic
is $\gq^{opp} = \gl + \gu_+$ where $\gu_+ = \gu_1 + \gu_2$.  Note that
\begin{equation}
\gs = \gu_{-1} + \gu_1 \text{, so } \gq \cap \gs = \gu_{-1} \text{, and }
\gk = \gu_{-2} + \gl + \gu_2.
\end{equation}

On the group level, we have the parabolic subgroup $Q \subset G$
where $Q$ has Lie algebra $\gq$.  The group $Q$ has Chevalley 
semidirect product decomposition $LU_-$ where $L$ is the reductive
component and $U_-$ is the unipotent radical.  Note that $G_0 \cap Q$
is a real form $L_0$ of $L$ and that $L_0$ is the centralizer in
$K_0$ of a circle subgroup of $T_0$.  The parabolic $Q$ defines a
complex flag manifold $Z = G/Q$, say with base point $z_0 = 1Q$,
and and open orbit $D = G_0(z_0) \cong G_0/L_0$.  The complex manifold
$D$ has maximal compact subvariety $Y = K_0(z_0) \cong K_0/L_0$, 
which is a smaller complex flag manifold $K/(K\cap Q)$.  
\medskip

Our choice of signs in (\ref{parabolic}) is such that
\begin{equation}\label{tanspaces}
\begin{aligned}
&\gu_+ \text{ is the holomorphic tangent space of } D \text{ at } z_0\, , \\
&\gu_2 \text{ is the holomorphic tangent space of } Y \text{ at } z_0\, , 
\text{ and } \\
&\gu_1 \text{ is the holomorphic normal space of } Y \text{ in } D
\text{ at } z_0\, .
\end{aligned}
\end{equation}
Since $G_0/K_0$ is irreducible but not hermitian we know that the action
of $\gk_0$ on $\gs_0$ is absolutely irreducible.  Thus the action of $\gk$ on
$\gs$ is irreducible.  From \cite[Theorem 8.13.3]{W2} we know that
the action of $\gl$ on each $\gu_i$ is irreducible.
It will be convenient to have the notation
\begin{equation}\label{repl}
\tau_i: \text{ representation of } L \text{ on the vector space } \gu_i.
\end{equation}
The contragredient (dual) of $\tau_i$ is $\tau_i^* = \tau_{-i}$.
Some obvious highest or lowest weight spaces of the $\tau_i$ are given by
\begin{equation}\label{wtspaces}
\begin{aligned}
& \tau_2 \text{ has highest weight space } \gg_\mu \text{ and }
	 \tau_{-2}  \text{ has lowest weight space } \gg_{-\mu}, \\
& \tau_1 \text{ has lowest weight space } \gg_\nu \text{ and }
	\tau_{-1} \text{ has highest weight space } \gg_{-\nu}\,.
\end{aligned}
\end{equation}
Note that the degree $\deg \tau_i = \dim_\C \gu_i$.  If $i \ne 0$ it is 
the number of roots $\alpha$ such that $n_\nu(\alpha) = i$.
\medskip

The basic tool in this paper is the real analytic fibration
\begin{equation}\label{fibration}
D \to S_0 \text{ with fiber } Y, \text{ in other words }
G_0/L_0 \to G_0/K_0 \text{ with fiber } K_0/L_0.
\end{equation}
The structure of the holomorphic tangent bundle and the holomorphic
normal bundle to $Y$ in $D$ is given by (\ref{tanspaces}),
(\ref{repl}) and (\ref{wtspaces}).  In the next section we will make
this explicit.  The fibration (\ref{fibration}) was first considered 
by W. Schmid in \cite{S1} and \cite{S2} for a related situation in 
which $L_0 = T_0$, and then somewhat later by R. O. Wells and one of us
\cite{WeW} without that restriction.  A much more general setting,
which drops the compactness assumption on $L_0$, is that of the
double fibration transform (see \cite{FHW} and the references there),
where $S_0$ is replaced by a complexification $S \subset G/K$.
Specialization of the double fibration transform to the Borel -- de Siebenthal
setting is carried out in \cite{EW}.
\medskip

The simple root system $\Psi = \{\psi_1, \dots , \psi_\ell\}$ of $\gg$ 
defines the system
\begin{equation}\label{fund-wts}
\Xi = \{\xi_1, \dots , \xi_\ell\} \text{ where }
        \tfrac{2\langle \xi_i,\psi_j\rangle}{\langle\psi_j,\psi_j\rangle}
        = \delta_{i,j}
\end{equation}
of fundamental simple weights.  Let $\gamma$ be the highest weight of 
an irreducible representation of $L_0$.  For our discussion of the
Borel -- de Siebenthal discrete series in Section \ref{sec5} 
we will need to know exactly when
$\langle \gamma + \rho_\gg, \alpha \rangle < 0$ for all positive complementary 
roots $\alpha$ (roots $\alpha$ that are not roots of $\gl$) where $\rho_\gg$
denotes half the sum of the positive roots of $\gg$.  The condition is 
Theorem \ref{delta} below.
\medskip

Define $\nu^*$ by $\langle \nu^*, \psi_j \rangle = 0$ for $\psi_j \ne \nu$ and
2$\langle \nu^*, \nu \rangle = \langle \nu, \nu \rangle$ (i.e. the
fundamental weight dual to $\nu$).  Then $\gamma\in i\gt_0^*$ 
decomposes as
\begin{equation}\label{lambdanought}
\gamma = \gamma_0 + t\nu^* \text{ where } 
	\langle \gamma_0 , \nu \rangle = 0 \text{ and } t \in \R.
\end{equation}
Define $\Delta_i = \{\alpha \in \Delta_G \mid n_\nu(\alpha) = i\}$,
so $\gu_i = \sum_{\alpha \in \Delta_i} \gg_\alpha$ for $i \in \{\pm 1, \pm 2\}$.
Thus the positive root system decomposes as 
$\Delta^+ = (\Delta_0 \cap \Delta^+) \cup \Delta_1 \cup \Delta_2$.
The highest weight of $\tau_2$, representation of $\gl$ on
$\gu_2 = \sum_{\Delta_2} \gg_\alpha$ is $\mu$.  If we subtract a positive
combination of roots of $\Psi \setminus \{\nu\}$ from $\mu$ we decrease
the inner product with $\gamma + \rho_\gg$.  Thus
\addtocounter{equation}{1}
\begin{equation}\tag{\theequation$a$} \label{delta1}
\langle \gamma + \rho_\gg, \alpha \rangle < 0 \text{ for all }
\alpha \in \Delta_2 \text{ if and only if }
\langle \gamma + \rho_\gg, \mu \rangle < 0.
\end{equation}
The highest weight of $\tau_{-1}$ is $-\nu$, so $\tau_1$ has highest weight
$w^0_\gl(\nu)$ where $w^0_\gl$ is the longest element of the Weyl
group of $\gl$.  Thus
\begin{equation}\tag{\theequation$b$}\label{delta2}
\langle \gamma + \rho_\gg, \alpha \rangle < 0 \text{ for all }
\alpha \in \Delta_1\text{ if and only if }
\langle \gamma + \rho_\gg, w^0_\gl(\nu)\rangle < 0.
\end{equation}
As $\nu^*$ is orthogonal to the roots of $\gl$ it is fixed
by the inverse of $w^0_\gl$, so
$\langle \nu^*,w^0_\gl(\nu)\rangle = \langle \nu^*,\nu\rangle = 1$.  
Using the decomposition (\ref{lambdanought}), and combining
(\ref{delta1}) and (\ref{delta2}), we have
\begin{theorem}\label{delta}
The following conditions are equivalent.

{\rm 1.} The inequality
$\langle \gamma + \rho_\gg, \alpha \rangle < 0$ holds for every root
$\alpha \in \Delta_1 \cup \Delta_2$  $($i.e. every positive complementary
root$)$ 

{\rm 2.} Both
$t < -\tfrac{1}{2}\langle \gamma_0 + \rho_\gg, \mu \rangle$  and 
$t < -\langle \gamma_0 + \rho_\gg, w^0_\gl(\nu)\rangle$.
\end{theorem}

\begin{remark} \label{charbds} {\rm In Theorem \ref{delta} it is automatic 
that $\langle \gamma + \rho_\gg, \beta \rangle > 0$ for every positive root 
of $\gl$, so the conditions of Theorem \ref{delta} ensure that
$\langle \gamma + \rho_\gg, \alpha \rangle \ne 0$ for every root $\alpha$,
in other words that $\gamma + \rho_\gg$ is the Harish--Chandra parameter
of a discrete series representation of $G_0$.  Specifically, in our setting, 
the conditions of Theorem \ref{delta} will characterize the 
Borel -- de Siebenthal discrete series.}
\hfill $\diamondsuit$
\end{remark}

\section{Classification}\label{sec3}
\setcounter{equation}{0}

In this section we give a complete list the simple Lie algebras $\gg_0$
for which the hypotheses (\ref{basic1}) hold.  We then specify the
complex parabolic subalgebra $\gq \subset \gg$ and the real subalgebras 
$\gk_0$ and $\gl_0$ of $\gg_0$.  Next, we give precise descriptions of the
representations $\tau_i$ and their representation spaces $\gu_i$.  Much
of this is done using the Dynkin diagrams to indicate highest weights
of representations.  There the special cases, where $G_0/K_0$ is a
quaternionic symmetric space, are visible at a glance: they are the ones
where $-\mu$ connects directly to $\nu$ in the extended Dynkin diagram of 
$\gg$. 
\medskip

We will denote highest weights of representations as follows.  In Dynkin
diagrams with two root lengths we denote short root nodes by black dots 
$\bullet$ and long roots by the usual circles 
\setlength{\unitlength}{.5 mm}
\begin{picture}(1,2)
\put(0,2.5){\circle{2}}
\end{picture}.  
Extended
diagrams are those with the negative of the maximal root $\mu$ attached by
the usual rules.  Recall the system $\Xi$ of fundamental simple
weights from (\ref{fund-wts}).  The (irreducible finite
dimensional) representation of $G$ and $\gg$ of highest weight
$\sum n_i\xi_i$ is indicated by the Dynkin diagram of $\gg$ with
$n_i$ written next to the $i^{th}$ node, except that we omit writing
zeroes.  So for example the adjoint representations are indicated by
$$
\begin{tabular}{|c|l||c|l|}\hline
$A_\ell\,, \ell \geqq 1$ &
\setlength{\unitlength}{.5 mm}
\begin{picture}(100,12)
\put(5,2){\circle{2}}
\put(4,5){$2$}
\put(10,2){or}
\put(20,2){\circle{2}}
\put(19,5){$1$}
\put(21,2){\line(1,0){13}}
\put(35,2){\circle{2}}
\put(36,2){\line(1,0){13}}
\put(52,2){\circle*{1}}
\put(55,2){\circle*{1}}
\put(58,2){\circle*{1}}
\put(66,2){\line(1,0){13}}
\put(80,2){\circle{2}}
\put(79,5){$1$}
\end{picture} 
&
$B_\ell$\,, $\ell \geqq 3$ &
\setlength{\unitlength}{.5 mm}
\begin{picture}(90,12)
\put(5,2){\circle{2}}
\put(6,2){\line(1,0){13}}
\put(20,2){\circle{2}}
\put(19,5){$1$}
\put(21,2){\line(1,0){13}}
\put(37,2){\circle*{1}}
\put(40,2){\circle*{1}}
\put(43,2){\circle*{1}}
\put(46,2){\line(1,0){13}}
\put(60,2){\circle{2}}
\put(61,2.5){\line(1,0){13}}
\put(61,1.5){\line(1,0){13}}
\put(75,2){\circle*{2}}
\end{picture} 
\\
\hline 
$C_\ell$\,, $\ell \geqq 2$ &
\setlength{\unitlength}{.5 mm}
\begin{picture}(100,10)
\put(5,2){\circle*{2}}
\put(4,5){$2$}
\put(6,2){\line(1,0){13}}
\put(20,2){\circle*{2}}
\put(21,2){\line(1,0){13}}
\put(37,2){\circle*{1}}
\put(40,2){\circle*{1}}
\put(43,2){\circle*{1}}
\put(46,2){\line(1,0){13}}
\put(60,2){\circle*{2}}
\put(61,2.5){\line(1,0){13}}
\put(61,1.5){\line(1,0){13}}
\put(75,2){\circle{2}}
\end{picture} 
&
$D_\ell$\,, $\ell \geqq 4$ &
\setlength{\unitlength}{.5 mm}
\begin{picture}(90,20)
\put(5,9){\circle{2}}
\put(6,9){\line(1,0){13}}
\put(20,9){\circle{2}}
\put(19,12){$1$}
\put(21,9){\line(1,0){13}}
\put(37,9){\circle*{1}}
\put(40,9){\circle*{1}}
\put(43,9){\circle*{1}}
\put(46,9){\line(1,0){13}}
\put(60,9){\circle{2}}
\put(61,8.5){\line(2,-1){13}}
\put(75,2){\circle{2}}
\put(61,9.5){\line(2,1){13}}
\put(75,16){\circle{2}}
\end{picture} 
\\
\hline
$G_2$ &
\setlength{\unitlength}{.5 mm}
\begin{picture}(100,10)
\put(5,2){\circle*{2}}
\put(6,1){\line(1,0){13}}
\put(6,2){\line(1,0){13}}
\put(6,3){\line(1,0){13}}
\put(20,2){\circle{2}}
\put(19,5){$1$} 
\end{picture} 
&
$F_4$ &
\setlength{\unitlength}{.5 mm}
\begin{picture}(90,11)
\put(5,2){\circle{2}}
\put(4,4){$1$}
\put(6,2){\line(1,0){13}}
\put(20,2){\circle{2}}
\put(21,1.5){\line(1,0){13}}
\put(21,2.5){\line(1,0){13}}
\put(35,2){\circle*{2}}
\put(36,2){\line(1,0){13}}
\put(50,2){\circle*{2}}
\end{picture} 
\\
\hline
$E_6$ &
\setlength{\unitlength}{.5 mm}
\begin{picture}(100,20)
\put(5,16){\circle{2}}
\put(6,16){\line(1,0){13}}
\put(20,16){\circle{2}}
\put(21,16){\line(1,0){13}}
\put(35,16){\circle{2}}
\put(36,16){\line(1,0){13}}
\put(50,16){\circle{2}}
\put(51,16){\line(1,0){13}}
\put(65,16){\circle{2}}
\put(35,15){\line(0,-1){13}}
\put(35,1){\circle{2}}
\put(38,1){$1$}
\end{picture}  
&
$E_7$ & 
\setlength{\unitlength}{.5 mm}
\begin{picture}(90,20)
\put(5,16){\circle{2}}
\put(4,8){$1$}
\put(6,16){\line(1,0){13}}
\put(20,16){\circle{2}}
\put(21,16){\line(1,0){13}}
\put(35,16){\circle{2}}
\put(36,16){\line(1,0){13}}
\put(50,16){\circle{2}}
\put(51,16){\line(1,0){13}}
\put(65,16){\circle{2}}
\put(66,16){\line(1,0){13}}
\put(80,16){\circle{2}}
\put(35,15){\line(0,-1){13}}
\put(35,1){\circle{2}}
\end{picture} 
\\
\hline
$E_8$ &
\setlength{\unitlength}{.5 mm}
\begin{picture}(100,20)
\put(5,16){\circle{2}}
\put(6,16){\line(1,0){13}}
\put(20,16){\circle{2}}
\put(21,16){\line(1,0){13}}
\put(35,16){\circle{2}}
\put(35,16){\line(1,0){13}}
\put(50,16){\circle{2}}
\put(51,16){\line(1,0){13}}
\put(65,16){\circle{2}}
\put(66,16){\line(1,0){13}}
\put(80,16){\circle{2}}
\put(81,16){\line(1,0){13}}
\put(95,16){\circle{2}}
\put(94,8){$1$}
\put(35,15){\line(0,-1){13}}
\put(35,1){\circle{2}}
\end{picture}  & &
\\
\hline
\end{tabular}
$$

We will use this notation for $\gk$ as well.  It can be identified from
its simple 
root system $\Psi_\gk = (\Psi \setminus \{\nu\})\cup\{-\mu\}$.  
We'll follow the
notation of \cite{BE}, except that we won't darken the dots.  Thus the
diagram of $\gl$ consists of the diagram of $\gg$, except that 
the {\footnotesize o} (resp. $\bullet$) at the node
for the noncompact simple root $\nu$ is replaced by an $\times$ 
(resp. {\scriptsize $\boxtimes$})\footnote{Of course one can also 
look at $L_0$ as a subgroup of $K_0$, and from that viewpoint the
diagram of $\gl$ is obtained from that of $\gk$ on replacing the
{\footnotesize o} at the node for $-\mu$ with a $\times$.  However it is
more convenient to look at $L_0$ as a subgroup of $G_0$, and the diagram
of $\gl$ from that viewpoint, when we consider the action of $L$ on
the subspaces $\gu_i$.}.  In the diagram of $\gl$, a
symbol $\times$ or {\scriptsize $\boxtimes$} indicates the $1$--dimensional 
center of $L$.
The irreducible representation of $L$ with highest weight 
$\sum n_i\xi_i$ now indicated by the Dynkin diagram of $\gl$ with
$n_i$ written next to the $i^{th}$ node, for $n_i \ne 0$, whether that node 
is {\footnotesize o}, $\bullet$, $\times$ or {\scriptsize $\boxtimes$}. 
\medskip

If $\nu^*$ is the fundamental simple weight corresponding to the noncompact
simple root $\nu$, and $x \in \gt$ by 
$\alpha(x) = \langle\xi,\alpha\rangle$ for $\alpha \in \gt^*$, then
$\gl_0$ has center $i\R\nu^*$.
\medskip

Now we use the fact that $K_0$ is connected, simply connected and semisimple.
The simple root system $\Psi_\gk$ decomposes into 
$$
\begin{aligned}
&\Psi_{\gk_1}: \text{ connected component that contains } -\mu, \text{ and } \\
&\Psi_{\gk_2}: \text{ the complement of } \Psi_{\gk_1} \text{ in } \Psi_\gk
\end{aligned}
$$
This results in decompositions $K_0 = K_1 \times K_2$ and 
$L_0 = L_1 \times L_2$, which we make explicit in each case.
\medskip

In the following we list the Dynkin diagrams, with the possibilities of
the noncompact simple root $\nu$ among the simple roots $\psi_i$.  Also
in the picture one finds the extended Dynkin diagram node for $-\mu$
where $\mu$ is the maximal root.  Diagrams of Type $A$ do not occur because
$\nu$ must have coefficient $2$ in the expression of 
$\mu$ as a linear combination of simple roots.  We now consider the cases 
where $\gg$ is of type $B$.
\medskip

\begin{noname} \label{so-odd-1}{\sc  Case $Spin(4,2\ell - 3)$.} {\rm
Here $G_0$ is the $2$--sheeted cover of the group $SO(4,2\ell - 3)$
which is a real analytic subgroup of the complex simply connected
group $Spin(2\ell + 1;\C)$.  Its extended Dynkin diagram is
\begin{equation}\tag{\theequation$a$}
\setlength{\unitlength}{.75 mm}
\begin{picture}(140,15)
\put(20,10){\circle{2}}
\put(25,0){\circle{2}}
\put(18,5){$\psi_1$}
\put(22,-3){$-\mu$}
\put(26,1){.}
\put(27,2){.}
\put(28,3){.}
\put(29,4){.}
\put(30,5){.}
\put(31,6){.}
\put(32,7){.}
\put(33,8){.}
\put(21,10){\line(1,0){13}}
\put(35,10){\circle{2}}
\put(33,5){$\psi_2$}
\put(33,12){$\nu$}
\put(36,10){\line(1,0){13}}
\put(52,10){\circle*{1}}
\put(55,10){\circle*{1}}
\put(58,10){\circle*{1}}
\put(61,10){\line(1,0){13}}
\put(75,10){\circle{2}}
\put(73,5){$\psi_{\ell - 1}$}
\put(76,10.5){\line(1,0){13}}
\put(76,9.5){\line(1,0){13}}
\put(90,10){\circle*{2}}
\put(89,5){$\psi_\ell$}
\put(110,10){(type $B_\ell$\,, $\ell > 2$)}
\end{picture}
\end{equation}

\noindent Thus $\gk$ is \setlength{\unitlength}{.5 mm}
\begin{picture}(92,15)
\put(10,10){\circle{2}}
\put(15,0){\circle{2}}
\put(8,5){$\psi_1$}
\put(12,-5){$-\mu$}
\put(25,10){\circle{2}}
\put(23,5){$\psi_3$}
\put(26,10){\line(1,0){13}}
\put(42,10){\circle*{1}}
\put(45,10){\circle*{1}}
\put(48,10){\circle*{1}}
\put(51,10){\line(1,0){13}}
\put(65,10){\circle{2}}
\put(63,5){$\psi_{\ell - 1}$}
\put(66,10.5){\line(1,0){13}}
\put(66,9.5){\line(1,0){13}}
\put(80,10){\circle*{2}}
\put(79,5){$\psi_\ell$}
\end{picture}
and $\gl$ is
\setlength{\unitlength}{.5 mm}
\begin{picture}(100,15)
\put(5,5){\circle{2}}
\put(2,0){$\psi_1$}
\put(6,5){\line(1,0){13}}
\put(17,3){$\times$}
\put(21,5){\line(1,0){13}}
\put(35,5){\circle{2}}
\put(32,0){$\psi_3$}
\put(36,5){\line(1,0){13}}
\put(52,5){\circle*{1}}
\put(55,5){\circle*{1}}
\put(58,5){\circle*{1}}
\put(61,5){\line(1,0){13}}
\put(75,5){\circle{2}}
\put(73,0){$\psi_{\ell - 1}$}
\put(76,5.5){\line(1,0){13}}
\put(76,4.5){\line(1,0){13}}
\put(90,5){\circle*{2}}
\put(89,0){$\psi_\ell$}
\end{picture}.  
\medskip

\noindent
Now the decompositions $\gk_0 = \gk_1 \oplus \gk_2$ and 
$\gl_0 = \gl_1 \oplus \gl_2$ are
\begin{equation}\tag{\theequation$b$}
\gk_0 = \gsp(1) \oplus \left ( \gsp(1) \oplus \gso(2\ell - 3) \right ) 
\text{ and }
\gl_0 = i\R\nu^* \oplus \left ( \gsp(1) \times \gso(2\ell - 3) \right )
\end{equation}
where $\nu^*$ the fundamental simple weight corresponding to $\nu$.
The representation of $\gk$ on $\gs$ has highest weight $-\nu = -\psi_2$
so its diagram is
\begin{picture}(85,15)
\put(10,9){\circle{2}}
\put(10,-1){\circle{2}}
\put(12,7){$1$}
\put(12,-3){$1$}
\put(25,6){\circle{2}}
\put(23,-1){$1$}
\put(26,6){\line(1,0){13}}
\put(42,6){\circle*{1}}
\put(45,6){\circle*{1}}
\put(48,6){\circle*{1}}
\put(51,6){\line(1,0){13}}
\put(65,6){\circle{2}}
\put(66,6.5){\line(1,0){13}}
\put(66,5.5){\line(1,0){13}}
\put(80,6){\circle*{2}}
\end{picture}.
Using (\ref{wtspaces}), the representation 
\begin{equation}\tag{\theequation$c$}
\tau_2:\,\, \gl \text{ on } \gu_2 \text{ is }
\setlength{\unitlength}{.5 mm}
\begin{picture}(100,15)
\put(5,5){\circle{2}}
\put(6,5){\line(1,0){13}}
\put(17,3){$\times$}
\put(17,-2){$1$}
\put(21,5){\line(1,0){13}}
\put(35,5){\circle{2}}
\put(36,5){\line(1,0){13}}
\put(52,5){\circle*{1}}
\put(55,5){\circle*{1}}
\put(58,5){\circle*{1}}
\put(61,5){\line(1,0){13}}
\put(75,5){\circle{2}}
\put(76,5.5){\line(1,0){13}}
\put(76,4.5){\line(1,0){13}}
\put(90,5){\circle*{2}}
\end{picture}
\end{equation}
Also, the action $\tau_{-1} \text{ of } \gl \text{ on } \gu_{-1} \text{ is }
\setlength{\unitlength}{.5 mm}
\begin{picture}(100,15)
\put(5,5){\circle{2}}
\put(3,-2){$1$}
\put(6,5){\line(1,0){13}}
\put(17,3){$\times$}
\put(15,-2){$-2$}
\put(21,5){\line(1,0){13}}
\put(35,5){\circle{2}}
\put(33,-2){$1$}
\put(36,5){\line(1,0){13}}
\put(52,5){\circle*{1}}
\put(55,5){\circle*{1}}
\put(58,5){\circle*{1}}
\put(61,5){\line(1,0){13}}
\put(75,5){\circle{2}}
\put(76,5.5){\line(1,0){13}}
\put(76,4.5){\line(1,0){13}}
\put(90,5){\circle*{2}}
\end{picture}$, 
so the dualizing diagram method of \cite{EW} shows that the representation
\begin{equation}\tag{\theequation$d$}
\tau_1:\,\, \gl \text{ on } \gu_1 \text{ is }
\setlength{\unitlength}{.5 mm}
\begin{picture}(100,10)
\put(5,5){\circle{2}}
\put(3,-2){$1$}
\put(6,5){\line(1,0){13}}
\put(17,3){$\times$}
\put(15,-2){$-1$}
\put(21,5){\line(1,0){13}}
\put(35,5){\circle{2}}
\put(33,-2){$1$}
\put(36,5){\line(1,0){13}}
\put(52,5){\circle*{1}}
\put(55,5){\circle*{1}}
\put(58,5){\circle*{1}}
\put(61,5){\line(1,0){13}}
\put(75,5){\circle{2}}
\put(76,5.5){\line(1,0){13}}
\put(76,4.5){\line(1,0){13}}
\put(90,5){\circle*{2}}
\end{picture}.
\end{equation}
}\end{noname}
\medskip
Here $\dim \gu_2 = 1$ and $\dim \gu_1 = (2\ell - 3)(2\ell - 4)$.

\begin{noname} \label{so-odd-2}{\sc  Case $Spin(2p,2\ell-2p+1), 2<p<\ell$.} 
{\rm Here $G_0$ is the $2$--sheeted cover of $SO(2p,2\ell - 2p + 1)$,
$2 < p < \ell$, contained in $Spin(2\ell + 1;\C)$.
Its extended Dynkin diagram is
\begin{equation}\tag{\theequation$a$}
\setlength{\unitlength}{.75 mm}
\begin{picture}(140,15)
\put(20,10){\circle{2}}
\put(20,0){\circle{2}}
\put(18,5){$\psi_1$}
\put(17,-3){$-\mu$}
\put(21,1){.}
\put(22,2){.}
\put(23,3){.}
\put(24,4){.}
\put(25,5){.}
\put(26,6){.}
\put(27,7){.}
\put(28,8){.}
\put(21,10){\line(1,0){8}}
\put(30,10){\circle{2}}
\put(28,5){$\psi_2$}
\put(31,10){\line(1,0){8}}
\put(42,10){\circle*{1}}
\put(44,10){\circle*{1}}
\put(46,10){\circle*{1}}
\put(49,10){\line(1,0){2}}
\put(52,10){\circle{2}}
\put(53,10){\line(1,0){2}}
\put(50,5){$\psi_p$}
\put(50,12){$\nu$}
\put(59,10){\circle*{1}}
\put(61,10){\circle*{1}}
\put(63,10){\circle*{1}}
\put(66,10){\line(1,0){8}}
\put(75,10){\circle{2}}
\put(73,5){$\psi_{\ell - 1}$}
\put(76,10.5){\line(1,0){13}}
\put(76,9.5){\line(1,0){13}}
\put(90,10){\circle*{2}}
\put(89,5){$\psi_\ell$}
\put(110,10){(type $B_\ell$\,, $\ell > 3$)}
\end{picture}
\end{equation}

\noindent Thus $\gk$ is 
\setlength{\unitlength}{.5 mm}
\begin{picture}(112,15)
\put(10,10){\circle{2}}
\put(10,0){\circle{2}}
\put(08,5){$\psi_1$}
\put(07,-5){$-\mu$}
\put(11,1){\line(1,1){8}}
\put(11,10){\line(1,0){8}}
\put(20,10){\circle{2}}
\put(20,5){$\psi_2$}
\put(21,10){\line(1,0){8}}
\put(32,10){\circle*{1}}
\put(34,10){\circle*{1}}
\put(36,10){\circle*{1}}
\put(39,10){\line(1,0){2}}
\put(42,10){\circle{2}}
\put(37,5){$\psi_{p-1}$}
\put(42,10){\circle{2}}
\put(60,5){$\psi_{p+1}$}
\put(62,10){\circle{2}}
\put(63,10){\line(1,0){2}}
\put(65,10){\circle*{1}}
\put(67,10){\circle*{1}}
\put(69,10){\circle*{1}}
\put(72,10){\line(1,0){8}}
\put(81,10){\circle{2}}
\put(79,5){$\psi_{\ell - 1}$}
\put(82,10.5){\line(1,0){13}}
\put(82,9.5){\line(1,0){13}}
\put(96,10){\circle*{2}}
\put(94,5){$\psi_\ell$}
\end{picture}
and $\gl$ is
\setlength{\unitlength}{.5 mm}
\begin{picture}(120,15)
\put(14,5){\circle{2}}
\put(14,0){$\psi_1$}
\put(15,5){\line(1,0){8}}
\put(24,5){\circle{2}}
\put(22,0){$\psi_2$}
\put(25,5){\line(1,0){8}}
\put(36,5){\circle*{1}}
\put(38,5){\circle*{1}}
\put(40,5){\circle*{1}}
\put(43,5){\line(1,0){8}}
\put(52,5){\circle{2}}
\put(47,0){$\psi_{p-1}$}
\put(53,5){\line(1,0){7}}
\put(59,3){$\times$}
\put(63,5){\line(1,0){7}}
\put(70,0){$\psi_{p+1}$}
\put(72,5){\circle{2}}
\put(73,5){\line(1,0){8}}
\put(81,5){\circle*{1}}
\put(83,5){\circle*{1}}
\put(85,5){\circle*{1}}
\put(88,5){\line(1,0){8}}
\put(97,5){\circle{2}}
\put(95,0){$\psi_{\ell - 1}$}
\put(98,5.5){\line(1,0){13}}
\put(98,4.5){\line(1,0){13}}
\put(112,5){\circle*{2}}
\put(110,0){$\psi_\ell$}
\end{picture}.
\medskip

\noindent
Now the decompositions $\gk_0 = \gk_1 \oplus \gk_2$ and 
$\gl_0 = \gl_1 \oplus \gl_2$ are
\begin{equation}\tag{\theequation$b$}
\gk_0 = \gso(2p) \oplus \gso(2\ell - 2p + 1)\text{ and }
\gl_0 = \gu(p) \oplus \gso(2\ell - 2p + 1)
= i\R\nu^* \oplus \gsu(p) \oplus \gso(2\ell - 2p + 1).
\end{equation}
The representation of $\gk$ on $\gs$ has highest weight $-\nu = -\psi_p$:
\setlength{\unitlength}{.5 mm}
\begin{picture}(112,25)
\put(10,10){\circle{2}}
\put(10,0){\circle{2}}
\put(11,1){\line(1,1){8}}
\put(11,10){\line(1,0){8}}
\put(20,10){\circle{2}}
\put(21,10){\line(1,0){8}}
\put(32,10){\circle*{1}}
\put(34,10){\circle*{1}}
\put(36,10){\circle*{1}}
\put(39,10){\line(1,0){2}}
\put(42,10){\circle{2}}
\put(40,2){$1$}
\put(62,10){\circle{2}}
\put(60,1){$1$}
\put(63,10){\line(1,0){2}}
\put(65,10){\circle*{1}}
\put(67,10){\circle*{1}}
\put(69,10){\circle*{1}}
\put(72,10){\line(1,0){8}}
\put(81,10){\circle{2}}
\put(82,10.5){\line(1,0){13}}
\put(82,9.5){\line(1,0){13}}
\put(96,10){\circle*{2}}
\end{picture}.
\noindent Using (\ref{wtspaces}), the representation 
\begin{equation}\tag{\theequation$c$}
\tau_2:\,\, \gl \text{ on } \gu_2 \text{ is }
\setlength{\unitlength}{.5 mm}
\begin{picture}(120,15)
\put(14,5){\circle{2}}
\put(15,5){\line(1,0){8}}
\put(24,5){\circle{2}}
\put(22,-2){$1$}
\put(25,5){\line(1,0){8}}
\put(36,5){\circle*{1}}
\put(38,5){\circle*{1}}
\put(40,5){\circle*{1}}
\put(43,5){\line(1,0){8}}
\put(52,5){\circle{2}}
\put(53,5){\line(1,0){7}}
\put(59,3){$\times$}
\put(63,5){\line(1,0){7}}
\put(72,5){\circle{2}}
\put(73,5){\line(1,0){8}}
\put(81,5){\circle*{1}}
\put(83,5){\circle*{1}}
\put(85,5){\circle*{1}}
\put(88,5){\line(1,0){8}}
\put(97,5){\circle{2}}
\put(98,5.5){\line(1,0){13}}
\put(98,4.5){\line(1,0){13}}
\put(112,5){\circle*{2}}
\end{picture}.
\end{equation}
Also, the action $\tau_{-1}$ of $\gl$ on $\gu_{-1}$ is 
\setlength{\unitlength}{.5 mm}
\begin{picture}(120,15)
\put(14,5){\circle{2}}
\put(15,5){\line(1,0){8}}
\put(24,5){\circle{2}}
\put(25,5){\line(1,0){8}}
\put(36,5){\circle*{1}}
\put(38,5){\circle*{1}}
\put(40,5){\circle*{1}}
\put(43,5){\line(1,0){8}}
\put(52,5){\circle{2}}
\put(50,-2){$1$}
\put(53,5){\line(1,0){7}}
\put(59,3){$\times$}
\put(57,-2){$-2$}
\put(63,5){\line(1,0){7}}
\put(72,5){\circle{2}}
\put(70,-2){$1$}
\put(73,5){\line(1,0){8}}
\put(81,5){\circle*{1}}
\put(83,5){\circle*{1}}
\put(85,5){\circle*{1}}
\put(88,5){\line(1,0){8}}
\put(97,5){\circle{2}}
\put(98,5.5){\line(1,0){13}}
\put(98,4.5){\line(1,0){13}}
\put(112,5){\circle*{2}}
\end{picture}
so the dualizing diagram method of \cite{EW} shows that the representation
\begin{equation}\tag{\theequation$d$}
\tau_1:\,\, \gl \text{ on } \gu_1 \text{ is }
\setlength{\unitlength}{.5 mm}
\begin{picture}(120,15)
\put(14,5){\circle{2}}
\put(12,-2){$1$}
\put(15,5){\line(1,0){8}}
\put(24,5){\circle{2}}
\put(25,5){\line(1,0){8}}
\put(36,5){\circle*{1}}
\put(38,5){\circle*{1}}
\put(40,5){\circle*{1}}
\put(43,5){\line(1,0){8}}
\put(52,5){\circle{2}}
\put(53,5){\line(1,0){7}}
\put(59,3){$\times$}
\put(57,-2){$-1$}
\put(63,5){\line(1,0){7}}
\put(72,5){\circle{2}}
\put(70,-2){$1$}
\put(73,5){\line(1,0){8}}
\put(81,5){\circle*{1}}
\put(83,5){\circle*{1}}
\put(85,5){\circle*{1}}
\put(88,5){\line(1,0){8}}
\put(97,5){\circle{2}}
\put(98,5.5){\line(1,0){13}}
\put(98,4.5){\line(1,0){13}}
\put(112,5){\circle*{2}}
\end{picture}.
\end{equation}
}\end{noname}
\medskip
Here $\dim \gu_2 = p(p-1)/2$ and $\dim \gu_1 = p(2\ell - 2p +1)$.

\begin{noname} \label{so-odd-3}{\sc  Case $Spin(4,1)$.} {\rm
Here $G_0$ is the (universal) double cover of the group $SO(4,1)$.
Its extended Dynkin diagram is
\begin{equation}\tag{\theequation$a$}
\setlength{\unitlength}{.75 mm}
\begin{picture}(90,10)
\put(15,5){\circle{2}}
\put(13,0){$\psi_{1}$}
\put(16,5.5){\line(1,0){13}}
\put(16,4.5){\line(1,0){13}}
\put(30,5){\circle*{2}}
\put(29,0){$\psi_2$}
\put(29,8){$\nu$}
\put(31,5.5){.}
\put(32,5.5){.}
\put(33,5.5){.}
\put(34,5.5){.}
\put(35,5.5){.}
\put(36,5.5){.}
\put(37,5.5){.}
\put(38,5.5){.}
\put(39,5.5){.}
\put(40,5.5){.}
\put(41,5.5){.}
\put(42,5.5){.}
\put(43,5.5){.}
\put(31,4.5){.}
\put(32,4.5){.}
\put(33,4.5){.}
\put(34,4.5){.}
\put(35,4.5){.}
\put(36,4.5){.}
\put(37,4.5){.}
\put(38,4.5){.}
\put(39,4.5){.}
\put(40,4.5){.}
\put(41,4.5){.}
\put(42,4.5){.}
\put(43,4.5){.}
\put(45,5){\circle{2}}
\put(42,0){$-\mu$}
\put(70,5){(type $B_2$)}
\end{picture}
\end{equation}

\noindent Thus $\gk$ is 
\setlength{\unitlength}{.5 mm}
\begin{picture}(28,10)
\put(05,5){\circle{2}}
\put(03,0){$\psi_{1}$}
\put(20,5){\circle{2}}
\put(15,0){$-\mu$}
\end{picture}
and $\gl$ is
\setlength{\unitlength}{.5 mm}
\begin{picture}(28,10)
\put(5,5){\circle{2}}
\put(3,0){$\psi_{1}$}
\put(6,5.5){\line(1,0){13}}
\put(6,4.5){\line(1,0){13}}
\put(18,3.1){\tiny $\boxtimes$}
\end{picture}.
Now the decompositions $\gk_0 = \gk_1 \oplus \gk_2$ and 
$\gl_0 = \gl_1 \oplus \gl_2$ are
\begin{equation}\tag{\theequation$b$}
\gk_0 = \gsp(1) \oplus \gsp(1) \text{ and } \gl_0 = i\R\nu^* \oplus \gsp(1).
\end{equation}
Here $\gk$ acts on $\gs$ with highest weight $-\nu = -\psi_2$:
\setlength{\unitlength}{.5 mm}
\begin{picture}(22,10)
\put(05,7){\circle{2}}
\put(03,0){$1$}
\put(16,7){\circle{2}}
\put(15,0){$1$}
\end{picture}.
\noindent Using (\ref{wtspaces}), the representation 
\begin{equation}\tag{\theequation$c$}
\tau_2:\,\, \gl \text{ on } \gu_2 \text{ is }
\setlength{\unitlength}{.5 mm}
\begin{picture}(28,10)
\put(5,8){\circle{2}}
\put(6,8.5){\line(1,0){13}}
\put(6,7.5){\line(1,0){13}}
\put(18,6.1){\tiny $\boxtimes$}
\put(16,0){$-2$}
\end{picture}.
\end{equation}
Also, the action $\tau_{-1} \text{ of } \gl \text{ on } \gu_{-1} \text{ is }$
\setlength{\unitlength}{.5 mm}
\begin{picture}(28,10)
\put(5,8){\circle{2}}
\put(5,0){$1$}
\put(6,8.5){\line(1,0){13}}
\put(6,7.5){\line(1,0){13}}
\put(18,6.1){\tiny $\boxtimes$}
\put(16,0){$-2$}
\end{picture}
so the dualizing diagram method of \cite{EW} shows that the representation
\begin{equation}\tag{\theequation $d$}
\tau_1:\,\, \gl \text{ on } \gu_1 \text{ is }
\setlength{\unitlength}{.5 mm}
\begin{picture}(28,10)
\put(5,8){\circle{2}}
\put(5,0){$1$}
\put(6,8.5){\line(1,0){13}}
\put(6,7.5){\line(1,0){13}}
\put(18,6.1){\tiny $\boxtimes$}
\put(16,0){$-1$}
\end{picture}.
\end{equation}
}\end{noname}
\medskip
$\dim \gu_2 = 1$ and $\dim \gu_1 = 2$.

\begin{noname}\label{so-odd-4} {\sc Case $Spin(2\ell, 1),\; \ell > 1$.}  {\rm
Here $G_0$ is the universal ($2$--sheeted) cover 
of the group $SO(2\ell,1)$ with $\ell > 1$.  Its extended Dynkin diagram is
\begin{equation}\tag{\theequation$a$}
\setlength{\unitlength}{.75 mm}
\begin{picture}(140,15)
\put(20,10){\circle{2}}
\put(20,0){\circle{2}}
\put(18,5){$\psi_1$}
\put(17,-3){$-\mu$}
\put(21,1){.}
\put(22,2){.}
\put(23,3){.}
\put(24,4){.}
\put(25,5){.}
\put(26,6){.}
\put(27,7){.}
\put(28,8){.}
\put(21,10){\line(1,0){8}}
\put(30,10){\circle{2}}
\put(28,5){$\psi_2$}
\put(31,10){\line(1,0){8}}
\put(42,10){\circle*{1}}
\put(45,10){\circle*{1}}
\put(48,10){\circle*{1}}
\put(51,10){\line(1,0){8}}
\put(60,10){\circle{2}}
\put(58,5){$\psi_{\ell - 1}$}
\put(61,10.5){\line(1,0){13}}
\put(61,9.5){\line(1,0){13}}
\put(75,10){\circle*{2}}
\put(74,5){$\psi_\ell$}
\put(74,12){$\nu$}
\put(95,10){(type $B_\ell$\,, $\ell > 2$)}
\end{picture}
\end{equation}

\noindent Thus $\gk$ is 
\setlength{\unitlength}{.6 mm}
\begin{picture}(65,15)
\put(10,10){\circle{2}}
\put(10,0){\circle{2}}
\put(08,5){$\psi_1$}
\put(09,-3){$-\mu$}
\put(11,1){\line(1,1){8}}
\put(11,10){\line(1,0){8}}
\put(20,10){\circle{2}}
\put(18,5){$\psi_2$}
\put(21,10){\line(1,0){8}}
\put(32,10){\circle*{1}}
\put(35,10){\circle*{1}}
\put(38,10){\circle*{1}}
\put(41,10){\line(1,0){8}}
\put(50,10){\circle{2}}
\put(48,5){$\psi_{\ell - 1}$}
\end{picture}
and $\gl$ is
\setlength{\unitlength}{.6 mm}
\begin{picture}(70,15)
\put(10,7){\circle{2}}
\put(08,2){$\psi_1$}
\put(11,7){\line(1,0){8}}
\put(20,7){\circle{2}}
\put(18,2){$\psi_2$}
\put(21,7){\line(1,0){8}}
\put(32,7){\circle*{1}}
\put(35,7){\circle*{1}}
\put(38,7){\circle*{1}}
\put(41,7){\line(1,0){8}}
\put(50,7){\circle{2}}
\put(48,2){$\psi_{\ell - 1}$}
\put(51,7.5){\line(1,0){13}}
\put(51,6.5){\line(1,0){13}}
\put(63,5.6){\tiny $\boxtimes$}
\end{picture}.
\medskip

\noindent
Now the decompositions $\gk_0 = \gk_1 \oplus \gk_2$ and 
$\gl_0 = \gl_1 \oplus \gl_2$ are
\begin{equation}\tag{\theequation$b$}
\gk_0 = \gso(2\ell) \text{ and } \gl_0 = \gu(\ell) = i\R\nu^* \oplus \gsu(\ell).
\end{equation}
The representation of $\gk$ on $\gs$ has highest weight $-\nu = -\psi_\ell$:
\setlength{\unitlength}{.5 mm}
\begin{picture}(70,15)
\put(10,10){\circle{2}}
\put(10,0){\circle{2}}
\put(11,1){\line(1,1){8}}
\put(11,10){\line(1,0){8}}
\put(20,10){\circle{2}}
\put(21,10){\line(1,0){8}}
\put(32,10){\circle*{1}}
\put(35,10){\circle*{1}}
\put(38,10){\circle*{1}}
\put(41,10){\line(1,0){8}}
\put(50,10){\circle{2}}
\put(48,2){$1$}
\end{picture}.
\noindent Using (\ref{wtspaces}), the representation 
\begin{equation}\tag{\theequation$c$}
\tau_2:\,\, \gl \text{ on } \gu_2 \text{ is }
\setlength{\unitlength}{.5 mm}
\begin{picture}(70,15)
\put(10,7){\circle{2}}
\put(11,7){\line(1,0){8}}
\put(20,7){\circle{2}}
\put(19,0){$1$}
\put(21,7){\line(1,0){8}}
\put(32,7){\circle*{1}}
\put(35,7){\circle*{1}}
\put(38,7){\circle*{1}}
\put(41,7){\line(1,0){8}}
\put(50,7){\circle{2}}
\put(51,7.5){\line(1,0){13}}
\put(51,6.5){\line(1,0){13}}
\put(63,5.6){\tiny $\boxtimes$}
\end{picture}.
\end{equation}
Also, the action $\tau_{-1} \text{ of } \gl \text{ on } \gu_{-1} \text{ is }
\setlength{\unitlength}{.5 mm}
\begin{picture}(70,15)
\put(10,7){\circle{2}}
\put(11,7){\line(1,0){8}}
\put(20,7){\circle{2}}
\put(21,7){\line(1,0){8}}
\put(32,7){\circle*{1}}
\put(35,7){\circle*{1}}
\put(38,7){\circle*{1}}
\put(41,7){\line(1,0){8}}
\put(50,7){\circle{2}}
\put(49,0){$1$}
\put(51,7.5){\line(1,0){13}}
\put(51,6.5){\line(1,0){13}}
\put(63,5.6){\tiny $\boxtimes$}
\put(60,0){$-2$}
\end{picture}$
so the dualizing diagram method of \cite{EW} shows that the representation
\begin{equation}\tag{\theequation$d$}
\tau_1:\,\, \gl \text{ on } \gu_1 \text{ is }
\setlength{\unitlength}{.5 mm}
\begin{picture}(70,15)
\put(10,7){\circle{2}}
\put(11,7){\line(1,0){8}}
\put(20,7){\circle{2}}
\put(9,0){$1$}
\put(21,7){\line(1,0){8}}
\put(32,7){\circle*{1}}
\put(35,7){\circle*{1}}
\put(38,7){\circle*{1}}
\put(41,7){\line(1,0){8}}
\put(50,7){\circle{2}}
\put(51,7.5){\line(1,0){13}}
\put(51,6.5){\line(1,0){13}}
\put(63,5.6){\tiny $\boxtimes$}
\end{picture}.
\end{equation}
}\end{noname}

\noindent
Here  $\dim \gu_2 = \ell(\ell - 1)/2$ and $\dim \gu_1 = \ell$.
This exhausts the cases where $\gg$ is of type $B$, and we go on to
consider the cases where $\gg$ is of type $C$.
\medskip

\begin{noname}\label{sp-1} {\sc Case $Sp(p,\ell - p), 1 < p < \ell$.} {\rm  
Here $G_0$ is simply connected, and its extended Dynkin diagram is
\begin{equation}\tag{\theequation$a$}
\setlength{\unitlength}{.75 mm}
\begin{picture}(160,10)
\put(5,5){\circle{2}}
\put(2,0){$-\mu$}
\put(6,5.5){.}
\put(7,5.5){.}
\put(8,5.5){.}
\put(9,5.5){.}
\put(10,5.5){.}
\put(11,5.5){.}
\put(12,5.5){.}
\put(13,5.5){.}
\put(6,4.5){.}
\put(7,4.5){.}
\put(8,4.5){.}
\put(9,4.5){.}
\put(10,4.5){.}
\put(11,4.5){.}
\put(12,4.5){.}
\put(13,4.5){.}
\put(15,5){\circle*{2}}
\put(13,0){$\psi_1$}
\put(16,5){\line(1,0){8}}
\put(25,5){\circle*{2}}
\put(23,0){$\psi_2$}
\put(26,5){\line(1,0){8}}
\put(37,5){\circle*{1}}
\put(40,5){\circle*{1}}
\put(43,5){\circle*{1}}
\put(46,5){\line(1,0){8}}
\put(55,5){\circle*{2}}
\put(56,5){\line(1,0){8}}
\put(54,7){$\nu$}
\put(53,0){$\psi_p$}
\put(67,5){\circle*{1}}
\put(70,5){\circle*{1}}
\put(73,5){\circle*{1}}
\put(76,5){\line(1,0){8}}
\put(85,5){\circle*{2}}
\put(83,0){$\psi_{\ell - 1}$}
\put(86,5.5){\line(1,0){8}}
\put(86,4.5){\line(1,0){8}}
\put(95,5){\circle{2}}
\put(95,0){$\psi_\ell$}
\put(110,3){(type $C_\ell$\,, $\ell > 1$)}
\end{picture}
\end{equation}
Thus $\gk$ is
{\footnotesize
\setlength{\unitlength}{.5 mm}
\begin{picture}(125,10)
\put(5,5){\circle{2}}
\put(0,0){$-\mu$}
\put(6,5.5){\line(1,0){8}}
\put(6,4.5){\line(1,0){8}}
\put(15,5){\circle{2}}
\put(13,0){$\psi_1$}
\put(16,5){\line(1,0){8}}
\put(25,5){\circle*{2}}
\put(23,0){$\psi_2$}
\put(26,5){\line(1,0){8}}
\put(37,5){\circle*{1}}
\put(40,5){\circle*{1}}
\put(43,5){\circle*{1}}
\put(46,5){\line(1,0){8}}
\put(55,5){\circle*{2}}
\put(53,0){$\psi_{p-1}$}
\put(75,5){\circle*{2}}
\put(76,5){\line(1,0){8}}
\put(73,0){$\psi_{p+1}$}
\put(87,5){\circle*{1}}
\put(90,5){\circle*{1}}
\put(93,5){\circle*{1}}
\put(96,5){\line(1,0){8}}
\put(105,5){\circle*{2}}
\put(103,0){$\psi_{\ell - 1}$}
\put(106,5.5){\line(1,0){8}}
\put(106,4.5){\line(1,0){8}}
\put(115,5){\circle{2}}
\put(115,0){$\psi_\ell$}
\end{picture}}
with $\Psi_{\gk_1} = \{-\mu, \psi_1, \psi_2, \dots , \psi_{p-1}\}$ and
$\Psi_{\gk_2} = \{\psi_{p+1}, \psi_{p+2}, \dots , \psi_\ell\}$,
and $\gl$ is
{\footnotesize
\setlength{\unitlength}{.5 mm}
\begin{picture}(100,10)
\put(5,5){\circle*{2}}
\put(3,0){$\psi_1$}
\put(6,5){\line(1,0){8}}
\put(15,5){\circle*{2}}
\put(13,0){$\psi_2$}
\put(16,5){\line(1,0){8}}
\put(27,5){\circle*{1}}
\put(30,5){\circle*{1}}
\put(33,5){\circle*{1}}
\put(36,5){\line(1,0){8}}
\put(43,3.5){\tiny $\boxtimes$}
\put(46,5){\line(1,0){8}}
\put(57,5){\circle*{1}}
\put(60,5){\circle*{1}}
\put(63,5){\circle*{1}}
\put(66,5){\line(1,0){8}}
\put(75,5){\circle*{2}}
\put(73,0){$\psi_{\ell - 1}$}
\put(76,5.5){\line(1,0){8}}
\put(76,4.5){\line(1,0){8}}
\put(85,5){\circle{2}}
\put(88,0){$\psi_\ell$}
\end{picture}}.
\medskip

\noindent
Now the decompositions $\gk_0 = \gk_1 \oplus \gk_2$ and
$\gl_0 = \gl_1 \oplus \gl_2$ are
\begin{equation}\tag{\theequation$b$}
\gk_0 = \gsp(p) \oplus \gsp(\ell - p) \text{ and } 
\gl_0 = \gu(p) \oplus \gsp(\ell - p) = 
i\R\nu^* \oplus \gsu(p) \oplus \gsp(\ell - p).
\end{equation}
The representation of $\gk$ on $\gs$ has highest weight $-\nu = -\psi_p$:
{\footnotesize
\setlength{\unitlength}{.5 mm}
\begin{picture}(115,10)
\put(0,5){\circle{2}}
\put(1,5.5){\line(1,0){8}}
\put(1,4.5){\line(1,0){8}}
\put(10,5){\circle*{2}}
\put(11,5){\line(1,0){8}}
\put(20,5){\circle*{2}}
\put(21,5){\line(1,0){8}}
\put(32,5){\circle*{1}}
\put(35,5){\circle*{1}}
\put(38,5){\circle*{1}}
\put(41,5){\line(1,0){8}}
\put(50,5){\circle*{2}}
\put(48,-2){$1$}
\put(68,-2){$1$}
\put(70,5){\circle*{2}}
\put(71,5){\line(1,0){8}}
\put(82,5){\circle*{1}}
\put(85,5){\circle*{1}}
\put(88,5){\circle*{1}}
\put(91,5){\line(1,0){8}}
\put(100,5){\circle*{2}}
\put(101,5.5){\line(1,0){8}}
\put(101,4.5){\line(1,0){8}}
\put(110,5){\circle{2}}
\end{picture}}.
\noindent Using (\ref{wtspaces}), the representation
\begin{equation}\tag{\theequation$c$}
\tau_2:\,\, \gl \text{ on } \gu_2 \text{ is }
\setlength{\unitlength}{.5 mm}
\begin{picture}(110,10)
\put(5,5){\circle*{2}}
\put(3,-2){$2$}
\put(6,5){\line(1,0){8}}
\put(15,5){\circle*{2}}
\put(16,5){\line(1,0){8}}
\put(27,5){\circle*{1}}
\put(30,5){\circle*{1}}
\put(33,5){\circle*{1}}
\put(36,5){\line(1,0){8}}
\put(45,5){\circle*{2}}
\put(46,5){\line(1,0){8}}
\put(54,3){\tiny $\boxtimes$}
\put(57,5){\line(1,0){8}}
\put(66,5){\circle*{2}}
\put(67,5){\line(1,0){8}}
\put(78,5){\circle*{1}}
\put(81,5){\circle*{1}}
\put(84,5){\circle*{1}}
\put(87,5){\line(1,0){8}}
\put(96,5){\circle*{2}}
\put(97,5.5){\line(1,0){8}}
\put(97,4.5){\line(1,0){8}}
\put(106,5){\circle{2}}
\end{picture}
\end{equation}
The action $\tau_{-1}$ of $\gl$ on $\gu_{-1}$ is
\setlength{\unitlength}{.5 mm}
\begin{picture}(110,10)
\put(5,5){\circle*{2}}
\put(6,5){\line(1,0){8}}
\put(15,5){\circle*{2}}
\put(16,5){\line(1,0){8}}
\put(27,5){\circle*{1}}
\put(30,5){\circle*{1}}
\put(33,5){\circle*{1}}
\put(36,5){\line(1,0){8}}
\put(45,5){\circle*{2}}
\put(43,-2){$1$}
\put(46,5){\line(1,0){8}}
\put(53.5,3){\tiny $\boxtimes$}
\put(50,-2){$-2$}
\put(57,5){\line(1,0){8}}
\put(66,5){\circle*{2}}
\put(64,-2){$1$}
\put(67,5){\line(1,0){8}}
\put(78,5){\circle*{1}}
\put(81,5){\circle*{1}}
\put(84,5){\circle*{1}}
\put(87,5){\line(1,0){8}}
\put(96,5){\circle*{2}}
\put(97,5.5){\line(1,0){8}}
\put(97,4.5){\line(1,0){8}}
\put(106,5){\circle{2}}
\end{picture}, 
so the dualizing diagram method of \cite{EW} shows that the representation
\begin{equation}\tag{\theequation$d$}
\tau_1:\,\, \gl \text{ on } \gu_1 \text{ is }
\setlength{\unitlength}{.5 mm}
\begin{picture}(110,10)
\put(5,5){\circle*{2}}
\put(3,-2){$1$}
\put(6,5){\line(1,0){8}}
\put(15,5){\circle*{2}}
\put(16,5){\line(1,0){8}}
\put(27,5){\circle*{1}}
\put(30,5){\circle*{1}}
\put(33,5){\circle*{1}}
\put(36,5){\line(1,0){8}}
\put(45,5){\circle*{2}}
\put(46,5){\line(1,0){8}}
\put(53.5,3){\tiny $\boxtimes$}
\put(50,-2){$-1$}
\put(57,5){\line(1,0){8}}
\put(66,5){\circle*{2}}
\put(64,-2){$1$}
\put(67,5){\line(1,0){8}}
\put(78,5){\circle*{1}}
\put(81,5){\circle*{1}}
\put(84,5){\circle*{1}}
\put(87,5){\line(1,0){8}}
\put(96,5){\circle*{2}}
\put(97,5.5){\line(1,0){8}}
\put(97,4.5){\line(1,0){8}}
\put(106,5){\circle{2}}
\end{picture}
\end{equation}
}\end{noname}
Here $\dim \gu_2 = (p-1)(p+2)/2$ and $\dim \gu_1 = 2p(\ell - p)$.
\medskip

\begin{noname}\label{sp-2} {\sc Case $Sp(1,\ell - 1)$}.  
{\rm Here $G_0$ is simply connected, and its extended Dynkin diagram is
\begin{equation}\tag{\theequation$a$}
\setlength{\unitlength}{.75 mm}
\begin{picture}(160,10)
\put(5,5){\circle{2}}
\put(2,0){$-\mu$}
\put(6,5.5){.}
\put(7,5.5){.}
\put(8,5.5){.}
\put(9,5.5){.}
\put(10,5.5){.}
\put(11,5.5){.}
\put(12,5.5){.}
\put(13,5.5){.}
\put(6,4.5){.}
\put(7,4.5){.}
\put(8,4.5){.}
\put(9,4.5){.}
\put(10,4.5){.}
\put(11,4.5){.}
\put(12,4.5){.}
\put(13,4.5){.}
\put(15,5){\circle*{2}}
\put(14,7){$\nu$}
\put(13,0){$\psi_1$}
\put(16,5){\line(1,0){8}}
\put(25,5){\circle*{2}}
\put(23,0){$\psi_2$}
\put(26,5){\line(1,0){8}}
\put(37,5){\circle*{1}}
\put(40,5){\circle*{1}}
\put(43,5){\circle*{1}}
\put(46,5){\line(1,0){8}}
\put(55,5){\circle*{2}}
\put(56,5){\line(1,0){8}}
\put(53,0){$\psi_p$}
\put(67,5){\circle*{1}}
\put(70,5){\circle*{1}}
\put(73,5){\circle*{1}}
\put(76,5){\line(1,0){8}}
\put(85,5){\circle*{2}}
\put(83,0){$\psi_{\ell - 1}$}
\put(86,5.5){\line(1,0){8}}
\put(86,4.5){\line(1,0){8}}
\put(95,5){\circle{2}}
\put(95,0){$\psi_\ell$}
\put(110,3){(type $C_\ell$\,, $\ell > 1$)}
\end{picture}
\end{equation}
Thus $\gk$ is
\setlength{\unitlength}{.5 mm}
\begin{picture}(90,10)
\put(5,5){\circle{2}}
\put(0,0){$-\mu$}
\put(25,5){\circle*{2}}
\put(26,5){\line(1,0){13}}
\put(23,0){$\psi_2$}
\put(42,5){\circle*{1}}
\put(45,5){\circle*{1}}
\put(48,5){\circle*{1}}
\put(51,5){\line(1,0){13}}
\put(65,5){\circle*{2}}
\put(63,0){$\psi_{\ell - 1}$}
\put(66,5.5){\line(1,0){13}}
\put(66,4.5){\line(1,0){13}}
\put(80,5){\circle{2}}
\put(80,0){$\psi_\ell$}
\end{picture}
and $\gl$ is
\setlength{\unitlength}{.5 mm}
\begin{picture}(100,10)
\put(2.5,3.5){\tiny $\boxtimes$}
\put(6,5){\line(1,0){13}}
\put(20,5){\circle*{2}}
\put(18,0){$\psi_2$}
\put(21,5){\line(1,0){13}}
\put(37,5){\circle*{1}}
\put(40,5){\circle*{1}}
\put(43,5){\circle*{1}}
\put(46,5){\line(1,0){13}}
\put(60,5){\circle*{2}}
\put(58,0){$\psi_{\ell - 1}$}
\put(61,5.5){\line(1,0){13}}
\put(61,4.5){\line(1,0){13}}
\put(75,5){\circle{2}}
\put(78,0){$\psi_\ell$}
\end{picture}.
\medskip

\noindent
Now the decompositions $\gk_0 = \gk_1 \oplus \gk_2$ and
$\gl_0 = \gl_1 \oplus \gl_2$ are
\begin{equation}\tag{\theequation$b$}
\gk_0 = \gsp(1)\oplus \gsp(\ell - 1) \text{ and } 
\gl_0 = i\R\nu^* \oplus \gsp(\ell - 1).
\end{equation}
The representation of $\gk$ on $\gs$ has highest weight $-\nu = -\psi_1$:
{\footnotesize
\setlength{\unitlength}{.5 mm}
\begin{picture}(90,10)
\put(5,5){\circle{2}}
\put(3,-2){$1$}
\put(23,-2){$1$}
\put(25,5){\circle*{2}}
\put(26,5){\line(1,0){13}}
\put(42,5){\circle*{1}}
\put(45,5){\circle*{1}}
\put(48,5){\circle*{1}}
\put(51,5){\line(1,0){13}}
\put(65,5){\circle*{2}}
\put(66,5.5){\line(1,0){13}}
\put(66,4.5){\line(1,0){13}}
\put(80,5){\circle{2}}
\end{picture}}.
\noindent Using (\ref{wtspaces}), the representation
\begin{equation}\tag{\theequation$c$}
\tau_2:\,\, \gl \text{ on } \gu_2 \text{ is }
\setlength{\unitlength}{.5 mm}
\begin{picture}(100,10)
\put(2.5,3){\tiny $\boxtimes$}
\put(3,-2.5){$2$}
\put(6,5){\line(1,0){13}}
\put(20,5){\circle*{2}}
\put(21,5){\line(1,0){13}}
\put(37,5){\circle*{1}}
\put(40,5){\circle*{1}}
\put(43,5){\circle*{1}}
\put(46,5){\line(1,0){13}}
\put(60,5){\circle*{2}}
\put(61,5.5){\line(1,0){13}}
\put(61,4.5){\line(1,0){13}}
\put(75,5){\circle{2}}
\end{picture}
\end{equation}
The action $\tau_{-1}$ of $\gl$ on $\gu_{-1}$ is
\setlength{\unitlength}{.5 mm}
\begin{picture}(80,10)
\put(2.5,3){\tiny $\boxtimes$}
\put(1,-2.5){$-2$}
\put(6,5){\line(1,0){13}}
\put(20,5){\circle*{2}}
\put(18,-2.5){$1$}
\put(21,5){\line(1,0){13}}
\put(37,5){\circle*{1}}
\put(40,5){\circle*{1}}
\put(43,5){\circle*{1}}
\put(46,5){\line(1,0){13}}
\put(60,5){\circle*{2}}
\put(61,5.5){\line(1,0){13}}
\put(61,4.5){\line(1,0){13}}
\put(75,5){\circle{2}}
\end{picture}, 
so the dualizing diagram method of \cite{EW} shows that the representation
\begin{equation}\tag{\theequation$d$}
\tau_1:\,\, \gl \text{ on } \gu_1 \text{ is }
\setlength{\unitlength}{.5 mm}
\begin{picture}(80,10)
\put(2.5,3){\tiny $\boxtimes$}
\put(6,5){\line(1,0){13}}
\put(20,5){\circle*{2}}
\put(18,-2.5){$1$}
\put(21,5){\line(1,0){13}}
\put(37,5){\circle*{1}}
\put(40,5){\circle*{1}}
\put(43,5){\circle*{1}}
\put(46,5){\line(1,0){13}}
\put(60,5){\circle*{2}}
\put(61,5.5){\line(1,0){13}}
\put(61,4.5){\line(1,0){13}}
\put(75,5){\circle{2}}
\end{picture}.
\end{equation}
\medskip
}\end{noname}
\noindent
Here $\dim \gu_2 = 1$ and $\dim \gu_1 = 2\ell$.
This exhausts the cases where $\gg$ is of type $C$, and we go on to
consider the cases where $\gg$ is of type $D$.
\medskip

\begin{noname}\label{so-even-1} {\sc Case $Spin(4,2\ell - 4), \ell > 4$}.  
{\rm Here $G_0$ is the $2$--sheeted cover of $SO(4,2\ell - 4)$
contained in $Spin(2\ell;\C)$.  Its extended Dynkin diagram is
\begin{equation}\tag{\theequation$a$}
\setlength{\unitlength}{.75 mm}
\begin{picture}(140,15)
\put(20,10){\circle{2}}
\put(25,0){\circle{2}}
\put(18,5){$\psi_1$}
\put(22,-3){$-\mu$}
\put(26,1){.}
\put(27,2){.}
\put(28,3){.}
\put(29,4){.}
\put(30,5){.}
\put(31,6){.}
\put(32,7){.}
\put(33,8){.}
\put(21,10){\line(1,0){13}}
\put(35,10){\circle{2}}
\put(33,5){$\psi_2$}
\put(33,12){$\nu$}
\put(36,10){\line(1,0){13}}
\put(52,10){\circle*{1}}
\put(55,10){\circle*{1}}
\put(58,10){\circle*{1}}
\put(61,10){\line(1,0){13}}
\put(75,10){\circle{2}}
\put(70,5){$\psi_{\ell - 2}$}
\put(76,9.5){\line(2,-1){13}}
\put(90,3){\circle{2}}
\put(92,2){$\psi_\ell$}
\put(76,10.5){\line(2,1){13}}
\put(90,17){\circle{2}}
\put(92,16){$\psi_{\ell - 1}$}
\put(110,10){(type $D_\ell$\,, $\ell > 4$)}
\end{picture}
\end{equation}

\noindent Thus $\gk$ is \setlength{\unitlength}{.5 mm}
\begin{picture}(92,15)
\put(10,10){\circle{2}}
\put(15,0){\circle{2}}
\put(8,4){$\psi_1$}
\put(12,-5){$-\mu$}
\put(25,10){\circle{2}}
\put(23,4){$\psi_3$}
\put(26,10){\line(1,0){13}}
\put(42,10){\circle*{1}}
\put(45,10){\circle*{1}}
\put(48,10){\circle*{1}}
\put(51,10){\line(1,0){13}}
\put(65,10){\circle{2}}
\put(58,4){$\psi_{\ell - 2}$}
\put(66,9.5){\line(2,-1){13}}
\put(80,3){\circle{2}}
\put(82,2){$\psi_\ell$}
\put(66,10.5){\line(2,1){13}}
\put(80,17){\circle{2}}
\put(82,16){$\psi_{\ell - 1}$}
\end{picture}
and $\gl$ is
\setlength{\unitlength}{.5 mm}
\begin{picture}(100,15)
\put(5,5){\circle{2}}
\put(2,0){$\psi_1$}
\put(6,5){\line(1,0){13}}
\put(17,3){$\times$}
\put(21,5){\line(1,0){13}}
\put(35,5){\circle{2}}
\put(32,0){$\psi_3$}
\put(36,5){\line(1,0){13}}
\put(52,5){\circle*{1}}
\put(55,5){\circle*{1}}
\put(58,5){\circle*{1}}
\put(61,5){\line(1,0){13}}
\put(75,5){\circle{2}}
\put(67,0){$\psi_{\ell - 2}$}
\put(76,4.5){\line(2,-1){13}}
\put(90,-2){\circle{2}}
\put(92,-3){$\psi_\ell$}
\put(76,5.5){\line(2,1){13}}
\put(90,12){\circle{2}}
\put(92,11){$\psi_{\ell - 1}$}
\end{picture}.  
\medskip

\noindent
Now the decompositions $\gk_0 = \gk_1 \oplus \gk_2$ and
$\gl_0 = \gl_1 \oplus \gl_2$ are
\begin{equation}\tag{\theequation$b$}
\gk_0 = \gsp(1) \oplus \left ( \gsp(1) \oplus \gso(2\ell -4)\right ) \text{ and }
\gl_0 = i\R\nu^* \oplus \left ( \gsp(1)\oplus \gso(2\ell -4)\right ).
\end{equation}
The representation of $\gk$ on $\gs$ has highest weight $-\nu = -\psi_2$
so its diagram is
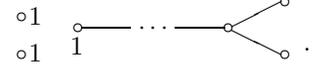
\begin{picture}(85,15)
\put(10,9){\circle{2}}
\put(10,-1){\circle{2}}
\put(12,7){$1$}
\put(12,-3){$1$}
\put(25,6){\circle{2}}
\put(23,-1){$1$}
\put(26,6){\line(1,0){13}}
\put(42,6){\circle*{1}}
\put(45,6){\circle*{1}}
\put(48,6){\circle*{1}}
\put(51,6){\line(1,0){13}}
\put(65,6){\circle{2}}
\put(66,5.5){\line(2,-1){13}}
\put(80,-1){\circle{2}}
\put(66,6.5){\line(2,1){13}}
\put(80,13){\circle{2}}

\end{picture}.
Using (\ref{wtspaces}), the representation 
\begin{equation}\tag{\theequation$c$}
\tau_2:\,\, \gl \text{ on } \gu_2 \text{ is }
\setlength{\unitlength}{.5 mm}
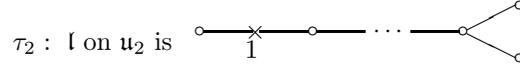
\begin{picture}(100,15)
\put(5,5){\circle{2}}
\put(6,5){\line(1,0){13}}
\put(17,3){$\times$}
\put(17,-2){$1$}
\put(21,5){\line(1,0){13}}
\put(35,5){\circle{2}}
\put(36,5){\line(1,0){13}}
\put(52,5){\circle*{1}}
\put(55,5){\circle*{1}}
\put(58,5){\circle*{1}}
\put(61,5){\line(1,0){13}}
\put(75,5){\circle{2}}
\put(76,4.5){\line(2,-1){13}}
\put(90,-2){\circle{2}}
\put(76,5.5){\line(2,1){13}}
\put(90,12){\circle{2}}
\end{picture}
\end{equation}
Also, the action $\tau_{-1} \text{ of } \gl \text{ on } \gu_{-1} \text{ is }
\setlength{\unitlength}{.5 mm}
\begin{picture}(100,15)
\put(5,5){\circle{2}}
\put(3,-2){$1$}
\put(6,5){\line(1,0){13}}
\put(17,3){$\times$}
\put(15,-2){$-2$}
\put(21,5){\line(1,0){13}}
\put(35,5){\circle{2}}
\put(33,-2){$1$}
\put(36,5){\line(1,0){13}}
\put(52,5){\circle*{1}}
\put(55,5){\circle*{1}}
\put(58,5){\circle*{1}}
\put(61,5){\line(1,0){13}}
\put(75,5){\circle{2}}
\put(76,4.5){\line(2,-1){13}}
\put(90,-2){\circle{2}}
\put(76,5.5){\line(2,1){13}}
\put(90,12){\circle{2}}
\end{picture}$, 
so the dualizing diagram method of \cite{EW} shows that the representation
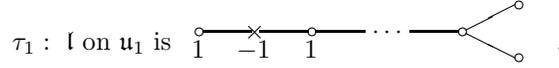
\begin{equation}\tag{\theequation$d$}
\tau_1:\,\, \gl \text{ on } \gu_1 \text{ is }
\setlength{\unitlength}{.5 mm}
\begin{picture}(100,10)
\put(5,5){\circle{2}}
\put(3,-2){$1$}
\put(6,5){\line(1,0){13}}
\put(17,3){$\times$}
\put(15,-2){$-1$}
\put(21,5){\line(1,0){13}}
\put(35,5){\circle{2}}
\put(33,-2){$1$}
\put(36,5){\line(1,0){13}}
\put(52,5){\circle*{1}}
\put(55,5){\circle*{1}}
\put(58,5){\circle*{1}}
\put(61,5){\line(1,0){13}}
\put(75,5){\circle{2}}
\put(76,4.5){\line(2,-1){13}}
\put(90,-2){\circle{2}}
\put(76,5.5){\line(2,1){13}}
\put(90,12){\circle{2}}
\end{picture}.
\end{equation}
Here $\dim \gu_2 = 1$ and $\dim \gu_1 = 4(\ell - 2)$.
}\end{noname}
\medskip

\begin{noname}\label{so-even-2} {\sc Case $SO(4,4)$.} {\rm  
Here $G_0$ is the $2$--sheeted cover of the group $SO(4,4)$
that is contained in $Spin(8;\C)$.  Its extended Dynkin diagram is
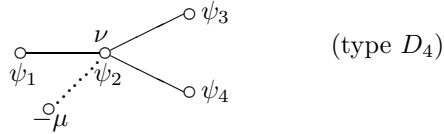
\begin{equation}\tag{\theequation$a$}
\setlength{\unitlength}{.75 mm}
\begin{picture}(90,15)
\put(10,10){\circle{2}}
\put(15,0){\circle{2}}
\put(8,5){$\psi_1$}
\put(12,-3){$-\mu$}
\put(16,1){.}
\put(17,2){.}
\put(18,3){.}
\put(19,4){.}
\put(20,5){.}
\put(21,6){.}
\put(22,7){.}
\put(23,8){.}
\put(11,10){\line(1,0){13}}
\put(23,5){$\psi_2$}
\put(23,12){$\nu$}
\put(25,10){\circle{2}}
\put(26,9.5){\line(2,-1){13}}
\put(40,3){\circle{2}}
\put(42,2){$\psi_4$}
\put(26,10.5){\line(2,1){13}}
\put(40,17){\circle{2}}
\put(42,16){$\psi_3$}
\put(65,10){(type $D_4$)}
\end{picture}
\end{equation}

\noindent Thus $\gk$ is \setlength{\unitlength}{.5 mm}
\begin{picture}(45,15)
\put(10,10){\circle{2}}
\put(10,0){\circle{2}}
\put(30,10){\circle{2}}
\put(30,0){\circle{2}}
\put(8,4){$\psi_1$}
\put(7,-5){$-\mu$}
\put(28,4){$\psi_3$}
\put(28,-5){$\psi_4$}
\end{picture}
and $\gl$ is
\setlength{\unitlength}{.5 mm}
\begin{picture}(50,15)
\put(10,5){\circle{2}}
\put(8,-1){$\psi_1$}
\put(11,5){\line(1,0){13}}
\put(22,3.3){$\times$}
\put(26,4.5){\line(2,-1){13}}
\put(40,-2){\circle{2}}
\put(42,-3){$\psi_4$}
\put(26,5.5){\line(2,1){13}}
\put(40,12){\circle{2}}
\put(42,11){$\psi_3$}
\end{picture}.  
\medskip

\noindent
Now the decompositions $\gk_0 = \gk_1 \oplus \gk_2$ and 
$\gl_0 = \gl_1 \oplus \gl_2$ are
\begin{equation}\tag{\theequation$b$}
\gk_0 = \gsp(1) \oplus \left ( \gsp(1)\oplus \gsp(1)\oplus \gsp(1) \right ) 
\text{ and } 
\gl_0 = i\R\nu^* \oplus \left ( \gsp(1)\oplus \gsp(1)\oplus \gsp(1) \right ).
\end{equation}
The representation of $\gk$ on $\gs$ has highest weight $-\nu = -\psi_2$
so its diagram is
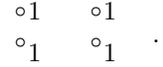
\begin{picture}(45,15)
\put(10,8){\circle{2}}
\put(10,0){\circle{2}}
\put(30,8){\circle{2}}
\put(30,0){\circle{2}}
\put(12,6){$1$}
\put(12,-5){$1$}
\put(32,6){$1$}
\put(32,-5){$1$}
\end{picture}.
Using (\ref{wtspaces}), the representation 
\begin{equation}\tag{\theequation$c$}
\tau_2:\,\, \gl \text{ on } \gu_2 \text{ is }
\setlength{\unitlength}{.5 mm}
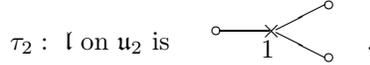
\begin{picture}(50,15)
\put(10,5){\circle{2}}
\put(11,5){\line(1,0){13}}
\put(22,3.3){$\times$}
\put(22,-2){$1$}
\put(26,4.5){\line(2,-1){13}}
\put(40,-2){\circle{2}}
\put(26,5.5){\line(2,1){13}}
\put(40,12){\circle{2}}
\end{picture}.
\end{equation}
Also, the action $\tau_{-1} \text{ of } \gl \text{ on } \gu_{-1} \text{ is }
\setlength{\unitlength}{.5 mm}
\begin{picture}(50,15)
\put(10,5){\circle{2}}
\put(5,3){$1$}
\put(11,5){\line(1,0){13}}
\put(22,3.3){$\times$}
\put(20,-2){$-2$}
\put(26,4.5){\line(2,-1){13}}
\put(40,-2){\circle{2}}
\put(44,-4){$1$}
\put(44,10){$1$}
\put(26,5.5){\line(2,1){13}}
\put(40,12){\circle{2}}
\end{picture}$
so the dualizing diagram method of \cite{EW} shows that the representation
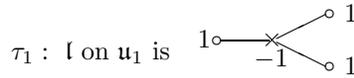
\begin{equation}\tag{\theequation$d$}
\tau_1:\,\, \gl \text{ on } \gu_1 \text{ is }
\setlength{\unitlength}{.5 mm}
\begin{picture}(50,15)
\put(10,5){\circle{2}}
\put(5,3){$1$}
\put(11,5){\line(1,0){13}}
\put(22,3.3){$\times$}
\put(20,-2){$-1$}
\put(26,4.5){\line(2,-1){13}}
\put(40,-2){\circle{2}}
\put(44,-4){$1$}
\put(44,10){$1$}
\put(26,5.5){\line(2,1){13}}
\put(40,12){\circle{2}}
\end{picture}
\end{equation}
}\end{noname}
Here $\dim \gu_2 = 1$ and $\dim \gu_1 = 8$.
\medskip

\begin{noname}\label{so-even-3} {\sc Case $Spin(2p,2\ell-2p), 2<p<\ell-2$}.
{\rm Here $G_0$ is the $2$--sheeted cover of the group $SO(2p,2\ell - 2p)$ 
that is contained in $Spin(2\ell;\C)$, with $2 < p < \ell - 2$.
Its extended Dynkin diagram is
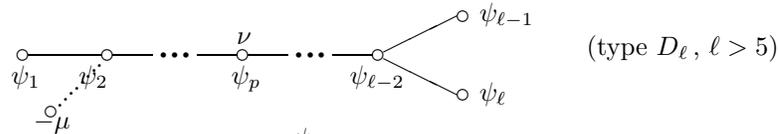
\begin{equation}\tag{\theequation$a$}
\setlength{\unitlength}{.75 mm}
\begin{picture}(140,15)
\put(10,10){\circle{2}}
\put(15,0){\circle{2}}
\put(8,5){$\psi_1$}
\put(12,-3){$-\mu$}
\put(16,1){.}
\put(17,2){.}
\put(18,3){.}
\put(19,4){.}
\put(20,5){.}
\put(21,6){.}
\put(22,7){.}
\put(23,8){.}
\put(11,10){\line(1,0){13}}
\put(25,10){\circle{2}}
\put(20,5){$\psi_2$}
\put(26,10){\line(1,0){7}}
\put(35,10){\circle*{1}}
\put(37,10){\circle*{1}}
\put(39,10){\circle*{1}}
\put(41,10){\line(1,0){7}}
\put(48,12){$\nu$}
\put(47,5){$\psi_p$}
\put(49,10){\circle{2}}
\put(50,10){\line(1,0){7}}
\put(59,10){\circle*{1}}
\put(61,10){\circle*{1}}
\put(63,10){\circle*{1}}
\put(65,10){\line(1,0){7}}
\put(73,10){\circle{2}}
\put(68,5){$\psi_{\ell - 2}$}
\put(74,9.5){\line(2,-1){13}}
\put(88,3){\circle{2}}
\put(91,2){$\psi_\ell$}
\put(74,10.5){\line(2,1){13}}
\put(88,17){\circle{2}}
\put(91,16){$\psi_{\ell - 1}$}
\put(110,10){(type $D_\ell$\,, $\ell > 5$)}
\end{picture}
\end{equation}

\noindent Thus $\gk$ is 
{\footnotesize
\setlength{\unitlength}{.4 mm}
\begin{picture}(142,15)
\put(10,10){\circle{2}}
\put(10,0){\circle{2}}
\put(08,5){$\psi_1$}
\put(07,-5){$-\mu$}
\put(11,1){\line(1,1){9}}
\put(11,10){\line(1,0){8}}
\put(20,10){\circle{2}}
\put(20,5){$\psi_2$}
\put(21,10){\line(1,0){8}}
\put(32,10){\circle*{1}}
\put(34,10){\circle*{1}}
\put(36,10){\circle*{1}}
\put(39,10){\line(1,0){8}}
\put(48,10){\circle{2}}
\put(43,5){$\psi_{p-1}$}
\put(48,10){\circle{2}}
\put(60,5){$\psi_{p+1}$}
\put(62,10){\circle{2}}
\put(63,10){\line(1,0){8}}
\put(73,10){\circle*{1}}
\put(75,10){\circle*{1}}
\put(77,10){\circle*{1}}
\put(80,10){\line(1,0){8}}
\put(89,10){\circle{2}}
\put(83,5){$\psi_{\ell - 1}$}
\put(90,10){\line(1,0){13}}
\put(104,10){\circle{2}}
\put(99,2){$\psi_{\ell - 2}$}
\put(105,9.5){\line(2,-1){13}}
\put(119,3){\circle{2}}
\put(121,2){$\psi_\ell$}
\put(105,10.5){\line(2,1){13}}
\put(119,17){\circle{2}}
\put(120,16){$\psi_{\ell - 1}$}
\end{picture}}
and $\gl$ is
{\footnotesize
\setlength{\unitlength}{.4 mm}
\begin{picture}(130,15)
\put(14,5){\circle{2}}
\put(10,0){$\psi_1$}
\put(15,5){\line(1,0){8}}
\put(24,5){\circle{2}}
\put(22,0){$\psi_2$}
\put(25,5){\line(1,0){8}}
\put(36,5){\circle*{1}}
\put(38,5){\circle*{1}}
\put(40,5){\circle*{1}}
\put(43,5){\line(1,0){8}}
\put(52,5){\circle{2}}
\put(47,0){$\psi_{p-1}$}
\put(53,5){\line(1,0){7}}
\put(59,3){$\times$}
\put(63,5){\line(1,0){7}}
\put(70,0){$\psi_{p+1}$}
\put(72,5){\circle{2}}
\put(73,5){\line(1,0){8}}
\put(84,5){\circle*{1}}
\put(86,5){\circle*{1}}
\put(88,5){\circle*{1}}
\put(90,5){\line(1,0){8}}
\put(99,5){\circle{2}}
\put(93,-2){$\psi_{\ell - 2}$}
\put(100,4.5){\line(2,-1){13}}
\put(114,-2){\circle{2}}
\put(116,-3){$\psi_\ell$}
\put(100,5.5){\line(2,1){13}}
\put(114,12){\circle{2}}
\put(116,11){$\psi_{\ell - 1}$}
\end{picture}}.
\medskip

\noindent
Now the decompositions $\gk_0 = \gk_1 \oplus \gk_2$ and 
$\gl_0 = \gl_1 \oplus \gl_2$ are
\begin{equation}\tag{\theequation$b$}
\gk_0 = \gso(2p) \oplus \gso(2\ell - 2p) \text{ and }
\gl_0 = \gu(p)  \oplus \gso(2\ell - 2p) =
i\R\nu^* \oplus \gsu(p)  \oplus \gso(2\ell - 2p).
\end{equation}
The representation of $\gk$ on $\gs$ has highest weight $-\nu = -\psi_p$:
{\footnotesize
\setlength{\unitlength}{.4 mm}
\begin{picture}(125,20)
\put(10,5){\circle{2}}
\put(10,-5){\circle{2}}
\put(11,-4){\line(1,1){9}}
\put(11,5){\line(1,0){8}}
\put(20,5){\circle{2}}
\put(21,5){\line(1,0){8}}
\put(32,5){\circle*{1}}
\put(34,5){\circle*{1}}
\put(36,5){\circle*{1}}
\put(39,5){\line(1,0){8}}
\put(48,5){\circle{2}}
\put(46,-3){$1$}
\put(60,-3){$1$}
\put(62,5){\circle{2}}
\put(63,5){\line(1,0){8}}
\put(73,5){\circle*{1}}
\put(75,5){\circle*{1}}
\put(77,5){\circle*{1}}
\put(80,5){\line(1,0){8}}
\put(89,5){\circle{2}}
\put(90,5){\line(1,0){13}}
\put(104,5){\circle{2}}
\put(105,4.5){\line(2,-1){13}}
\put(119,-2){\circle{2}}
\put(105,5.5){\line(2,1){13}}
\put(119,12){\circle{2}}
\end{picture}}.
\noindent Using (\ref{wtspaces}), the representation 
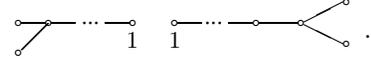
\begin{equation}\tag{\theequation$c$}
\tau_2:\,\, \gl \text{ on } \gu_2 \text{ is }
\setlength{\unitlength}{.5 mm}
\begin{picture}(130,15)
\put(14,5){\circle{2}}
\put(15,5){\line(1,0){8}}
\put(24,5){\circle{2}}
\put(22,-2){$1$}
\put(25,5){\line(1,0){8}}
\put(36,5){\circle*{1}}
\put(38,5){\circle*{1}}
\put(40,5){\circle*{1}}
\put(43,5){\line(1,0){8}}
\put(52,5){\circle{2}}
\put(53,5){\line(1,0){7}}
\put(59,3){$\times$}
\put(63,5){\line(1,0){7}}
\put(72,5){\circle{2}}
\put(73,5){\line(1,0){8}}
\put(84,5){\circle*{1}}
\put(86,5){\circle*{1}}
\put(88,5){\circle*{1}}
\put(91,5){\line(1,0){13}}
\put(105,5){\circle{2}}
\put(106,4.5){\line(2,-1){13}}
\put(120,-2){\circle{2}}
\put(106,5.5){\line(2,1){13}}
\put(120,12){\circle{2}}
\end{picture}.
\end{equation}
Also, the action $\tau_{-1}$ of $\gl$ on $\gu_{-1}$ is 
\setlength{\unitlength}{.5 mm}
\begin{picture}(120,15)
\put(4,5){\circle{2}}
\put(5,5){\line(1,0){8}}
\put(14,5){\circle{2}}
\put(15,5){\line(1,0){8}}
\put(26,5){\circle*{1}}
\put(28,5){\circle*{1}}
\put(30,5){\circle*{1}}
\put(33,5){\line(1,0){8}}
\put(42,5){\circle{2}}
\put(40,-2){$1$}
\put(43,5){\line(1,0){7}}
\put(49,3){$\times$}
\put(45,-2){$-2$}
\put(53,5){\line(1,0){7}}
\put(62,5){\circle{2}}
\put(60,-2){$1$}
\put(63,5){\line(1,0){8}}
\put(74,5){\circle*{1}}
\put(76,5){\circle*{1}}
\put(78,5){\circle*{1}}
\put(81,5){\line(1,0){13}}
\put(95,5){\circle{2}}
\put(96,4.5){\line(2,-1){13}}
\put(110,-2){\circle{2}}
\put(96,5.5){\line(2,1){13}}
\put(110,12){\circle{2}}
\end{picture}
so the dualizing diagram method of \cite{EW} shows that the representation
\begin{equation}\tag{\theequation$d$}
\tau_1:\,\, \gl \text{ on } \gu_1 \text{ is }
\setlength{\unitlength}{.5 mm}
\begin{picture}(130,15)
\put(14,5){\circle{2}}
\put(12,-2){$1$}
\put(15,5){\line(1,0){8}}
\put(24,5){\circle{2}}
\put(25,5){\line(1,0){8}}
\put(36,5){\circle*{1}}
\put(38,5){\circle*{1}}
\put(40,5){\circle*{1}}
\put(43,5){\line(1,0){8}}
\put(52,5){\circle{2}}
\put(53,5){\line(1,0){7}}
\put(59,3){$\times$}
\put(55,-2){$-1$}
\put(63,5){\line(1,0){7}}
\put(72,5){\circle{2}}
\put(72,-2){$1$}
\put(73,5){\line(1,0){8}}
\put(84,5){\circle*{1}}
\put(86,5){\circle*{1}}
\put(88,5){\circle*{1}}
\put(91,5){\line(1,0){13}}
\put(105,5){\circle{2}}
\put(106,4.5){\line(2,-1){13}}
\put(120,-2){\circle{2}}
\put(106,5.5){\line(2,1){13}}
\put(120,12){\circle{2}}
\end{picture}.
\end{equation}
\medskip
}\end{noname}
\noindent
Here $\dim \gu_2 = p(p-1)/2$ and $\dim \gu_1 = 2p(\ell -p)$.
This exhausts the cases where $\gg$ is of type $D$, and thus exhausts 
the classical cases.  We go on the the exceptional cases.
\medskip

\begin{noname}\label{g2} {\sc Case $G_{2,A_1A_1}$}. {\rm
Here $G_0$ is the split real Lie group of type $G_2$.
It has maximal compact subgroup $SO(4)$.  Its extended Dynkin diagram is
\begin{equation}\tag{\theequation$a$}
\setlength{\unitlength}{.75 mm}
\begin{picture}(50,8)
\put(5,3){\circle*{2}}
\put(3,-2){$\psi_1$}
\put(6,2){\line(1,0){13}}
\put(6,3){\line(1,0){13}}
\put(6,4){\line(1,0){13}}
\put(20,3){\circle{2}}
\put(18,-2){$\psi_2$}
\put(18,6){$\nu$}
\put(21,3){.}
\put(22,3){.}
\put(23,3){.}
\put(24,3){.}
\put(25,3){.}
\put(26,3){.}
\put(27,3){.}
\put(28,3){.}
\put(29,3){.}
\put(30,3){.}
\put(31,3){.}
\put(32,3){.}
\put(33,3){.}
\put(35,3){\circle{2}}
\put(32,-2){$-\mu$}
\end{picture}
\text{(type $G_2$)}.
\end{equation}
\noindent

Thus $\gk$ is 
\setlength{\unitlength}{.5 mm}
\begin{picture}(35,8)
\put(5,3){\circle*{2}}
\put(3,-2){$\psi_1$}
\put(25,3){\circle{2}}
\put(22,-2){$-\mu$}
\end{picture}
and $\gl$ is
\setlength{\unitlength}{.5 mm}
\begin{picture}(25,8)
\put(5,3){\circle*{2}}
\put(3,-2){$\psi_1$}
\put(6,2){\line(1,0){13}}
\put(6,3){\line(1,0){13}}
\put(6,4){\line(1,0){13}}
\put(19.2,1.6){\tiny $\boxtimes$}
\end{picture}.
\medskip

\noindent
Now the decompositions $\gk_0 = \gk_1 \oplus \gk_2$ and 
$\gl_0 = \gl_1 \oplus \gl_2$ are
\begin{equation}\tag{\theequation$b$}
\gk_0 = \gsp(1) \oplus \gsp(1) \text{ and } \gl_0 = i\R\nu^* \oplus \gsp(1).
\end{equation}
The representation of $\gk$ on $\gs$ has highest weight $-\nu = -\psi_2$:
\setlength{\unitlength}{.5 mm}
\begin{picture}(35,11)
\put(5,6){\circle*{2}}
\put(3,-2){$3$}
\put(25,6){\circle{2}}
\put(23,-2){$1$}
\end{picture}.
\noindent Using (\ref{wtspaces}), the representation
\begin{equation}\tag{\theequation$c$}
\tau_2:\,\, \gl \text{ on } \gu_2 \text{ is }
\setlength{\unitlength}{.75 mm}
\begin{picture}(25,8)
\put(5,3){\circle*{2}}
\put(6,2){\line(1,0){13}}
\put(6,3){\line(1,0){13}}
\put(6,4){\line(1,0){13}}
\put(19.2,1.8){\tiny $\boxtimes$}
\put(19,-2){$1$}
\end{picture}
\end{equation}
Also, the action $\tau_{-1}$ of $\gl$ on $\gu_{-1}$ is
\setlength{\unitlength}{.5 mm}
\begin{picture}(27,8)
\put(5,5){\circle*{2}}
\put(5,-2){$3$}
\put(6,4){\line(1,0){13}}
\put(6,5){\line(1,0){13}}
\put(6,6){\line(1,0){13}}
\put(19.2,3.3){\tiny $\boxtimes$}
\put(16,-2){$-2$}
\end{picture}
so the representation
\begin{equation}\tag{\theequation$d$}
\tau_1:\,\, \gl \text{ on } \gu_1 \text{ is }
\setlength{\unitlength}{.75 mm}
\begin{picture}(30,8)
\put(5,3){\circle*{2}}
\put(4.5,-2.5){$3$}
\put(6,2){\line(1,0){13}}
\put(6,3){\line(1,0){13}}
\put(6,4){\line(1,0){13}}
\put(19.2,1.8){\tiny $\boxtimes$}
\put(17.5,-2.5){$-1$}
\end{picture}.
\end{equation}
} \end{noname}
\medskip

\noindent Note that $\tau_1|_{[\gl,\gl]}$ has degree $4$, is self--dual, and 
has an antisymmetric bilinear invariant.  Also, $\dim \gu_2 = 1$ and
$\tau_2|_{[\gl,\gl]}$ is trivial, so that bilinear invariant is given by the
Lie algebra product $\gu_1 \times \gu_1 \to \gu_2$.
\medskip

\begin{noname}\label{f4-1} {\sc Case $F_{4,A_1C_3}$}.  {\rm
Here $G_0$ is the simply connected real Lie group of type $F_4$
whose maximal compact subgroup has $2$--sheeted cover $Sp(1)\times Sp(3)$.  
Its extended Dynkin diagram is
\begin{equation}\tag{\theequation$a$}
\setlength{\unitlength}{.75 mm}
\begin{picture}(90,12)
\put(5,5){\circle{2}}
\put(1,0){$-\mu$}
\put(6,5){.}
\put(7,5){.}
\put(8,5){.}
\put(9,5){.}
\put(10,5){.}
\put(11,5){.}
\put(12,5){.}
\put(13,5){.}
\put(14,5){.}
\put(15,5){.}
\put(16,5){.}
\put(17,5){.}
\put(18,5){.}
\put(20,5){\circle{2}}
\put(18,0){$\psi_1$}
\put(19,7){$\nu$}
\put(21,5){\line(1,0){13}}
\put(35,5){\circle{2}}
\put(33,0){$\psi_2$}
\put(36,4.5){\line(1,0){13}}
\put(36,5.5){\line(1,0){13}}
\put(50,5){\circle*{2}}
\put(48,0){$\psi_3$}
\put(51,5){\line(1,0){13}}
\put(65,5){\circle*{2}}
\put(63,0){$\psi_4$}
\put(80,2){(Type $F_4$)}
\end{picture}
\end{equation}
Thus $\gk$ is 
\setlength{\unitlength}{.5 mm}
\begin{picture}(65,12)
\put(5,5){\circle{2}}
\put(1,0){$-\mu$}
\put(25,5){\circle{2}}
\put(23,0){$\psi_2$}
\put(26,4.5){\line(1,0){13}}
\put(26,5.5){\line(1,0){13}}
\put(40,5){\circle*{2}}
\put(38,0){$\psi_3$}
\put(41,5){\line(1,0){13}}
\put(55,5){\circle*{2}}
\put(53,0){$\psi_4$}
\end{picture}
and $\gl$ is
\setlength{\unitlength}{.5 mm}
\begin{picture}(65,12)
\put(8,3){$\times$}
\put(11,5){\line(1,0){13}}
\put(25,5){\circle{2}}
\put(23,0){$\psi_2$}
\put(26,4.5){\line(1,0){13}}
\put(26,5.5){\line(1,0){13}}
\put(40,5){\circle*{2}}
\put(38,0){$\psi_3$}
\put(41,5){\line(1,0){13}}
\put(55,5){\circle*{2}}
\put(53,0){$\psi_4$}
\end{picture}.
\medskip

\noindent
Now the decompositions $\gk_0 = \gk_1 \oplus \gk_2$ and 
$\gl_0 = \gl_1 \oplus \gl_2$ are
\begin{equation}\tag{\theequation$b$}
\gk_0 = \gsp(1)\oplus \gsp(3) \text{ and } \gl_0 = i\R\nu^* \oplus \gsp(3).
\end{equation}
The representation of $\gk$ on $\gs$ has highest weight $-\nu = -\psi_1$:
\setlength{\unitlength}{.5 mm}
\begin{picture}(65,12)
\put(5,6){\circle{2}}
\put(4,-1){$1$}
\put(24,-1){$1$}
\put(25,6){\circle{2}}
\put(26,5.5){\line(1,0){13}}
\put(26,6.5){\line(1,0){13}}
\put(40,6){\circle*{2}}
\put(41,6){\line(1,0){13}}
\put(55,6){\circle*{2}}
\end{picture}.
\noindent Using (\ref{wtspaces}), the representation
\begin{equation}\tag{\theequation$c$}
\tau_2:\,\, \gl \text{ on } \gu_2 \text{ is }
\setlength{\unitlength}{.75 mm}
\begin{picture}(65,12)
\put(8,0){$1$}
\put(8,4){$\times$}
\put(11,5){\line(1,0){13}}
\put(25,5){\circle{2}}
\put(26,4.5){\line(1,0){13}}
\put(26,5.5){\line(1,0){13}}
\put(40,5){\circle*{2}}
\put(41,5){\line(1,0){13}}
\put(55,5){\circle*{2}}
\end{picture}
\end{equation}
Also, the action $\tau_{-1}$ of $\gl$ on $\gu_{-1}$ is
\setlength{\unitlength}{.5 mm}
\begin{picture}(65,12)
\put(7,-1){$-2$}
\put(8,4){$\times$}
\put(11,6){\line(1,0){13}}
\put(24,-1){$1$}
\put(25,6){\circle{2}}
\put(26,5.5){\line(1,0){13}}
\put(26,6.5){\line(1,0){13}}
\put(40,6){\circle*{2}}
\put(41,6){\line(1,0){13}}
\put(55,6){\circle*{2}}
\end{picture}
so the representation
\begin{equation}\tag{\theequation$d$}
\tau_1:\,\, \gl \text{ on } \gu_1 \text{ is }
\setlength{\unitlength}{.75 mm}
\begin{picture}(65,12)
\put(8,4){$\times$}
\put(7,-1){$-1$}
\put(11,5){\line(1,0){13}}
\put(24,-1){$1$}
\put(25,5){\circle{2}}
\put(26,4.5){\line(1,0){13}}
\put(26,5.5){\line(1,0){13}}
\put(40,5){\circle*{2}}
\put(41,5){\line(1,0){13}}
\put(55,5){\circle*{2}}
\end{picture}.
\end{equation}
}\end{noname}
\medskip

\noindent Note that $\tau_1|_{[\gl,\gl]}$ has degree $14$, is self--dual, and 
has an antisymmetric bilinear invariant.  Also, $\dim \gu_2 = 1$ and
$\tau_2|_{[\gl,\gl]}$ is trivial, so that bilinear invariant is given by the
Lie algebra product $\gu_1 \times \gu_1 \to \gu_2$.
\medskip

\begin{noname}\label{f4-2} {\sc Case $F_{4,B_4}$}.  {\rm
Here $G_0$ is the simply connected real Lie group of type $F_4$
with maximal compact subgroup $Spin(9)$.  Its extended 
Dynkin diagram is
\begin{equation}\tag{\theequation$a$}
\setlength{\unitlength}{.75 mm}
\begin{picture}(90,12)
\put(5,5){\circle{2}}
\put(1,0){$-\mu$}
\put(6,5){.}
\put(7,5){.}
\put(8,5){.}
\put(9,5){.}
\put(10,5){.}
\put(11,5){.}
\put(12,5){.}
\put(13,5){.}
\put(14,5){.}
\put(15,5){.}
\put(16,5){.}
\put(17,5){.}
\put(18,5){.}
\put(20,5){\circle{2}}
\put(18,0){$\psi_1$}
\put(21,5){\line(1,0){13}}
\put(35,5){\circle{2}}
\put(33,0){$\psi_2$}
\put(36,4.5){\line(1,0){13}}
\put(36,5.5){\line(1,0){13}}
\put(50,5){\circle*{2}}
\put(48,0){$\psi_3$}
\put(64,7){$\nu$}
\put(51,5){\line(1,0){13}}
\put(65,5){\circle*{2}}
\put(63,0){$\psi_4$}
\put(80,2){(Type $F_4$)}
\end{picture}
\end{equation}
Thus $\gk$ is 
\setlength{\unitlength}{.5 mm}
\begin{picture}(60,12)
\put(1,6){$-\mu$}
\put(5,2){\circle{2}}
\put(6,2){\line(1,0){13}}
\put(18,6){$\psi_1$}
\put(20,2){\circle{2}}
\put(21,2){\line(1,0){13}}
\put(33,6){$\psi_2$}
\put(35,2){\circle{2}}
\put(36,1.5){\line(1,0){13}}
\put(36,2.5){\line(1,0){13}}
\put(48,6){$\psi_3$}
\put(50,2){\circle*{2}}
\end{picture}
and $\gl$ is
\setlength{\unitlength}{.5 mm}
\begin{picture}(65,12)
\put(10,2){\circle{2}}
\put(8,6){$\psi_1$}
\put(11,2){\line(1,0){13}}
\put(25,2){\circle{2}}
\put(23,6){$\psi_2$}
\put(26,1.5){\line(1,0){13}}
\put(26,2.5){\line(1,0){13}}
\put(40,2){\circle*{2}}
\put(38,6){$\psi_3$}
\put(41,2){\line(1,0){13}}
\put(54.2,0.8){\tiny $\boxtimes$}
\end{picture}.
\medskip

\noindent
Now the decompositions $\gk_0 = \gk_1 \oplus \gk_2$ and 
$\gl_0 = \gl_1 \oplus \gl_2$ are
\begin{equation}\tag{\theequation$b$}
\gk_0 = \gso(9) \text{ and } \gl_0 = i\R\nu^* \oplus \gso(7).
\end{equation}
The representation of $\gk$ on $\gs$ has highest weight $-\nu = -\psi_4$:
\setlength{\unitlength}{.5 mm}
\begin{picture}(60,12)
\put(5,2){\circle{2}}
\put(6,2){\line(1,0){13}}
\put(20,2){\circle{2}}
\put(21,2){\line(1,0){13}}
\put(35,2){\circle{2}}
\put(36,1.5){\line(1,0){13}}
\put(36,2.5){\line(1,0){13}}
\put(50,2){\circle*{2}}
\put(48,5){$1$}
\end{picture}.
\noindent Using (\ref{wtspaces}), the representation
\begin{equation}\tag{\theequation$c$}
\tau_2:\,\, \gl \text{ on } \gu_2 \text{ is }
\setlength{\unitlength}{.75 mm}
\begin{picture}(65,12)
\put(10,5){\circle{2}}
\put(9,0){$1$}
\put(11,5){\line(1,0){13}}
\put(25,5){\circle{2}}
\put(26,4.5){\line(1,0){13}}
\put(26,5.5){\line(1,0){13}}
\put(40,5){\circle*{2}}
\put(41,5){\line(1,0){13}}
\put(54.2,3.8){\tiny $\boxtimes$}
\end{picture}
\end{equation}
Also, the action $\tau_{-1}$ of $\gl$ on $\gu_{-1}$ is
\setlength{\unitlength}{.5 mm}
\begin{picture}(65,12)
\put(10,2){\circle{2}}
\put(11,2){\line(1,0){13}}
\put(25,2){\circle{2}}
\put(26,1.5){\line(1,0){13}}
\put(26,2.5){\line(1,0){13}}
\put(40,2){\circle*{2}}
\put(39,6){$1$}
\put(41,2){\line(1,0){13}}
\put(54.2,0.8){\tiny $\boxtimes$}
\put(52,6){$-2$}
\end{picture}
so the representation
\begin{equation}\tag{\theequation$d$}
\tau_1:\,\, \gl \text{ on } \gu_1 \text{ is }
\setlength{\unitlength}{.75 mm}
\begin{picture}(65,12)
\put(10,5){\circle{2}}
\put(11,5){\line(1,0){13}}
\put(25,5){\circle{2}}
\put(26,4.5){\line(1,0){13}}
\put(26,5.5){\line(1,0){13}}
\put(40,5){\circle*{2}}
\put(39,0){$1$}
\put(41,5){\line(1,0){13}}
\put(54.2,3.8){\tiny $\boxtimes$}
\put(51,0){$-1$}
\end{picture}.
\end{equation}
}
\end{noname}
\noindent Note that $\tau_1|_{[\gl,\gl]}$ has degree $8$, is self--dual, 
and has a symmetric bilinear invariant.  In effect $\tau_1$ is the action
of $Spin(7)$ on the Cayley numbers, and $\tau_2$ is its action (factored
through $SO(7)$ on the pure imaginary Cayley numbers.  thus
$\dim \gu_2 = 7$ and $\dim \gu_1 = 8$.
\medskip

\begin{noname}\label{e6-1} {\sc Case $E_{6,A_1A_5,1}$}. {\rm
Here $G_0$ is the group of type $E_6$ whose maximal compact subgroup
is the $2$--sheeted cover of $SU(2)\times SU(6)$.  The noncompact simple 
root $\nu = \psi_3$, so 
$L$ is of type $T_1A_1A_4$.  We do not consider the case $\nu = \psi_5$
separately because the two differ only by an outer automorphism of $E_6$.
Here the extended Dynkin diagram (Bourbaki root order) is
\begin{equation}\tag{\theequation$a$}
\setlength{\unitlength}{.75 mm}
\begin{picture}(100,20)
\put(5,6){\circle{2}}
\put(3,9){$\psi_1$}
\put(6,6){\line(1,0){13}}
\put(20,6){\circle{2}}
\put(18,9){$\psi_3$}
\put(19,1){$\nu$}
\put(21,6){\line(1,0){13}}
\put(35,6){\circle{2}}
\put(33,9){$\psi_4$}
\put(50,13){\circle{2}}
\put(48,16){$\psi_5$}
\put(51,13){\line(1,0){13}}
\put(65,13){\circle{2}}
\put(63,16){$\psi_6$}
\put(36,5.5){\line(2,-1){13}}
\put(36,6.5){\line(2,1){13}}
\put(50,-1){\circle{2}}
\put(48,2){$\psi_2$}
\put(51,-1){.}
\put(52,-1){.}
\put(53,-1){.}
\put(54,-1){.}
\put(55,-1){.}
\put(56,-1){.}
\put(57,-1){.}
\put(58,-1){.}
\put(59,-1){.}
\put(60,-1){.}
\put(61,-1){.}
\put(62,-1){.}
\put(63,-1){.}
\put(65,-1){\circle{2}}
\put(61,2){$-\mu$}
\put(85,5){(Type $E_6$)}
\end{picture}
\end{equation}
Thus $\gk$ is 
\setlength{\unitlength}{.5 mm}
\begin{picture}(65,17)
\put(5,2){\circle{2}}
\put(3,5){$\psi_1$}
\put(25,2){\circle{2}}
\put(21,5){$\psi_4$}
\put(40,9){\circle{2}}
\put(38,12){$\psi_5$}
\put(41,9){\line(1,0){13}}
\put(55,9){\circle{2}}
\put(53,12){$\psi_6$}
\put(26,1.5){\line(2,-1){13}}
\put(26,2.5){\line(2,1){13}}
\put(40,-5){\circle{2}}
\put(38,-2){$\psi_2$}
\put(41,-5){\line(1,0){13}}
\put(55,-5){\circle{2}}
\put(51,-2){$-\mu$}
\end{picture}
and $\gl$ is
\setlength{\unitlength}{.5 mm}
\begin{picture}(75,20)
\put(5,2){\circle{2}}
\put(3,5){$\psi_1$}
\put(6,2){\line(1,0){13}}
\put(17.5,1){$\times$}
\put(21,2){\line(1,0){13}}
\put(35,2){\circle{2}}
\put(32,6){$\psi_4$}
\put(50,9){\circle{2}}
\put(48,12){$\psi_5$}
\put(51,9){\line(1,0){13}}
\put(65,9){\circle{2}}
\put(63,12){$\psi_6$}
\put(36,1.5){\line(2,-1){13}}
\put(36,2.5){\line(2,1){13}}
\put(50,-5){\circle{2}}
\put(48,-2){$\psi_2$}
\end{picture}.
\bigskip

\noindent
Now the decompositions $\gk_0 = \gk_1 \oplus \gk_2$ and 
$\gl_0 = \gl_1 \oplus \gl_2$ are
\begin{equation}\tag{\theequation$b$}
\gk_0 = \gsu(6) \oplus \gsu(2) \text{ and } \gl_0 = 
 \bigl ( \gsu(5)\oplus i\R\nu^*  \bigr )  \oplus \gsu(2).
\end{equation}
The representation of $\gk$ on $\gs$ 
\medskip

\noindent has highest weight $-\nu = -\psi_3$:
\setlength{\unitlength}{.5 mm}
\begin{picture}(65,17)
\put(5,2){\circle{2}}
\put(4,5){$1$}
\put(25,2){\circle{2}}
\put(24,5){$1$}
\put(40,9){\circle{2}}
\put(41,9){\line(1,0){13}}
\put(55,9){\circle{2}}
\put(26,1.5){\line(2,-1){13}}
\put(26,2.5){\line(2,1){13}}
\put(40,-5){\circle{2}}
\put(41,-5){\line(1,0){13}}
\put(55,-5){\circle{2}}
\end{picture}.  
\noindent Using (\ref{wtspaces}), the representation
\begin{equation}\tag{\theequation$c$}
\tau_2:\,\, \gl \text{ on } \gu_2 \text{ is }
\setlength{\unitlength}{.75 mm}
\begin{picture}(75,15)
\put(5,2){\circle{2}}
\put(6,2){\line(1,0){13}}
\put(18,1){$\times$}
\put(21,2){\line(1,0){13}}
\put(35,2){\circle{2}}
\put(50,9){\circle{2}}
\put(51,9){\line(1,0){13}}
\put(65,9){\circle{2}}
\put(36,1.5){\line(2,-1){13}}
\put(36,2.5){\line(2,1){13}}
\put(50,-5){\circle{2}}
\put(48,-2){$1$}
\end{picture}.
\end{equation}
Also, the action $\tau_{-1}$ of $\gl$ on $\gu_{-1}$ is
\setlength{\unitlength}{.5 mm}
\begin{picture}(75,15)
\put(5,2){\circle{2}}
\put(4,5){$1$}
\put(6,2){\line(1,0){13}}
\put(17.5,1){$\times$}
\put(16,5){$-2$}
\put(21,2){\line(1,0){13}}
\put(35,2){\circle{2}}
\put(34,5){$1$}
\put(50,9){\circle{2}}
\put(51,9){\line(1,0){13}}
\put(65,9){\circle{2}}
\put(36,1.5){\line(2,-1){13}}
\put(36,2.5){\line(2,1){13}}
\put(50,-5){\circle{2}}
\end{picture}
so the representation
\begin{equation}\tag{\theequation$d$}
\tau_1:\,\, \gl \text{ on } \gu_1 \text{ is }
\setlength{\unitlength}{.75 mm}
\begin{picture}(75,15)
\put(5,2){\circle{2}}
\put(4,4){$1$}
\put(6,2){\line(1,0){13}}
\put(18,1){$\times$}
\put(17,4){$-1$}
\put(21,2){\line(1,0){13}}
\put(35,2){\circle{2}}
\put(49,11){$1$}
\put(50,9){\circle{2}}
\put(51,9){\line(1,0){13}}
\put(65,9){\circle{2}}
\put(36,1.5){\line(2,-1){13}}
\put(36,2.5){\line(2,1){13}}
\put(50,-5){\circle{2}}
\end{picture}.
\end{equation}
}\end{noname}
\noindent Note that $\tau_1|_{[\gl,\gl]}$ has degree $20$ and is not self--dual.
We have $\dim \gu_2 = 5$ and $\dim \gu_1 = 20$.
\medskip

\begin{noname}\label{e6-2} {\sc Case $E_{6,A_1A_5,2}$}.  {\rm
Here $G_0$ is the group of type $E_6$ with maximal compact subgroup
$SU(2)\times SU(6)$.  The noncompact simple root $\nu = \psi_2$, so 
$L$ is of type $T_1A_5$.  The extended Dynkin diagram (Bourbaki root order) is
\begin{equation}\tag{\theequation$a$}
\setlength{\unitlength}{.75 mm}
\begin{picture}(100,20)
\put(5,6){\circle{2}}
\put(3,9){$\psi_1$}
\put(6,6){\line(1,0){13}}
\put(20,6){\circle{2}}
\put(18,9){$\psi_3$}
\put(21,6){\line(1,0){13}}
\put(35,6){\circle{2}}
\put(33,9){$\psi_4$}
\put(50,13){\circle{2}}
\put(48,16){$\psi_5$}
\put(51,13){\line(1,0){13}}
\put(65,13){\circle{2}}
\put(63,16){$\psi_6$}
\put(36,5.5){\line(2,-1){13}}
\put(36,6.5){\line(2,1){13}}
\put(50,-1){\circle{2}}
\put(48,2){$\psi_2$}
\put(49,-5){$\nu$}
\put(51,-1){.}
\put(52,-1){.}
\put(53,-1){.}
\put(54,-1){.}
\put(55,-1){.}
\put(56,-1){.}
\put(57,-1){.}
\put(58,-1){.}
\put(59,-1){.}
\put(60,-1){.}
\put(61,-1){.}
\put(62,-1){.}
\put(63,-1){.}
\put(65,-1){\circle{2}}
\put(61,2){$-\mu$}
\put(85,5){(Type $E_6$)}
\end{picture}
\end{equation}
Thus $\gk$ is 
\setlength{\unitlength}{.5 mm}
\begin{picture}(75,17)
\put(5,2){\circle{2}}
\put(3,5){$\psi_1$}
\put(6,2){\line(1,0){13}}
\put(20,2){\circle{2}}
\put(18,5){$\psi_3$}
\put(21,2){\line(1,0){13}}
\put(35,2){\circle{2}}
\put(31,5){$\psi_4$}
\put(50,9){\circle{2}}
\put(48,12){$\psi_5$}
\put(51,9){\line(1,0){13}}
\put(65,9){\circle{2}}
\put(63,12){$\psi_6$}
\put(36,2.5){\line(2,1){13}}
\put(65,-5){\circle{2}}
\put(61,-2){$-\mu$}
\end{picture}
and $\gl$ is
\setlength{\unitlength}{.5 mm}
\begin{picture}(75,20)
\put(5,2){\circle{2}}
\put(3,5){$\psi_1$}
\put(6,2){\line(1,0){13}}
\put(20,2){\circle{2}}
\put(18,5){$\psi_3$}
\put(21,2){\line(1,0){13}}
\put(35,2){\circle{2}}
\put(32,6){$\psi_4$}
\put(50,9){\circle{2}}
\put(48,12){$\psi_5$}
\put(51,9){\line(1,0){13}}
\put(65,9){\circle{2}}
\put(63,12){$\psi_6$}
\put(36,1.5){\line(2,-1){13}}
\put(36,2.5){\line(2,1){13}}
\put(48,-6.2){$\times$}
\end{picture}.
\bigskip

\noindent
Now the decompositions $\gk_0 = \gk_1 \oplus \gk_2$ and 
$\gl_0 = \gl_1 \oplus \gl_2$ are
\begin{equation}\tag{\theequation$b$}
\gk_0 = \gsp(1)\oplus \gsu(6) \text{ and } \gl_0 = i\R\nu^* \oplus \gsu(6).
\end{equation}
The representation of $\gk$ on $\gs$ has highest weight $-\nu = -\psi_2$:
\setlength{\unitlength}{.5 mm}
\setlength{\unitlength}{.5 mm}
\begin{picture}(75,17)
\put(5,2){\circle{2}}
\put(6,2){\line(1,0){13}}
\put(20,2){\circle{2}}
\put(21,2){\line(1,0){13}}
\put(35,2){\circle{2}}
\put(33,5){$1$}
\put(50,9){\circle{2}}
\put(51,9){\line(1,0){13}}
\put(65,9){\circle{2}}
\put(36,2.5){\line(2,1){13}}
\put(65,-5){\circle{2}}
\put(63,-2){$1$}
\end{picture}.  
\noindent Using (\ref{wtspaces}), the representation
\begin{equation}\tag{\theequation$c$}
\tau_2:\,\, \gl \text{ on } \gu_2 \text{ is }
\setlength{\unitlength}{.75 mm}
\begin{picture}(75,15)
\put(5,2){\circle{2}}
\put(6,2){\line(1,0){13}}
\put(20,2){\circle{2}}
\put(21,2){\line(1,0){13}}
\put(35,2){\circle{2}}
\put(50,9){\circle{2}}
\put(51,9){\line(1,0){13}}
\put(65,9){\circle{2}}
\put(36,1.5){\line(2,-1){13}}
\put(36,2.5){\line(2,1){13}}
\put(48,-6){$\times$}
\put(48,-2){$1$}
\end{picture}.
\end{equation}
Also, the action $\tau_{-1}$ of $\gl$ on $\gu_{-1}$ is
\setlength{\unitlength}{.5 mm}
\begin{picture}(75,15)
\put(5,2){\circle{2}}
\put(6,2){\line(1,0){13}}
\put(20,2){\circle{2}}
\put(21,2){\line(1,0){13}}
\put(35,2){\circle{2}}
\put(34,5){$1$}
\put(50,9){\circle{2}}
\put(51,9){\line(1,0){13}}
\put(65,9){\circle{2}}
\put(36,1.5){\line(2,-1){13}}
\put(36,2.5){\line(2,1){13}}
\put(48,-6){$\times$}
\put(46,-2){$-2$}
\end{picture}
so the representation
\begin{equation}\tag{\theequation$d$}
\tau_1:\,\, \gl \text{ on } \gu_1 \text{ is }
\setlength{\unitlength}{.75 mm}
\begin{picture}(75,15)
\put(5,2){\circle{2}}
\put(6,2){\line(1,0){13}}
\put(20,2){\circle{2}}
\put(21,2){\line(1,0){13}}
\put(35,2){\circle{2}}
\put(34,5){$1$}
\put(50,9){\circle{2}}
\put(51,9){\line(1,0){13}}
\put(65,9){\circle{2}}
\put(36,1.5){\line(2,-1){13}}
\put(36,2.5){\line(2,1){13}}
\put(48,-6){$\times$}
\put(46,-3){$-1$}
\end{picture}.
\end{equation}
}\end{noname}
\noindent Note that $\tau_1|_{[\gl,\gl]}$ has degree $20$, is self--dual, and 
has an antisymmetric bilinear invariant. Also, $\dim \gu_2 = 1$ and
$\tau_2|_{[\gl,\gl]}$ is trivial, so that bilinear invariant is given by the
Lie algebra product $\gu_1 \times \gu_1 \to \gu_2$.  In brief,
$\dim \gu_2 = 1$ and $\dim \gu_1 = 20$.
\medskip

\begin{noname}\label{e7-1} {\sc Case $E_{7,A_1D_6,1}$}.  
{\rm Here $G_0$ is the group of type $E_7$ with maximal compact subgroup 
that is the $2$--sheeted cover of $SU(2)\times Spin(12)$.  The noncompact 
simple root $\nu = \psi_1$, so $L$ is of type $T_1D_6$ and the extended Dynkin
diagram is

\begin{equation}\tag{\theequation$a$}
\setlength{\unitlength}{.75 mm}
\begin{picture}(140,18)
\put(5,16){\circle{2}}
\put(2,12){$-\mu$}
\put(6,16){.}
\put(7,16){.}
\put(8,16){.}
\put(9,16){.}
\put(10,16){.}
\put(11,16){.}
\put(12,16){.}
\put(13,16){.}
\put(14,16){.}
\put(15,16){.}
\put(16,16){.}
\put(17,16){.}
\put(18,16){.}
\put(20,16){\circle{2}}
\put(18,12){$\psi_1$}
\put(19,18){$\nu$}
\put(21,16){\line(1,0){13}}
\put(35,16){\circle{2}}
\put(33,12){$\psi_3$}
\put(36,16){\line(1,0){13}}
\put(50,16){\circle{2}}
\put(52,12){$\psi_4$}
\put(51,16){\line(1,0){13}}
\put(65,16){\circle{2}}
\put(63,12){$\psi_5$}
\put(66,16){\line(1,0){13}}
\put(80,16){\circle{2}}
\put(78,12){$\psi_6$}
\put(81,16){\line(1,0){13}}
\put(95,16){\circle{2}}
\put(93,12){$\psi_7$}
\put(50,15){\line(0,-1){13}}
\put(50,1){\circle{2}}
\put(52,1){$\psi_2$}
\put(110,10){(Type $E_7$)}
\end{picture}
\end{equation}
Thus $\gk$ is
\setlength{\unitlength}{.5 mm}
\begin{picture}(110,18)
\put(10,11){\circle{2}}
\put(7,6){$-\mu$}
\put(35,11){\circle{2}}
\put(33,5){$\psi_3$}
\put(36,11){\line(1,0){13}}
\put(50,11){\circle{2}}
\put(52,5){$\psi_4$}
\put(51,11){\line(1,0){13}}
\put(65,11){\circle{2}}
\put(63,5){$\psi_5$}
\put(66,11){\line(1,0){13}}
\put(80,11){\circle{2}}
\put(78,5){$\psi_6$}
\put(81,11){\line(1,0){13}}
\put(95,11){\circle{2}}
\put(93,5){$\psi_7$}
\put(50,10){\line(0,-1){13}}
\put(50,-4){\circle{2}}
\put(52,-6){$\psi_2$}
\end{picture}
and $\gl$ is
\setlength{\unitlength}{.5 mm}
\begin{picture}(95,18)
\put(3,9){$\times$}
\put(6,11){\line(1,0){13}}
\put(20,11){\circle{2}}
\put(18,5){$\psi_3$}
\put(21,11){\line(1,0){13}}
\put(35,11){\circle{2}}
\put(37,5){$\psi_4$}
\put(36,11){\line(1,0){13}}
\put(50,11){\circle{2}}
\put(48,5){$\psi_5$}
\put(51,11){\line(1,0){13}}
\put(65,11){\circle{2}}
\put(63,5){$\psi_6$}
\put(66,11){\line(1,0){13}}
\put(80,11){\circle{2}}
\put(78,5){$\psi_7$}
\put(35,10){\line(0,-1){13}}
\put(35,-4){\circle{2}}
\put(37,-6){$\psi_2$}
\end{picture}.
\medskip

\noindent
Now the decompositions $\gk_0 = \gk_1 \oplus \gk_2$ and 
$\gl_0 = \gl_1 \oplus \gl_2$ are
\begin{equation}\tag{\theequation$b$}
\gk_0 = \gsp(1) \oplus \gso(12) \text{ and } 
\gl_0 = i\R\nu^* \oplus \gso(12).
\end{equation}
The representation of $\gk$ on $\gs$ has highest weight $-\nu = -\psi_1$:
\setlength{\unitlength}{.4 mm}
\begin{picture}(95,20)
\put(5,11){\circle{2}}
\put(4,3){$1$}
\put(25,11){\circle{2}}
\put(24,3){$1$}
\put(26,11){\line(1,0){13}}
\put(40,11){\circle{2}}
\put(41,11){\line(1,0){13}}
\put(55,11){\circle{2}}
\put(56,11){\line(1,0){13}}
\put(70,11){\circle{2}}
\put(71,11){\line(1,0){13}}
\put(85,11){\circle{2}}
\put(40,10){\line(0,-1){13}}
\put(40,-4){\circle{2}}
\end{picture}.
\noindent Using (\ref{wtspaces}), the representation
\begin{equation}\tag{\theequation$c$}
\tau_2:\,\, \gl \text{ on } \gu_2 \text{ is }
\setlength{\unitlength}{.75 mm}
\begin{picture}(95,18)
\put(3,10){$\times$}
\put(4,6){$1$}
\put(6,11){\line(1,0){13}}
\put(20,11){\circle{2}}
\put(21,11){\line(1,0){13}}
\put(35,11){\circle{2}}
\put(36,11){\line(1,0){13}}
\put(50,11){\circle{2}}
\put(51,11){\line(1,0){13}}
\put(65,11){\circle{2}}
\put(66,11){\line(1,0){13}}
\put(80,11){\circle{2}}
\put(35,10){\line(0,-1){13}}
\put(35,-4){\circle{2}}
\end{picture}
\end{equation}
Also by (\ref{wtspaces}), the representation $\tau_{-1}$ of $\gl$
on $\gu_{-1}$ is
\setlength{\unitlength}{.5 mm}
\begin{picture}(95,18)
\put(3,9){$\times$}
\put(2,4){$-2$}
\put(6,11){\line(1,0){13}}
\put(20,11){\circle{2}}
\put(19,4){$1$}
\put(21,11){\line(1,0){13}}
\put(35,11){\circle{2}}
\put(36,11){\line(1,0){13}}
\put(50,11){\circle{2}}
\put(51,11){\line(1,0){13}}
\put(65,11){\circle{2}}
\put(66,11){\line(1,0){13}}
\put(80,11){\circle{2}}
\put(35,10){\line(0,-1){13}}
\put(35,-4){\circle{2}}
\end{picture}
so the representation
\begin{equation}\tag{\theequation$d$}
\tau_1:\,\, \gl \text{ on } \gu_1 \text{ is }
\setlength{\unitlength}{.75 mm}
\begin{picture}(95,18)
\put(3,10){$\times$}
\put(2,6){$-1$}
\put(6,11){\line(1,0){13}}
\put(20,11){\circle{2}}
\put(19,6){$1$}
\put(21,11){\line(1,0){13}}
\put(35,11){\circle{2}}
\put(36,11){\line(1,0){13}}
\put(50,11){\circle{2}}
\put(51,11){\line(1,0){13}}
\put(65,11){\circle{2}}
\put(66,11){\line(1,0){13}}
\put(80,11){\circle{2}}
\put(35,10){\line(0,-1){13}}
\put(35,-4){\circle{2}}
\end{picture}
\end{equation}.
}\end{noname}
\noindent Note that $\tau_1|_{[\gl,\gl]}$ has degree $32$, is self--dual, 
and has an antisymmetric bilinear invariant.  Also, $\dim \gu_2 = 1$ and
$\tau_2|_{[\gl,\gl]}$ is trivial, so that bilinear invariant is given by the
Lie algebra product $\gu_1 \times \gu_1 \to \gu_2$.  In brief,
$\dim \gu_2 = 1$ and $\dim \gu_1 = 32$.
\medskip

\begin{noname}\label{e7-2} {\sc Case $E_{7,A_1D_6,2}$}.  
{\rm Again $G_0$ is the group of type $E_7$ with maximal compact subgroup 
that is a $2$--sheeted cover of $SU(2)\times Spin(12)$, but now noncompact 
simple root $\nu = \psi_6$, so $L$ is of type $T_1A_1D_5$ and the 
extended Dynkin diagram is
\begin{equation}\tag{\theequation$a$}
\setlength{\unitlength}{.75 mm}
\begin{picture}(140,18)
\put(5,16){\circle{2}}
\put(2,12){$-\mu$}
\put(6,16){.}
\put(7,16){.}
\put(8,16){.}
\put(9,16){.}
\put(10,16){.}
\put(11,16){.}
\put(12,16){.}
\put(13,16){.}
\put(14,16){.}
\put(15,16){.}
\put(16,16){.}
\put(17,16){.}
\put(18,16){.}
\put(20,16){\circle{2}}
\put(18,12){$\psi_1$}
\put(21,16){\line(1,0){13}}
\put(35,16){\circle{2}}
\put(33,12){$\psi_3$}
\put(36,16){\line(1,0){13}}
\put(50,16){\circle{2}}
\put(52,12){$\psi_4$}
\put(51,16){\line(1,0){13}}
\put(65,16){\circle{2}}
\put(63,12){$\psi_5$}
\put(66,16){\line(1,0){13}}
\put(80,16){\circle{2}}
\put(78,12){$\psi_6$}
\put(79,18){$\nu$}
\put(81,16){\line(1,0){13}}
\put(95,16){\circle{2}}
\put(93,12){$\psi_7$}
\put(50,15){\line(0,-1){13}}
\put(50,1){\circle{2}}
\put(52,1){$\psi_2$}
\put(110,10){(Type $E_7)$}
\end{picture}
\end{equation}
Thus $\gk$ is
\setlength{\unitlength}{.5 mm}
\begin{picture}(100,18)
\put(5,11){\circle{2}}
\put(2,6){$-\mu$}
\put(6,11){\line(1,0){13}}
\put(20,11){\circle{2}}
\put(18,5){$\psi_1$}
\put(21,11){\line(1,0){13}}
\put(35,11){\circle{2}}
\put(33,5){$\psi_3$}
\put(36,11){\line(1,0){13}}
\put(50,11){\circle{2}}
\put(52,5){$\psi_4$}
\put(51,11){\line(1,0){13}}
\put(65,11){\circle{2}}
\put(63,5){$\psi_5$}
\put(85,11){\circle{2}}
\put(83,5){$\psi_7$}
\put(50,10){\line(0,-1){13}}
\put(50,-4){\circle{2}}
\put(52,-6){$\psi_2$}
\end{picture}
and $\gl$ is
\setlength{\unitlength}{.5 mm}
\begin{picture}(95,18)
\put(5,11){\circle{2}}
\put(3,5){$\psi_1$}
\put(6,11){\line(1,0){13}}
\put(20,11){\circle{2}}
\put(18,5){$\psi_3$}
\put(21,11){\line(1,0){13}}
\put(35,11){\circle{2}}
\put(37,5){$\psi_4$}
\put(36,11){\line(1,0){13}}
\put(50,11){\circle{2}}
\put(48,5){$\psi_5$}
\put(51,11){\line(1,0){13}}
\put(62.5,9){$\times$}
\put(66,11){\line(1,0){13}}
\put(80,11){\circle{2}}
\put(78,5){$\psi_7$}
\put(35,10){\line(0,-1){13}}
\put(35,-4){\circle{2}}
\put(37,-6){$\psi_2$}
\end{picture}.
\medskip

\noindent
Now the decompositions $\gk_0 = \gk_1 \oplus \gk_2$ and 
$\gl_0 = \gl_1 \oplus \gl_2$ are
\begin{equation}\tag{\theequation$b$}
\gk_0 = \gso(12) \oplus \gsp(1) \text{ and } 
\gl_0 = \left ( \gso(2) \oplus \gso(10) \right ) \oplus \gsp(1) =
 \left ( i\R\nu^* \oplus \gso(10) \right ) \oplus \gsp(1).
\end{equation}
The representation of $\gk$ on $\gs$ has highest weight $-\nu = -\psi_6$:
\setlength{\unitlength}{.4 mm}
\begin{picture}(100,25)
\put(5,11){\circle{2}}
\put(6,11){\line(1,0){13}}
\put(20,11){\circle{2}}
\put(21,11){\line(1,0){13}}
\put(35,11){\circle{2}}
\put(36,11){\line(1,0){13}}
\put(50,11){\circle{2}}
\put(51,11){\line(1,0){13}}
\put(65,11){\circle{2}}
\put(64,3){$1$}
\put(85,11){\circle{2}}
\put(84,3){$1$}
\put(50,10){\line(0,-1){13}}
\put(50,-4){\circle{2}}
\end{picture}.
\noindent Using (\ref{wtspaces}), the representation
\begin{equation}\tag{\theequation$c$}
\tau_2:\,\, \gl \text{ on } \gu_2 \text{ is }
\setlength{\unitlength}{.75 mm}
\begin{picture}(95,18)
\put(5,11){\circle{2}}
\put(4,6){$1$}
\put(6,11){\line(1,0){13}}
\put(20,11){\circle{2}}
\put(21,11){\line(1,0){13}}
\put(35,11){\circle{2}}
\put(36,11){\line(1,0){13}}
\put(50,11){\circle{2}}
\put(51,11){\line(1,0){13}}
\put(63,9.5){$\times$}
\put(66,11){\line(1,0){13}}
\put(80,11){\circle{2}}
\put(35,10){\line(0,-1){13}}
\put(35,-4){\circle{2}}
\end{picture}.
\end{equation}
Also by (\ref{wtspaces}), the representation $\tau_{-1}$ of $\gl$
on $\gu_{-1}$ is
\setlength{\unitlength}{.5 mm}
\begin{picture}(95,18)
\put(5,11){\circle{2}}
\put(6,11){\line(1,0){13}}
\put(20,11){\circle{2}}
\put(21,11){\line(1,0){13}}
\put(35,11){\circle{2}}
\put(36,11){\line(1,0){13}}
\put(50,11){\circle{2}}
\put(49,4.5){$1$}
\put(51,11){\line(1,0){13}}
\put(62.5,9.5){$\times$}
\put(62,4.5){$-2$}
\put(66,11){\line(1,0){13}}
\put(80,11){\circle{2}}
\put(79,4.5){$1$}
\put(35,10){\line(0,-1){13}}
\put(35,-4){\circle{2}}
\end{picture}
so the representation
\begin{equation}\tag{\theequation$d$}
\tau_1:\,\, \gl \text{ on } \gu_1 \text{ is }
\setlength{\unitlength}{.75 mm}
\begin{picture}(95,18)
\put(5,11){\circle{2}}
\put(6,11){\line(1,0){13}}
\put(20,11){\circle{2}}
\put(21,11){\line(1,0){13}}
\put(35,11){\circle{2}}
\put(36,11){\line(1,0){13}}
\put(50,11){\circle{2}}
\put(51,11){\line(1,0){13}}
\put(63.2,10){$\times$}
\put(62,6){$-1$}
\put(66,11){\line(1,0){13}}
\put(80,11){\circle{2}}
\put(79,6){$1$}
\put(35,10){\line(0,-1){13}}
\put(35,-4){\circle{2}}
\put(37,-5){$1$}
\end{picture}.
\end{equation}
} \end{noname}
\noindent Note that $\tau_1|_{[\gl,\gl]}$ has degree $16$ and is not self--dual.
In brief, $\dim \gu_2 = 10$ and $\dim \gu_1 = 16$.
\medskip

\begin{noname}\label{e7-3} {\sc Case $E_{7,A_7}$}.  
{\rm Here $G_0$ is the group of type $E_7$ with maximal 
compact subgroup $SU(8)/\{\pm 1\}$.  The noncompact simple root 
$\nu = \psi_2$, so $L$ is of type $T_1E_6$ and the extended Dynkin
diagram is
\begin{equation}\tag{\theequation$a$}
\setlength{\unitlength}{.75 mm}
\begin{picture}(140,18)
\put(5,16){\circle{2}}
\put(2,12){$-\mu$}
\put(6,16){.}
\put(7,16){.}
\put(8,16){.}
\put(9,16){.}
\put(10,16){.}
\put(11,16){.}
\put(12,16){.}
\put(13,16){.}
\put(14,16){.}
\put(15,16){.}
\put(16,16){.}
\put(17,16){.}
\put(18,16){.}
\put(20,16){\circle{2}}
\put(18,12){$\psi_1$}
\put(21,16){\line(1,0){13}}
\put(35,16){\circle{2}}
\put(33,12){$\psi_3$}
\put(36,16){\line(1,0){13}}
\put(50,16){\circle{2}}
\put(52,12){$\psi_4$}
\put(51,16){\line(1,0){13}}
\put(65,16){\circle{2}}
\put(63,12){$\psi_5$}
\put(66,16){\line(1,0){13}}
\put(80,16){\circle{2}}
\put(78,12){$\psi_6$}
\put(81,16){\line(1,0){13}}
\put(95,16){\circle{2}}
\put(93,12){$\psi_7$}
\put(50,15){\line(0,-1){13}}
\put(50,1){\circle{2}}
\put(45,0){$\nu$}
\put(52,0){$\psi_2$}
\put(110,10){(Type $E_7$)}
\end{picture}
\end{equation}
Thus $\gk$ is
\setlength{\unitlength}{.5 mm}
\begin{picture}(105,18)
\put(5,6){\circle{2}}
\put(2,0){$-\mu$}
\put(6,6){\line(1,0){13}}
\put(20,6){\circle{2}}
\put(18,0){$\psi_1$}
\put(21,6){\line(1,0){13}}
\put(35,6){\circle{2}}
\put(33,0){$\psi_3$}
\put(36,6){\line(1,0){13}}
\put(50,6){\circle{2}}
\put(48,0){$\psi_4$}
\put(51,6){\line(1,0){13}}
\put(65,6){\circle{2}}
\put(63,0){$\psi_5$}
\put(66,6){\line(1,0){13}}
\put(80,6){\circle{2}}
\put(78,0){$\psi_6$}
\put(81,6){\line(1,0){13}}
\put(95,6){\circle{2}}
\put(93,0){$\psi_7$}
\end{picture}
and $\gl$ is
\setlength{\unitlength}{.5 mm}
\begin{picture}(95,18)
\put(5,11){\circle{2}}
\put(3,5){$\psi_1$}
\put(6,11){\line(1,0){13}}
\put(20,11){\circle{2}}
\put(18,5){$\psi_3$}
\put(21,11){\line(1,0){13}}
\put(35,11){\circle{2}}
\put(37,5){$\psi_4$}
\put(36,11){\line(1,0){13}}
\put(50,11){\circle{2}}
\put(48,5){$\psi_5$}
\put(51,11){\line(1,0){13}}
\put(65,11){\circle{2}}
\put(63,5){$\psi_6$}
\put(66,11){\line(1,0){13}}
\put(80,11){\circle{2}}
\put(78,5){$\psi_7$}
\put(35,10){\line(0,-1){13}}
\put(33,-6){$\times$}
\end{picture}.
\medskip

\noindent
Now the decompositions $\gk_0 = \gk_1 \oplus \gk_2$ and 
$\gl_0 = \gl_1 \oplus \gl_2$ are
\begin{equation}\tag{\theequation$b$}
\gk_0 = \gsu(8) \text{ and } 
\gl_0 = \left ( \gu(1) \oplus \gu(7) \right ) \cap  \gsu(8)
= i\R\nu^* \oplus \gsu(7).
\end{equation}
The representation of $\gk$ on $\gs$ has highest weight $-\nu = -\psi_2$:
\setlength{\unitlength}{.5 mm}
\begin{picture}(105,18)
\put(5,6){\circle{2}}
\put(6,6){\line(1,0){13}}
\put(20,6){\circle{2}}
\put(21,6){\line(1,0){13}}
\put(35,6){\circle{2}}
\put(36,6){\line(1,0){13}}
\put(50,6){\circle{2}}
\put(49,0){$1$}
\put(51,6){\line(1,0){13}}
\put(65,6){\circle{2}}
\put(66,6){\line(1,0){13}}
\put(80,6){\circle{2}}
\put(81,6){\line(1,0){13}}
\put(95,6){\circle{2}}
\end{picture}.
\noindent Using (\ref{wtspaces}), the representation
\begin{equation}\tag{\theequation$c$}
\tau_2:\,\, \gl \text{ on } \gu_2 \text{ is }
\setlength{\unitlength}{.75 mm}
\begin{picture}(95,18)
\put(5,11){\circle{2}}
\put(4,6){$1$}
\put(6,11){\line(1,0){13}}
\put(20,11){\circle{2}}
\put(21,11){\line(1,0){13}}
\put(35,11){\circle{2}}
\put(36,11){\line(1,0){13}}
\put(50,11){\circle{2}}
\put(51,11){\line(1,0){13}}
\put(65,11){\circle{2}}
\put(66,11){\line(1,0){13}}
\put(80,11){\circle{2}}
\put(35,10){\line(0,-1){13}}
\put(33.5,-5){$\times$}
\end{picture}.
\end{equation}
Also by (\ref{wtspaces}), the representation $\tau_{-1}$ of $\gl$
on $\gu_{-1}$ is
\setlength{\unitlength}{.5 mm}
\begin{picture}(95,18)
\put(5,11){\circle{2}}
\put(6,11){\line(1,0){13}}
\put(20,11){\circle{2}}
\put(21,11){\line(1,0){13}}
\put(35,11){\circle{2}}
\put(36,5){$1$}
\put(36,11){\line(1,0){13}}
\put(50,11){\circle{2}}
\put(51,11){\line(1,0){13}}
\put(65,11){\circle{2}}
\put(66,11){\line(1,0){13}}
\put(80,11){\circle{2}}
\put(35,10){\line(0,-1){13}}
\put(32.5,-6){$\times$}
\put(38,-5){$-2$}
\end{picture}
so the representation
\begin{equation}\tag{\theequation$d$}
\tau_1:\,\, \gl \text{ on } \gu_1 \text{ is }
\setlength{\unitlength}{.75 mm}
\begin{picture}(95,18)
\put(5,11){\circle{2}}
\put(6,11){\line(1,0){13}}
\put(20,11){\circle{2}}
\put(21,11){\line(1,0){13}}
\put(35,11){\circle{2}}
\put(36,11){\line(1,0){13}}
\put(50,11){\circle{2}}
\put(49,6){$1$}
\put(51,11){\line(1,0){13}}
\put(65,11){\circle{2}}
\put(66,11){\line(1,0){13}}
\put(80,11){\circle{2}}
\put(35,10){\line(0,-1){13}}
\put(33.5,-5){$\times$}
\put(37,-5){$-1$}
\end{picture}.
\end{equation}
} \end{noname}
\medskip

\noindent Note that $\tau_1|_{[\gl,\gl]}$ has degree $35$ and is not self--dual.
In brief, $\dim \gu_2 = 7$ and $\dim \gu_1 = 35$
\medskip

\begin{noname} \label{e8-1} {\sc Case $E_{8,D_8}$}.  
{\rm Here $G_0$ is the group of type $E_8$ with maximal 
compact subgroup locally isomorphic to $Spin(16)$.  The noncompact simple root 
$\nu = \psi_1$, so $L$ is of type $T_1D_7$ and the extended Dynkin
diagram is
\begin{equation}\tag{\theequation$a$}
\setlength{\unitlength}{.75 mm}
\begin{picture}(155,18)
\put(5,16){\circle{2}}
\put(1,16){$\nu$}
\put(6,12){$\psi_1$}
\put(6,16){\line(1,0){13}}
\put(20,16){\circle{2}}
\put(18,12){$\psi_3$}
\put(21,16){\line(1,0){13}}
\put(35,16){\circle{2}}
\put(37,12){$\psi_4$}
\put(36,16){\line(1,0){13}}
\put(50,16){\circle{2}}
\put(48,12){$\psi_5$}
\put(51,16){\line(1,0){13}}
\put(65,16){\circle{2}}
\put(63,12){$\psi_6$}
\put(66,16){\line(1,0){13}}
\put(80,16){\circle{2}}
\put(78,12){$\psi_7$}
\put(81,16){\line(1,0){13}}
\put(95,16){\circle{2}}
\put(93,12){$\psi_8$}
\put(96,16){.}
\put(97,16){.}
\put(98,16){.}
\put(99,16){.}
\put(100,16){.}
\put(101,16){.}
\put(102,16){.}
\put(103,16){.}
\put(104,16){.}
\put(105,16){.}
\put(106,16){.}
\put(107,16){.}
\put(108,16){.}
\put(110,16){\circle{2}}
\put(108,12){$-\mu$}
\put(35,15){\line(0,-1){13}}
\put(35,1){\circle{2}}
\put(37,0){$\psi_2$}
\put(125,10){(Type $E_8$)}
\end{picture}
\end{equation}
Thus $\gk$ is
\setlength{\unitlength}{.5 mm}
\begin{picture}(115,18)
\put(10,11){\circle{2}}
\put(8,5){$\psi_3$}
\put(11,11){\line(1,0){13}}
\put(25,11){\circle{2}}
\put(27,5){$\psi_4$}
\put(26,11){\line(1,0){13}}
\put(40,11){\circle{2}}
\put(38,5){$\psi_5$}
\put(41,11){\line(1,0){13}}
\put(55,11){\circle{2}}
\put(53,5){$\psi_6$}
\put(56,11){\line(1,0){13}}
\put(70,11){\circle{2}}
\put(68,5){$\psi_7$}
\put(71,11){\line(1,0){13}}
\put(85,11){\circle{2}}
\put(83,5){$\psi_8$}
\put(86,11){\line(1,0){13}}
\put(100,11){\circle{2}}
\put(98,5){$-\mu$}
\put(25,10){\line(0,-1){13}}
\put(25,-4){\circle{2}}
\put(27,-5){$\psi_2$}
\end{picture}
and $\gl$ is
\setlength{\unitlength}{.5 mm}
\begin{picture}(105,18)
\put(3,9){$\times$}
\put(6,11){\line(1,0){13}}
\put(20,11){\circle{2}}
\put(18,6){$\psi_3$}
\put(21,11){\line(1,0){13}}
\put(35,11){\circle{2}}
\put(37,6){$\psi_4$}
\put(36,11){\line(1,0){13}}
\put(50,11){\circle{2}}
\put(48,6){$\psi_5$}
\put(51,11){\line(1,0){13}}
\put(65,11){\circle{2}}
\put(63,6){$\psi_6$}
\put(66,11){\line(1,0){13}}
\put(80,11){\circle{2}}
\put(78,6){$\psi_7$}
\put(81,11){\line(1,0){13}}
\put(95,11){\circle{2}}
\put(93,6){$\psi_8$}
\put(35,10){\line(0,-1){13}}
\put(35,-4){\circle{2}}
\put(37,-5){$\psi_2$}
\end{picture}.
\medskip

\noindent
Now the decompositions $\gk_0 = \gk_1 \oplus \gk_2$ and 
$\gl_0 = \gl_1 \oplus \gl_2$ are`
\begin{equation}\tag{\theequation$b$}
\gk_0 = \gso(16) \text{ and } 
\gl_0 =  \gso(2)\oplus \gso(14) = i\R\nu^* \oplus \gso(14).
\end{equation}
The representation of $\gk$ on $\gs$ has highest weight $-\nu = -\psi_1$:
\setlength{\unitlength}{.5 mm}
\begin{picture}(105,15)
\put(10,6){\circle{2}}
\put(9,-1){$1$}
\put(11,6){\line(1,0){13}}
\put(25,6){\circle{2}}
\put(26,6){\line(1,0){13}}
\put(40,6){\circle{2}}
\put(41,6){\line(1,0){13}}
\put(55,6){\circle{2}}
\put(56,6){\line(1,0){13}}
\put(70,6){\circle{2}}
\put(71,6){\line(1,0){13}}
\put(85,6){\circle{2}}
\put(86,6){\line(1,0){13}}
\put(100,6){\circle{2}}
\put(25,5){\line(0,-1){13}}
\put(25,-9){\circle{2}}
\end{picture}.
\noindent Using (\ref{wtspaces}), the representation
\begin{equation}\tag{\theequation$c$}
\tau_2:\,\, \gl \text{ on } \gu_2 \text{ is }
\setlength{\unitlength}{.75 mm}
\begin{picture}(115,18)
\put(3,10){$\times$}
\put(6,11){\line(1,0){13}}
\put(20,11){\circle{2}}
\put(21,11){\line(1,0){13}}
\put(35,11){\circle{2}}
\put(36,11){\line(1,0){13}}
\put(50,11){\circle{2}}
\put(51,11){\line(1,0){13}}
\put(65,11){\circle{2}}
\put(66,11){\line(1,0){13}}
\put(80,11){\circle{2}}
\put(81,11){\line(1,0){13}}
\put(95,11){\circle{2}}
\put(94,5){$1$}
\put(35,10){\line(0,-1){13}}
\put(35,-4){\circle{2}}
\end{picture}.
\end{equation}
Also by (\ref{wtspaces}), the representation $\tau_{-1}$ of $\gl$
on $\gu_{-1}$ is
\setlength{\unitlength}{.5 mm}
\begin{picture}(110,18)
\put(3,9){$\times$}
\put(2,4){$-2$}
\put(6,11){\line(1,0){13}}
\put(20,11){\circle{2}}
\put(19,4){$1$}
\put(21,11){\line(1,0){13}}
\put(35,11){\circle{2}}
\put(36,11){\line(1,0){13}}
\put(50,11){\circle{2}}
\put(51,11){\line(1,0){13}}
\put(65,11){\circle{2}}
\put(66,11){\line(1,0){13}}
\put(80,11){\circle{2}}
\put(81,11){\line(1,0){13}}
\put(95,11){\circle{2}}
\put(35,-4){\circle{2}}
\put(35,10){\line(0,-1){13}}
\end{picture}
so the representation
\begin{equation}\tag{\theequation$d$}
\tau_1:\,\, \gl \text{ on } \gu_1 \text{ is }
\setlength{\unitlength}{.75 mm}
\begin{picture}(110,18)
\put(3,10){$\times$}
\put(2,6){$-2$}
\put(6,11){\line(1,0){13}}
\put(20,11){\circle{2}}
\put(21,11){\line(1,0){13}}
\put(35,11){\circle{2}}
\put(36,11){\line(1,0){13}}
\put(50,11){\circle{2}}
\put(51,11){\line(1,0){13}}
\put(65,11){\circle{2}}
\put(66,11){\line(1,0){13}}
\put(80,11){\circle{2}}
\put(81,11){\line(1,0){13}}
\put(95,11){\circle{2}}
\put(35,-4){\circle{2}}
\put(36,-3){$1$}
\put(35,10){\line(0,-1){13}}
\end{picture}.
\end{equation}
} \end{noname}
\noindent Note that $\tau_1|_{[\gl,\gl]}$ has degree $64$ and is not 
self--dual; $\dim \gu_2 = 14$ and $\dim \gu_1 = 64$.
\medskip

\begin{noname} \label{e8-2} {\sc Case $E_{8,A_1E_7}$}.  
{\em Here $G_0$ is the group of type $E_8$ with maximal compact subgroup 
that has $SU(2)\times E_7$ as a double cover..  The noncompact simple root 
$\nu = \psi_8$, so $L$ is of type $T_1E_7$ and the extended Dynkin
diagram is
\begin{equation}\tag{\theequation$a$}
\setlength{\unitlength}{.75 mm}
\begin{picture}(155,18)
\put(5,16){\circle{2}}
\put(6,12){$\psi_1$}
\put(6,16){\line(1,0){13}}
\put(20,16){\circle{2}}
\put(18,12){$\psi_3$}
\put(21,16){\line(1,0){13}}
\put(35,16){\circle{2}}
\put(37,12){$\psi_4$}
\put(36,16){\line(1,0){13}}
\put(50,16){\circle{2}}
\put(48,12){$\psi_5$}
\put(51,16){\line(1,0){13}}
\put(65,16){\circle{2}}
\put(63,12){$\psi_6$}
\put(66,16){\line(1,0){13}}
\put(80,16){\circle{2}}
\put(78,12){$\psi_7$}
\put(81,16){\line(1,0){13}}
\put(95,16){\circle{2}}
\put(93,12){$\psi_8$}
\put(93,18){$\nu$}
\put(96,16){.}
\put(97,16){.}
\put(98,16){.}
\put(99,16){.}
\put(100,16){.}
\put(101,16){.}
\put(102,16){.}
\put(103,16){.}
\put(104,16){.}
\put(105,16){.}
\put(106,16){.}
\put(107,16){.}
\put(108,16){.}
\put(110,16){\circle{2}}
\put(108,12){$-\mu$}
\put(35,15){\line(0,-1){13}}
\put(35,1){\circle{2}}
\put(37,0){$\psi_2$}
\put(125,10){(Type $E_8$)}
\end{picture}
\end{equation}
Thus $\gk$ is
\setlength{\unitlength}{.5 mm}
\begin{picture}(115,18)
\put(5,11){\circle{2}}
\put(3,5){$\psi_1$}
\put(6,11){\line(1,0){13}}
\put(20,11){\circle{2}}
\put(18,5){$\psi_3$}
\put(21,11){\line(1,0){13}}
\put(35,11){\circle{2}}
\put(37,5){$\psi_4$}
\put(36,11){\line(1,0){13}}
\put(50,11){\circle{2}}
\put(48,5){$\psi_5$}
\put(51,11){\line(1,0){13}}
\put(65,11){\circle{2}}
\put(63,5){$\psi_6$}
\put(66,11){\line(1,0){13}}
\put(80,11){\circle{2}}
\put(78,5){$\psi_7$}
\put(100,11){\circle{2}}
\put(98,5){$-\mu$}
\put(35,10){\line(0,-1){13}}
\put(35,-4){\circle{2}}
\put(37,-5){$\psi_2$}
\end{picture}
and $\gl$ is
\setlength{\unitlength}{.5 mm}
\begin{picture}(105,18)
\put(5,11){\circle{2}}
\put(3,5){$\psi_1$}
\put(6,11){\line(1,0){13}}
\put(20,11){\circle{2}}
\put(18,5){$\psi_3$}
\put(21,11){\line(1,0){13}}
\put(35,11){\circle{2}}
\put(37,5){$\psi_4$}
\put(36,11){\line(1,0){13}}
\put(50,11){\circle{2}}
\put(48,5){$\psi_5$}
\put(51,11){\line(1,0){13}}
\put(65,11){\circle{2}}
\put(63,5){$\psi_6$}
\put(66,11){\line(1,0){13}}
\put(80,11){\circle{2}}
\put(78,5){$\psi_7$}
\put(81,11){\line(1,0){13}}
\put(93,9){$\times$}
\put(35,10){\line(0,-1){13}}
\put(35,-4){\circle{2}}
\put(37,-5){$\psi_2$}
\end{picture}.
\medskip

\noindent
Now the decompositions $\gk_0 = \gk_1 \oplus \gk_2$ and 
$\gl_0 = \gl_1 \oplus \gl_2$ are
\begin{equation}\tag{\theequation$b$}
\gk_0 = \gsp(1) \oplus \ge_7 \text{ and } \gl_0 = i\R\nu^* \oplus \ge7.
\end{equation}
The representation of $\gk$ on $\gs$ has highest weight $-\nu = -\psi_8$:
\setlength{\unitlength}{.5 mm}
\begin{picture}(105,20)
\put(5,9){\circle{2}}
\put(6,9){\line(1,0){13}}
\put(20,9){\circle{2}}
\put(21,9){\line(1,0){13}}
\put(35,9){\circle{2}}
\put(36,9){\line(1,0){13}}
\put(50,9){\circle{2}}
\put(51,9){\line(1,0){13}}
\put(65,9){\circle{2}}
\put(66,9){\line(1,0){13}}
\put(80,9){\circle{2}}
\put(79,2){$1$}
\put(90,9){\circle{2}}
\put(89,2){$1$}
\put(35,8){\line(0,-1){13}}
\put(35,-6){\circle{2}}
\end{picture}.
\noindent Using (\ref{wtspaces}), the representation
\begin{equation}\tag{\theequation$c$}
\tau_2:\,\, \gl \text{ on } \gu_2 \text{ is }
\setlength{\unitlength}{.75 mm}
\begin{picture}(115,18)
\put(5,11){\circle{2}}
\put(6,11){\line(1,0){13}}
\put(20,11){\circle{2}}
\put(21,11){\line(1,0){13}}
\put(35,11){\circle{2}}
\put(36,11){\line(1,0){13}}
\put(50,11){\circle{2}}
\put(51,11){\line(1,0){13}}
\put(65,11){\circle{2}}
\put(66,11){\line(1,0){13}}
\put(80,11){\circle{2}}
\put(81,11){\line(1,0){13}}
\put(93,10){$\times$}
\put(94,5){$1$}
\put(35,10){\line(0,-1){13}}
\put(35,-4){\circle{2}}
\end{picture}.
\end{equation}
Also by (\ref{wtspaces}), the representation $\tau_{-1}$ of $\gl$
on $\gu_{-1}$ is
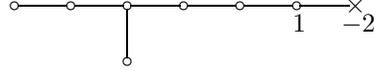
\setlength{\unitlength}{.5 mm}
\begin{picture}(110,18)
\put(5,11){\circle{2}}
\put(6,11){\line(1,0){13}}
\put(20,11){\circle{2}}
\put(21,11){\line(1,0){13}}
\put(35,11){\circle{2}}
\put(36,11){\line(1,0){13}}
\put(50,11){\circle{2}}
\put(51,11){\line(1,0){13}}
\put(65,11){\circle{2}}
\put(66,11){\line(1,0){13}}
\put(80,11){\circle{2}}
\put(79,4){$1$}
\put(81,11){\line(1,0){13}}
\put(93,9){$\times$}
\put(92,4){$-2$}
\put(35,-4){\circle{2}}
\put(35,10){\line(0,-1){13}}
\end{picture}
so the representation
\begin{equation}\tag{\theequation$d$}
\tau_1:\,\, \gl \text{ on } \gu_1 \text{ is }
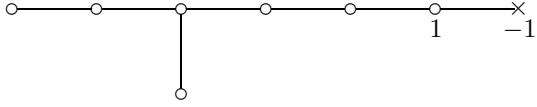
\setlength{\unitlength}{.75 mm}
\begin{picture}(110,18)
\put(5,11){\circle{2}}
\put(6,11){\line(1,0){13}}
\put(20,11){\circle{2}}
\put(21,11){\line(1,0){13}}
\put(35,11){\circle{2}}
\put(36,11){\line(1,0){13}}
\put(50,11){\circle{2}}
\put(51,11){\line(1,0){13}}
\put(65,11){\circle{2}}
\put(66,11){\line(1,0){13}}
\put(80,11){\circle{2}}
\put(79,6){$1$}
\put(81,11){\line(1,0){13}}
\put(93,10){$\times$}
\put(92,6){$-1$}
\put(35,-4){\circle{2}}
\put(35,10){\line(0,-1){13}}
\end{picture}.
\end{equation}
} \end{noname}
Note that $\tau_1|_{[\gl,\gl]}$, has degree $56$, is self-dual, and has an 
antisymmetric bilinear invariant.  Also, $\dim \gu_2 = 1$ and
$\tau_2|_{[\gl,\gl]}$ is trivial, so that bilinear invariant is given by the
Lie algebra product $\gu_1 \times \gu_1 \to \gu_2$.

\vskip 1 mm 
This completes our run through the exceptional cases.
\medskip

\section{Prehomogeneity and Relative Invariants for $(L,\gu_1)$}\label{sec4}
\setcounter{equation}{0}

Consider a connected linear algebraic group with a rational representation
on a complex vector space.  We say that the triple consisting of the 
group, the representation and the vector space is {\em prehomogeneous} if
there is a Zariski--dense orbit.  When no confusion is possible we omit the
representation.  A general theorem of Vinberg on graded Lie algebras 
(see \cite[Theorem 10.19]{Kn1}) shows that $(L,\gu_1)$ is
prehomogeneous.  Or one can verify that fact by running through the
lists of Section \ref{sec3} and the classification of \cite{SK}.  In
fact we will do the latter in order to describe the algebra of
relative--invariant polynomials on $\gu_1$ and the $L$--orbit structure
of $\gu_1$ for each instance of $(L,\gu_1)$.
We shall also use the notation $V$ for the fundamental space $\gu_1$
and $V^* = \gu_{-1}$ for its dual space.
\medskip

We recall some material on prehomogeneous spaces as it
applies to $(L,\gu_1)$.  There is no nonconstant
$L$--invariant rational function $f: \gu_1 \to \C$ because the $L$--invariance
would force it to be constant on the Zariski--dense $L$--orbit
\cite[Proposition 3 in \S 2]{SK}.  By {\em relative invariant} for 
$(L,\gu_1)$ we mean a nonconstant 
polynomial function $f: \gu_1 \to \C$ such that $f(\ell\xi) = \chi(\ell)f(\xi)$
for some rational character $\chi : L \to \C^\times$.  The quotient of two
relative invariants with the same character would be an
$L$--invariant rational function of $\gu_1$, hence constant, so a relative
invariant $f$ is determined up to scalar multiple by its character $\chi$
\cite[Proposition 3 in \S 4]{SK}.  In particular all the 
$f_c: \xi \mapsto f(c\xi)$ are proportional, so $f$ is a
{\em homogeneous} polynomial.  It will be convenient to denote
\begin{equation}\label{relinv}
\cA(L,\gu_1): \text{ the associative algebra of all relative invariants
of } (L,\gu_1).
\end{equation}
\medskip

The {\em regular set} for $(L,\gu_1)$ is the open $L$--orbit
$\cO_0 := \Ad(L)\xi_0 \subset \gu_1$ and the {\em singular set}
is its complement $\gu_1 \setminus \cO_0$.  Let $V_1, \dots , V_e$
be those components of the singular set that are of codimension $1$ 
in $\gu_1$.  For each $i$, $V_i$ is the zero set of an irreducible
polynomial $f_i$.  The algebra of relative invariants for $(L,\gu_1)$
is the polynomial algebra $\C[f_1, \dots , f_e]$
\cite[Proposition 5 in \S 4]{SK}.  In particular $(L,\gu_1)$ has a
relative invariant if and only if the its singular set has a component 
of codimension $1$.  So far we haven't used irreducibility of $L$ on
$\gu_1$, but now we use it to see \cite[Proposition 12 in \S 4]{SK}
that $e \leqq 1$, i.e. that either $\cA(L,\gu_1) = \C$ (in other words
$(L.\gu_1)$ has no relative invariant) or $\cA(L,\gu_1)$ has form $\C[f]$.
\medskip

\addtocounter{equation}{1} \noindent \underline{{\bf \theequation.} \phantom{x}
Cases $SO(2p,r)$}.
We first consider the various cases where $G_0$ is the universal covering
group of the indefinite orthogonal group $SO(2p,2q)$ or $SO(2p,2q+1)$.  For
convenience we write that as $SO(2p,r)$.  Then $L_0$ consists of all
$\left ( \begin{smallmatrix} a & 0 \\ 0 & b \end{smallmatrix}\right )$
where $a$ is in the image of the standard embedding 
$\iota: U(p) \hookrightarrow SO(2p)$ and where $b \in SO(r)$.  Here
$\gs_0 = 
\left \{ \left ( \left . \begin{smallmatrix} 0 & x \\ {}^tx & 0 
\end{smallmatrix}\right ) \right | x \in \R^{2p\times r}
\right \} \cong \R^{2p\times r}$ and the (conjugation) action 
$\left ( \begin{smallmatrix} a & 0 \\ 0 & b \end{smallmatrix}\right )\in L_0$
on $\gs_0$ is given by $x \mapsto axb^{-1}$.  Now $\gu_1   
\cong \C^{p\times r}$ with the action of $L \simeq GL(p;\C)\times SO(r;\C)$ 
given by $\ell = 
\left ( \begin{smallmatrix} \ell_1 & 0 \\ 0 & \ell_2 \end{smallmatrix}\right ):
z \mapsto \ell_1z\ell_2^{-1}$. Then $f(z) := \det(z\cdot {}^tz)$
transforms by $f(\ell(z)) = \det(\ell)^2f(z)$.  However it is a
relative invariant only when it is not identically zero, i.e. when
$p \leqq r$.
\medskip

On the other hand, if $p > r$ then the $(SL(p;\C)\times SO(r;\C))$--orbit
of $\left ( \begin{smallmatrix} I_r \\ 0_{p-r} \end{smallmatrix}\right )$  
is open in $\C^{p\times r}$, so there is no nonconstant 
$(SL(p;\C)\times SO(r;\C))$--invariant.  It follows that there is no 
relative invariant for $L$.
\medskip

Summary: if $p \leqq r$ then $\cA(L,\gu_1) = \C[f]$ where 
$f(z) := \det(z\cdot {}^tz)$, polynomial of degree $2p$.
If $p > r$ then $\cA(L,\gu_1) = \C$.
\medskip 

\addtocounter{equation}{1} \noindent \underline{{\bf \theequation.} \phantom{x}
Cases $Sp(p,q)$}.  We next consider the cases where $G_0 = Sp(p,q)$.  Then
$L_0$ consists of all
$\left ( \begin{smallmatrix} a & 0 \\ 0 & b \end{smallmatrix}\right )$
where $a$ is in the image of the standard embedding
$\iota: U(p) \hookrightarrow Sp(p)$ and where $b \in Sp(q)$. Here
$\gs_0 =
\left \{ \left ( \left . \begin{smallmatrix} 0 & x \\ x^* & 0
\end{smallmatrix}\right ) \right | x \in \H^{p\times q}
\right \} \cong \H^{p\times q}$ and the (conjugation) action
$\left ( \begin{smallmatrix} a & 0 \\ 0 & b \end{smallmatrix}\right )\in L_0$
on $\gs_0$ is given by $x \mapsto axb^{-1}$.  Now $\gu_1
\cong \C^{p\times 2q}$ with the action of $L \simeq GL(p;\C)\times Sp(q;\C)$
given by $\ell =
\left ( \begin{smallmatrix} \ell_1 & 0 \\ 0 & \ell_2 \end{smallmatrix}\right ):
z \mapsto \ell_1z\ell_2^{-1}$.   Let 
$J =  \left ( \begin{smallmatrix} 0 & I_q \\ -I_q & 0\end{smallmatrix}\right )$,
so $Sp(q;\C)$ is characterized by $b\cdot J\cdot {}^tb = J$.  Let $\Pf$ 
denote the Pfaffian polynomial on the space of antisymmetric 
$2q \times 2q$ matrices, so $\Pf(m)^2 = \det(m)$.  Then
$f(z) := \Pf(z\cdot J\cdot {}^tz)$ transforms by $f(\ell z) = \det(\ell)f(z)$.
And of course $f$ is a relative invariant only when it is not identically zero.
For that we must have the possibility that the $p\times p$ antisymmetric matrix
$z\cdot J\cdot {}^tz$ is nonsingular, which is the case just when 
both that $p \leqq 2q$ and $p$ is even.
\medskip

As before, if $p > 2q$ then the $(SL(p;\C)\times Sp(q;\C))$--orbit
of $\left ( \begin{smallmatrix} I_{2q} \\ 0_{p-2q} \end{smallmatrix}\right )$
is open in $\C^{p\times 2q}$, so there is no nonconstant
$(SL(p;\C)\times Sp(q;\C))$--invariant, and thus no relative invariant
for $L$.  
\medskip

If $p \leqq 2q$ but $p$ is odd we define $m \in \C^{p\times 2q}$ by
$m_{i,i} = 1$ if $i \leqq p$ and $i$ is odd,
$m_{i,q+i} = 1$ if $i \leqq p$ and $i$ is even, all other entries zero.  
The point is that the row space of $m$ has dimension $p$ and has
nullity $1$ relative to the bilinear form $J$ that
defines $Sp(q;\C)$.  Then $(SL(p;\C)\times Sp(q;\C))(m)$ consists of all 
elements of $\C^{p\times 2q}$ whose row space has dimension $p$ and nullity 
$1$ relative to $J$, and that is open in $\C^{p\times 2q}$.
As above, it follows that there is no relative invariant for $L$.
\medskip

Summary: if $p \leqq 2q$ and $p$ is even then $\cA(L,\gu_1) = \C[f]$ where
$f(z) := \Pf(z\cdot J\cdot {}^tz)$, polynomial of degree $p$.
If $p > 2q$ of if $p$ is odd then $\cA(L,\gu_1) = \C$.
\medskip

We have completely described $\cA(L\gu_1)$ when $G$ is classical.  Except
for a few extreme cases, the representation $\tau_1|_{L'}$ of the derived
group $L' = [L,L]$ failed to be self--dual because of a tensor factor
\setlength{\unitlength}{.5 mm}
\begin{picture}(70,10)
\put(14,5){\circle{2}}
\put(12,-2){$1$}
\put(15,5){\line(1,0){8}}
\put(24,5){\circle{2}}
\put(25,5){\line(1,0){8}}
\put(36,5){\circle*{1}}
\put(38,5){\circle*{1}}
\put(40,5){\circle*{1}}
\put(43,5){\line(1,0){8}}
\put(52,5){\circle{2}}
\put(53,5){\line(1,0){7}}
\put(61,5){\circle{2}}
\end{picture}.
In many of the exceptional group cases, $\tau_1|_{L'}$ is self--dual, so
it has a bilinear invariant, and when that bilinear invariant is symmetric 
it generates $\cA(L,\gu_1)$.  Also, if that bilinear invariant is
antisymmetric, then $\tau_1(L')$ is contained in the symplectic group $J$
of the bilinear invariant, and if $\sigma$ is the representation of $J$
on $\gu_1$ then $S^2(S^2(\sigma))$ contains the trivial representation with
multiplicity $1$, and that gives a quartic invariant that generates
$\cA(L,\gu_1)$.  Before formalizing these statements we look back at
Section \ref{sec3} to see the first five columns of
\begin{equation}\label{invariants-exceptionals}
\begin{tabular}{|l|l|c|c|c|c|}\hline
Case & $G_0$ & $\deg \tau_1$ & self--dual? & bilinear invariant & relative
invariant \\
\hline \hline
(\ref{g2}) & $G_{2,A_1A_1}$ & $4$ & yes & antisymmetric & degree $4$ \\
\hline
(\ref{f4-1}) & $F_{4,A_1C_3}$ & $14$ & yes & antisymmetric & degree $4$ \\
\hline
(\ref{f4-2}) & $F_{4,B_4}$ & $8$ & yes & symmetric & degree $2$\\
\hline
(\ref{e6-1}) & $E_{6,A_1A_5,1}$ & $20$ & no & none & none \\
\hline
(\ref{e6-2}) & $E_{6,A_1A_5,2}$ & $20$ & yes & antisymmetric & degree $4$ \\
\hline
(\ref{e7-1}) & $E_{7,A_1D_6,1}$ & $32$ & yes & antisymmetric & degree $4$ \\
\hline
(\ref{e7-2}) & $E_{7,A_1D_6,2}$ & $16$ & no & none & none \\
\hline
(\ref{e7-3}) & $E_{7,A_7}$ & $35$ & no & none &  degree $7$ \\
\hline
(\ref{e8-1}) & $E_{8,D_8}$ & $64$ & no & none &  degree $8$ \\
\hline
(\ref{e8-2}) & $E_{8,A_1E_7}$ & $56$ & yes & antisymmetric & degree $4$  \\
\hline
\end{tabular}\phantom{X}.
\end{equation}

The information of the last column of Table \ref{invariants-exceptionals}
is contained in \cite{SK}, but we can give a short direct proof of the
cases where there is a relative invariant, as follows.

\begin{lemma}\label{bilinv}
If $\tau_1|_{L'}$ is self--dual there are two possibilities.  Either
it has a nonzero symmetric bilinear invariant $b$ and $\cA(L,\gu_1)
= \C[b]$, or it has a nonzero antisymmetric bilinear invariant and
$\cA(L,\gu_1) = \C[f]$ where $f$ has degree $4$.  In the non self--dual
case {\rm (\ref{e7-3})} we have $\cA(L,\gu_1) = \C[f]$ where $f$ has degree $7$,
and in the non self--dual case {\rm (\ref{e8-1})} we have $\cA(L,\gu_1) = \C[f]$
where $f$ has degree $8$. 
\end{lemma}

\noindent {\bf Proof.}  If the bilinear invariant $b$ is symmetric,
then since it has degree $2$ it must generate $\cA(L,\gu_1)$.  If
$b$ is antisymmetric, then in each of the five relevant cases of Table
\ref{invariants-exceptionals} we compute symmetric powers $S^2(\tau_1|_{L'})$,
$S^3(\tau_1|_{L'})$ and $S^4(\tau_1|_{L'})$ to see that we first 
encounter a $\tau_1(L')$--invariant in degree $4$.  (This degree $4$
semiinvariant can also be seen by a classification free argument 
\cite[Proposition 1.4]{P}.)
\medskip

Consider the two non self--dual cases of Table \ref{invariants-exceptionals} 
for which we claim a $\tau_1|_{L'}$--invariant.
In case (\ref{e7-3}) we compute the
$S^r(\tau_1|_{L'})$ for $2 \leqq r \leqq 7$ to see that we first
encounter a $\tau_1(L')$--invariant in degree $7$, and
and in case (\ref{e8-1}) we compute the
$S^r(\tau_1|_{L'})$ for $2 \leqq r \leqq 8$ to see that we first
encounter a $\tau_1(L')$--invariant in degree $8$.
\hfill $\square$
\medskip

\section{Negativity and $K_0$--types}\label{sec5}

In this section we discuss negativity of a homogeneous holomorphic
vector bundle over $G_0/L_0$ and the $K_0$--types of the resulting
discrete series representations.
\medskip

Recall some notation from Section \ref{sec2}.  The flag domain
$D = G_0(z_0) \cong G_0/L_0$ is an open $G_0$--orbit in the complex
flag manifold $Z = G/Q$, where $z_0 = 1Q$ is the base point and
$L_0 = G_0\cap Q$.  The parabolic subgroup $Q$ of $G$ has Lie algebra
$\gq = \gl + \gu_-$ and its nilradical $\gu_-$ is opposite to
$\gu_+$, which in turn represents the holomorphic tangent space to
$D$ at $z_0$.  According the the multiplicity of the noncompact simple
root, $\gu_+ = \gu_1 + \gu_2$.  The maximal compact subvariety 
$Y = K_0(z_0) = K(z_0)$ has holomorphic tangent space at $z_0$ represented
by $\gu_2$ and has holomorphic normal space represented by $V = \gu_1$.  The
group $L$ acts irreducibly on both of them, and those representations were
derived explicitly in Section \ref{sec3}.  The variety $Y$ is a complex
flag manifold $K/(K\cap Q)$ in its own right, and is the fiber of
the basic fibration (\ref{fibration}) $D = G_0/L_0 \to G_0/K_0$.
\medskip

Fix an irreducible representation $\tau_\gamma$ of $L$.  Here $\gamma$
is the highest weight, $E_\gamma$ is the representation space, 
$\E_\gamma \to D$ is the associated homogeneous holomorphic vector
bundle. and $\cO(\E_\gamma) \to D$ is the sheaf of germs of holomorphic
sections.  
\medskip

By $\cO(\E_\gamma)|_Y \to D$ we mean the pull--back sheaf of 
$\cO(\E_\gamma) \to D$ under $Y \hookrightarrow D$.
It is a sheaf on $D$ supported on $Y$.
We filter it by order of vanishing in directions transverse to $Y$:
\begin{equation} \label{filtration}
\cF^n(\E_\gamma) = \bigl \{f \in \cO(\E_\gamma)|_Y \,\,\bigl | \,\, f
\text{ vanishes to order } \geqq n \text{ in directions transverse to } 
Y \bigr \}\bigr . .
\end{equation}
We also need the notation
\begin{equation}
\begin{aligned}
&\N_Y \to Y: \text{ holomorphic normal bundle to } Y \text{ in } D, \\
&\N_Y^* \to Y: \text{ holomorphic conormal bundle to } Y \text{ in } 
	D \text{ and } \\
&\cS^n(\N_Y^*) = \cO(S^n(\N_Y^*)) \text{ where } S^n(\N_Y^*) \to Y
\text{ is the } n^{th}\text{ symmetric power of } \N_Y^* \to Y.
\end{aligned}
\end{equation}
Then $\N_Y \to Y$ is the homogeneous holomorphic vector bundle over $Y$ with 
fiber represented by $V = \gu_1$, its dual $\N_Y^* \to Y$ is the homogeneous 
holomorphic vector bundle with fiber $V^* = \gu_{-1}$, similarly
for the third bundle with fibers $S^n(V^*)$, and we view $\cS^n(\N_Y^*)$
as a sheaf on $D$ supported on $Y$.
Now we have short exact sequences
\begin{equation}\label{short}
0 \to \cF^{n+1}(\E_\gamma) \to \cF^n(\E_\gamma) \to 
	\cO(\E_\gamma|_Y \otimes S^n(\N_Y^*)) \to 0
\end{equation}
of sheaves on $D$ supported in $Y$.  This leads to the long
exact sequences
\begin{equation}\label{longD}
\begin{aligned}
0 \to &H^0(D;\cF^{n+1}(\E_\gamma)) \xrightarrow{a}
       H^0(D;\cF^n(\E_\gamma)) \xrightarrow{b} 
       H^0(D; \cO(\E_\gamma|_Y \otimes S^n(\N_Y^*))) \xrightarrow{\delta} \\
&H^1(D;\cF^{n+1}(\E_\gamma)) \xrightarrow{a}
       H^1(D;\cF^n(\E_\gamma)) \xrightarrow{b} 
       H^1(D; \cO(\E_\gamma|_Y \otimes S^n(\N_Y^*))) \xrightarrow{\delta} \\
&\dots\dots\dots\dots\dots\dots\dots\dots\dots\dots\dots\dots\dots\dots\dots
	\dots\dots\dots\dots\dots\dots\dots  \\
&H^{s-1}(D;\cF^{n+1}(\E_\gamma)) \xrightarrow{a}
       H^{s-1}(D;\cF^n(\E_\gamma)) \xrightarrow{b}
       H^{s-1}(D; \cO(\E_\gamma|_Y \otimes S^n(\N_Y^*)))\xrightarrow{\delta} \\
&H^s(D;\cF^{n+1}(\E_\gamma)) \xrightarrow{a}
       H^s(D;\cF^n(\E_\gamma)) \xrightarrow{b}
       H^s(D; \cO(\E_\gamma|_Y \otimes S^n(\N_Y^*))) \xrightarrow{\delta} 0.
\end{aligned}
\end{equation}
where $a$ and $b$ are coefficient morphisms from (\ref{short}),
$\delta$ is the coboundary, and $s = \dim_{_\C}Y$.  If a sheaf on a locally
compact space (such as $D$) is supported on a closed subspace (such as $Y$)
then the inclusion induces a natural isomorphism of cohomologies 
\cite[Corollary to Lemma 4.9.2]{G}.  So we can rewrite (\ref{longD}) as
\begin{equation}\label{longY}
\begin{aligned}
0 \to &H^0(Y;\cF^{n+1}(\E_\gamma)) \xrightarrow{a}
       H^0(Y;\cF^n(\E_\gamma)) \xrightarrow{b} 
       H^0(Y; \cO(\E_\gamma|_Y \otimes S^n(\N_Y^*))) \xrightarrow{\delta} \\
&H^1(Y;\cF^{n+1}(\E_\gamma)) \xrightarrow{a}
       H^1(Y;\cF^n(\E_\gamma)) \xrightarrow{b} 
       H^1(Y; \cO(\E_\gamma|_Y \otimes S^n(\N_Y^*))) \xrightarrow{\delta} \\
&\dots\dots\dots\dots\dots\dots\dots\dots\dots\dots\dots\dots\dots\dots\dots
        \dots\dots\dots\dots\dots\dots\dots  \\
&H^{s-1}(Y;\cF^{n+1}(\E_\gamma)) \xrightarrow{a}
       H^{s-1}(Y;\cF^n(\E_\gamma)) \xrightarrow{b}
       H^{s-1}(Y; \cO(\E_\gamma|_Y \otimes S^n(\N_Y^*)))\xrightarrow{\delta} \\
&H^s(Y;\cF^{n+1}(\E_\gamma)) \xrightarrow{a}
       H^s(Y;\cF^n(\E_\gamma)) \xrightarrow{b}
       H^s(Y; \cO(\E_\gamma|_Y \otimes S^n(\N_Y^*))) \xrightarrow{\delta} 0.
\end{aligned}
\end{equation}
Note that (\ref{longY}) is an exact sequence of $K$--modules.
\medskip

\begin{remark}\label{factorout} {\rm Let $n, j \geqq 0$. Then
$H^j(Y;\cO(\E_\gamma|_Y \otimes S^n(\N_Y^*))) = 
H^j(Y;\cO(\E_\gamma|_Y)) \otimes S^n(\gu_{-1})$  
as $K_2$--module.  If\, $\E_\gamma \to Y$ is a line bundle then 
$K_2$ acts trivially on the first factor $H^j(Y;\cO(\E_\gamma|_Y))$.}
\end{remark}
\noindent {\bf Proof.} The group $K_2$ acts trivially on $Y$, so
its action on $S^n(\N_Y^*)$ factors out of the cohomology.  
Recall that $\N_Y^* \to Y$ is the $K_0$--homogeneous vector bundle based on the 
$L_0$--module $\gu_{-1}$.
If $\E_\gamma \to Y$ is a line bundle then $K_2$ acts
trivially on each $H^j(Y;\cO(\E_\gamma|_Y))$ because it is semisimple.  
\phantom{XXXXXXXXXXXX}\hfill $\diamondsuit$
\medskip

Recall that the positive compact roots are those for which the coefficient of
$\nu$, as a linear combination from $\Psi = \Psi_G$, is $0$ or $2$.  The
ones of coefficient $0$ are roots of $(\gl,\gt)$.  The others, forming the
set $\Delta_2$ of the discussion after (\ref{lambdanought}), are the 
complementary compact positive roots.  They give the holomorphic tangent space 
of $Y$.  Let $\rho_\gk$ denote half the sum of the positive compact roots
(positive roots of $\gk$).  Then the proof of (\ref{delta1}) gives us
\begin{equation}\label{delta1k}
\langle \gamma + \rho_\gk, \alpha \rangle < 0 \text{ for all }
\alpha \in \Delta_2 \text{ if and only if }
\langle \gamma + \rho_\gk, \mu \rangle < 0.
\end{equation}

If $\alpha_1 \in \Delta_1$ and $\alpha_2 \in \Delta_2$ then 
$\alpha_1 + \alpha_2$ is not a root, because it would have coefficient 
$3$ at $\nu$.  Thus $\langle \alpha_1, \alpha_2 \rangle \geqq 0$.
That gives us
\begin{lemma}\label{compare-negativity}
If $\alpha_2 \in \Delta_2$ then 
$\langle \gamma + \rho_\gk, \alpha_2  \rangle \leqq
\langle \gamma + \rho_\gg, \alpha_2  \rangle$.  In particular
if $\alpha_2 \in \Delta_2$ then
$\langle \gamma + \rho_\gg, \alpha_2 \rangle < 0 \text{ implies }
\langle \gamma + \rho_\gk, \alpha_2 \rangle < 0$.
Thus the $G_0$--negativity condition\hfill\newline
\centerline{$\langle \gamma + \rho_\gg, \alpha \rangle < 0$ for all
$\alpha \in \Delta_1 \cup \Delta_2$}
implies the $K_0$--negativity condition\hfill\newline
\centerline{$\langle \gamma + \rho_\gk, \alpha \rangle < 0$ for all
$\alpha \in \Delta_2$.}
\end{lemma}

We are going to need the following fact about tensor products of
irreducible finite dimensional representations.  It appears in \cite{H}
as Exercise 12 to Section 24, based on \cite{Ko1}.
\begin{lemma}\label{weights-tensor}
Let $E_{\gamma_1}$ and $E_{\gamma_2}$ be irreducible $L_0$--modules, where
$\gamma_i$ is the highest weight of $E_{\gamma_i}$.  Then every
irreducible summand of $E_{\gamma_1}\otimes E_{\gamma_2}$ has highest
weight of the form $\gamma_1 + \varphi$ for some weight $\varphi$ of
$E_{\gamma_2}$.
\end{lemma}

Now the $K_0$--negativity condition gives a vanishing result in
(\ref{longY}), as follows, where we take (\ref{delta1k}) into account.

\begin{theorem}\label{vanishing}
Suppose that $\langle \gamma + \rho_\gk, \mu \rangle < 0$.  Then
$H^j(Y;\cO(\E_\gamma|_Y \otimes S^n(\N_Y^*))) = 0$ whenever $j \ne s$
and $n \geqq 0$.
\end{theorem}
\noindent {\bf Proof.}  
Note that $\E_\gamma|_Y \otimes S^n(\N_Y^*) \to Y$ is the $K_0$--homogeneous
bundle based on the representation of $L_0$ on 
$E_\gamma \otimes S^n(\gu_{-1})$.  In view of Lemma \ref{weights-tensor}
that $L_0$--module is the sum of
irreducibles with highest weights of the form 
$\gamma + \varphi$ where $\varphi$ is a weight of $S^n(\gu_{-1})$.  Thus,
as a homogeneous 
holomorphic vector bundle, $\E_\gamma|_Y \otimes S^n(\N_Y^*)$ has
composition series with composition factors of the form
$\E_{\gamma + \alpha_1 + \dots + \alpha_n}$ where the 
$\alpha_i \in \Delta_{-1}$.
\medskip

Let $\alpha \in \Delta_2$.  Then (\ref{delta1k}) shows that 
$\langle \gamma + \rho_\gk, \alpha \rangle < 0$.
The coefficient of $\nu$ in $\alpha$ is $2$, so $\alpha - \alpha_i$ cannot
be a root. This forces $\langle \alpha_i, \alpha \rangle \leqq 0$.
Now $\langle \gamma + \rho_\gk, \alpha \rangle < 0$ forces
$\langle \gamma + \alpha_1 + \dots + \alpha_n + \rho_\gk, \alpha \rangle < 0$.
The Bott--Borel--Weil Theorem now tells us that 
$H^j(Y;\cO(\E_\gamma|_Y \otimes S^n(\N_Y^*))) = 0$ for $j \ne s$.
\hfill $\square$
\medskip

\begin{corollary}\label{lg-vanishing}
Suppose that $\langle \gamma + \rho_\gk, \mu \rangle < 0$.  Then
$H^j(Y;\cO(\E_\gamma|_Y)) = 0$ whenever $j \ne s$.
\end{corollary}

Following \cite{S2} and \cite{W4}, with the result \cite{S3}
that $\gamma + \rho_\gg$ need only be nonsingular (instead of
``sufficiently nonsingular''), one has the following vanishing theorem.

\begin{theorem}\label{g-vanishing}
If $\langle \gamma + \rho_\gg, \alpha \rangle < 0$ whenever 
$\alpha \in \Delta_1 \cup \Delta_2$, then
$H^j(D;\cO(\E_\gamma)) = 0$ for $j \ne s$.
\end{theorem}

Recall the decomposition of (\ref{lambdanought}): 
$\gamma = \gamma_0 + t\nu^* \text{ where }
        \langle \gamma_0 , \nu \rangle = 0 \text{ and } t \in \R$.
In view of (\ref{delta1}), (\ref{delta2}) and Theorem \ref{delta}, 
we reformulate Theorem \ref{g-vanishing} as follows.

\begin{theorem}\label{g-vanish-range}
Let $\gamma = \gamma_0 + t\nu^*$ as in {\rm (\ref{lambdanought})}.
If $\langle \gamma + \rho_\gg, \mu\rangle < 0$ and
$\langle \gamma + \rho_\gg, w_\gl^0(\nu)\rangle < 0$, in other words if
$t < -\tfrac{1}{2}\langle \gamma_0 + \rho_\gg, \mu \rangle
\text{ and } t < -\langle \gamma_0 + \rho_\gg, w^0_\gl(\nu)\rangle$,
then $H^j(D;\cO(\E_\gamma)) = 0$ for $j \ne s$.
\end{theorem}

\begin{definition}\label{def-suff-neg} {\rm
To facilitate use of these vanishing theorems we will say that
$\E_\gamma \to D$ is {\em sufficiently negative} if
$\langle \gamma + \rho_\gg, \alpha \rangle < 0$ whenever
$\alpha$ is a complementary positive root, i.e. whenever
$\alpha \in \Delta_1 \cup \Delta_2$. 
This means that $\E_\gamma \otimes \K^{1/2} \to D$ is 
negative in the sense of differential or algebraic geometry, where
$\K \to D$ is the canonical line bundle.
}\hfill $\diamondsuit$
\end{definition}
In the presence of sufficient negativity Theorem \ref{g-vanish-range}
trivializes the long exact sequences (\ref{longD}) and (\ref{longY}) 
as follows.

\begin{proposition}\label{long_suff_neg}
Suppose that $\E_\gamma \to D$ is sufficiently negative.  Then 
$$
\begin{aligned}
& H^q(Y;\cF^{n+1}(\E_\gamma)) \cong H^q(Y;\cF^n(\E_\gamma))
\phantom{XX} \text{ for } 0 \leqq q < s, \phantom{XX} \text{ and }  \\
&H^s(Y; \cO(\E_\gamma|_Y \otimes S^n(\N_Y^*))) \cong
H^s(Y;\cF^n(\E_\gamma))/H^s(Y;\cF^{n+1}(\E_\lambda)).
\end{aligned}
$$
\end{proposition}
\bigskip

Now we may apply the above case-by-case analysis and diagrams
to understand what amounts to the structure 
and geometric quantization of the coadjoint
elliptic orbits corresponding to the particular discrete series of
representations we have in mind, the so--called {\it Borel -- de
Siebenthal discrete series}. Also, the results above on
the filtration are crucial for the construction of the
cohomology groups carrying these representations; they are    
the analytic counterparts of the Vogan--Zuckerman derived functor
modules that are constructed purely algebraically. Thus  
we wish to construct the Borel -- de Siebenthal discrete
series by direct analysis on orbits,
and using the above results analyze the $K_0$--types
explicitly (without subscript $K$
denotes the complexified group); after this we shall end the paper with some
remarks and immediate consequences, and treat the
(interesting) analytic continuation of this particular
discrete series in
a later paper.
\medskip

As is clear from the above discussion there is some variation in the meaning 
of ``discrete series''.  Initially the discrete series of $G_0$ meant the
family of (equivalence classes of) irreducible unitary representation $\pi$
of $G_0$ that are discrete summands of the left regular representation.  This
is equivalent to the condition that the matrix coefficients
$f_{u,v}(g) = \langle u, \pi(v)\rangle$ of $\pi$ belong to $L^2(G)$.  That
is how they are treated in the work of Harish--Chandra, and there 
the discrete series representations are also treated as Harish--Chandra 
modules.  Later one had the construction of discrete series representation 
as the action of $G_0$ on cohomology spaces $H^q(D;\E)$ both as nuclear
Fr\' echet spaces (\cite{S2}, \cite{SW2}) and as Hilbert spaces \cite{W4}, 
and still later they appeared algebraically as
Zuckerman derived functor modules.  The underlying Harish--Chandra module
is the same for all these constructions, and we will use the
cohomology constructions.  
\medskip

We first recall some results about the discrete series
representations in general. See \cite[Theorem 9.20]{Kn1} 
where Harish-Chandra's parameterization is recalled: Here
there is given a {\it standard root order} (which is not the same
as we are working with in the diagrams above) of the
root system as follows:
\begin{equation} \label{standard_order}
\Delta^+_\lambda = \{\alpha \in \Delta \mid \langle \lambda, \alpha \rangle > 0 \} 
\end{equation}
where $\Delta$ is the root system, and $\lambda$ is the {\it Harish-Chandra
parameter} for the  discrete series representation $\pi_{\lambda}$.
The Harish--Chandra parameter $\lambda \in (i \gt)'$ satisfies the integrality
condition that $\lambda + \rho_\gg$ is analytically
integral, in other words that $\exp(\lambda + \rho_\gg)$ is a well defined
character on the maximal torus of $K_0$.  It also satisfies the 
nonsingularity condition that
$\langle \lambda, \alpha \rangle \neq 0$ for all $\alpha$ in
$\Delta$. Two such representations are equivalent if and only
if their parameters are conjugate under the compact Weyl group $W_\gk$.
Thus one could normalize the Harish--Chandra parameter by the condition
that $\langle \lambda, \alpha \rangle <0$ for all compact positive roots.
\medskip

The parameter $\lambda$ of course determines the positive root system
$\Delta^+_\lambda$ of the standard root order
(\ref{standard_order}), and conversely to
each Weyl chamber of $\gg$, modulo the action of $W_\gk$,
we associate a family of
discrete series representations. The family we are interested in
is in some sense the smallest possible kind of discrete series
representations of $G_0$.
\medskip

Of particular interest is the lowest $K_0$--type contained in the 
(Harish-Chandra module for) $\pi_{\lambda}$ given by its highest
weight (in the standard root order (\ref{standard_order}))
$$ 
\Lambda = \lambda + \rho_\gg - 2\rho_\gk
$$      
in terms of the usual half sums of positive roots. This $K_0$--type
has multiplicity one, and other $K_0$--types have highest weights of the form
$$
\Lambda' = \Lambda + \sum_{\alpha \in \Delta^+} n_{\alpha} \alpha
$$
for integers $n_{\alpha} \geqq 0$.
In the general theory of discrete series this statement
about the $K_0$--types only amounts to an inclusion, whereas our
results above analyzing the cohomology groups in terms of
restriction and Taylor expansion in the normal direction ($V$)
gives a concrete list of the $K_0$--types. We shall formulate this
precisely below.
\medskip
 
Let us first see how these parameters fit with the description
in \cite{GW} of the quaternionic discrete series $\pi^q_{\lambda}$. 
They write $\beta$ for the maximal root, but we translate that to
our notation of $\mu$ in describing their results.  Thus the  
Harish--Chandra (and infinitesimal character) parameter of their 
$\pi^q_{\lambda}$ is of the form
$\lambda = - \frac{k}{2} \mu + \rho_\gg$
where the integer $k \geq 2d + 1$ and $\dim_\R G_0/K_0 = 4d$.  We
consider the corresponding standard root order $\Delta^+_\lambda$.
Dividing as usual into compact and noncompact roots we have
$\rho_\gg = \rho_\gk + \rho_{\gg/\gk}$, 
$\rho_\gk = \rho_\gl + \frac{\mu}{2}$ and
$\rho_{\gg/\gk} = \frac{d}{2} \mu$,
where $\rho_{\gg/\gk}$ is half the sum of the noncompact positive roots
and $\rho_\gl$ is half sum of positive roots of $\gl$.  Similarly for the
standard root order we have
$\rho'_\gk = \rho_\gl - \frac{\mu}{2}$ and
$\rho'_{\gg/\gk} = -\frac{d}{2} \mu$,
so the lowest $K_0$--type in the standard root order $\Delta^+_\lambda$ is
$\Lambda = -\frac{k}{2} \mu + \rho_\gg + \rho'_\gg - 2 \rho'_\gk$.  That
simplifies to $\Lambda = \frac{-k+2}{2} \mu$.
This is exactly the highest weight for the $(k-1)$--dimensional
representation of the simple $SU(2)$ factor in $K$ found
as the lowest $K_0$--type by Gross and Wallach. In the following
we shall find the analogous lowest $K_0$--type for the
Borel -- de Siebenthal discrete series, and at the same time
realize it (and in fact all $K_0$--types) as cohomology groups
on the compact Hermitian symmetric space $Y$.        
\medskip
   
Now recall the noncompact simple root $\nu$ from Section \ref{sec2}.
As before, $\nu^*$ denote the dual to $\nu$ in the system of
fundamental simple weights (\ref{fund-wts}).
The parabolic subalgebra $\gq$ of $\gg$ may also be defined by
means of $\nu^*$, and the centralizer of $\nu^*$ is $\gl$. Thus
the coadjoint orbit $\Ad^*(G_0)(\nu^*)$ is our space $G_0/L_0$ 
and is fibered by $Y$.   
Multiples of this $\nu^*$ will define the line bundles we shall
need, and the representations in the Borel -- de Siebenthal
discrete series are then the cohomology groups in degree
$s = \dim_\C Y$ with coefficients in the bundle.
\medskip

Recall the maximal compact subgroup $K_0 = K_1 \times K_2$
explicit in the classification of Section \ref{sec3}, where 
the ``small'' factor $K_1$ corresponds to the component of
the simple root system $\Psi_\gk = (\Psi \setminus \{\nu\}) \cup \{-\mu\}$
that contains $\{-\mu\}$.  In the quaternionic case  $L_1 = Sp(1)$.
Now $Y = K_0/L_0 = (K_1 \times K_2)/(L_1 \times K_2) = K_1/L_1$.
Thus, as far as induced representations and cohomology, the
action of the $K_2$ factor will be rather simple.
This we will make explicit below.  Also, it is important that
the factor $L_1$ in $L_0$ contains the center of $L_0$, and that the
action of that center on the holomorphic normal space $V = \gu_1$
is given explicitly in the case by case diagrams of Section \ref{sec2}.
\medskip
               
Let $\E_{\gamma_k} \to D$ be the holomorphic vector bundle induced from the
representation of $L_0$ with highest weight
$\gamma_k = \gamma_0 - k\nu^*, \, k \in \N$.  (As $G$ is simply 
connected $\exp(\gamma_k)$ is the highest weight of a representation
of $L_0$).  Denote
\begin{equation}\label{defpik}
\pi_{\lambda_k}: \text{ representation of } G_0 \text{ on } 
H^s(D, \cO(\E_{\gamma_k})) \text{ where } \lambda_k = \gamma_k + \rho_\gg = 
\gamma_0 -k \nu^* + \rho_\gg\, .
\end{equation}
We will say that the integer $k$ is {\em sufficiently positive} if
$\lambda_k = \gamma_k + \rho_\gg = \gamma_0 -k \nu^* + \rho_\gg$ is 
sufficiently negative in the sense of Definition \ref{def-suff-neg}
and Theorem \ref{delta}.
Using the filtration (\ref{filtration}) and arguments analogous to those
of \cite{GW} we characterize the line bundle valued Borel -- de Siebenthal
discrete series as follows.

\begin{theorem}\label{bds_line_ds}
The Borel -- de Siebenthal discrete series representations
of $G_0$ are the $\pi_{\lambda_k}$ of {\rm (\ref{defpik})} for which 
$k \in \N$ is sufficiently positive.  As a $K$--module, the underlying
Harish--Chandra module is 
$$
\sum_{m \geqq 0} H^s(Y, \cO(\E_{\gamma_k } \otimes S^m(V^*))),
$$
and it not only is $K_0$--admissible but is $K_1$--admissible.  Further,
the lowest $K_0$--type {\rm (corresponding to $m = 0$)} is given by
$$
W_{\lambda_k} = H^s(Y, \cO(\E_{\gamma_k})),
$$
and it has multiplicity $1$ in $\pi_{\lambda_k}$.
\end{theorem}

\noindent {\bf Remark.}  The $L_0$--modules  $S^m(V^*)$ are not 
always multiplicity free, though they are multiplicity free in
many cases.  For example for the group of type $D_9$ and $m = 6$, 
calculation with the computer program LiE produces multiplicities,
while there are none for $F_4$.  Thus even in the scalar case,
where $\E_{\gamma_n} \to D$ is a line bundle, i.e. when $\gamma_0 = 0$,
$\pi_{\lambda_k}$ need not be $K_0$--multiplicity free.  This is of course 
in contrast the the $K_0$--multiplicity free property of the 
line bundle holomorphic discrete series. \hfill $\diamondsuit$
\medskip

\noindent{\bf Proof.} We use the filtration (\ref{filtration}) and the
exact sequences (\ref{longD}) and (\ref{longY}), together with the fact
that $Y$ is a compact hermitian symmetric space for $K_1$.  The action of $K_2$
is part of the holomorphically induced representation, and the action of 
$L_1$ on $V$ and its dual $V^*$  is given as above. Finally the
admissibility can be read off from the $K_0$--types directly:
each $H^s(Y, \cO(\E_{\gamma_k} \otimes S^m(V^*)))$ is a 
sum of irreducible representations 
of $K_1$, disjoint for different $m$, and the $S^m(V^*)$ are
finite dimensional representations of $K_2$ and also
of $L_0$.  
\medskip

Consider the parabolic subgroup $Q\cap K = LU_{-2}$ of $K$.  Whenever
$M$ is a finite dimensional $(Q\cap K)$--module, the space of
$K$--finite vectors in the induced representation
$\Ind_{Q\cap K}^K(M)$ is $\sum_{\delta \in \widehat{K}} 
	V_{\delta} \otimes (V_{\delta}^* \otimes M)^{Q\cap K}$.
In particular the multiplicity of a $K_0$--type $\delta$ is
equal to the number of times the highest weight vector of $M$ occurs
as a highest weight vector for $L$ in $V_{\delta}$. In our case 
the highest weights of $M$ will grow with $m$ in $S^m(V^*)$, and
they are distinguished by the action of the center of $L_0$. Thus
each $K_1$--type only occurs finitely many times, so $\pi_{\lambda_k}$
is $K_1$--admissible.  That, of course, implies admissibility for $K_0$.
\hfill $\square$
\medskip

\noindent{\bf Remark.}  We compare our parameter for the lowest $K_0$--type
with the general description mentioned for the scalar quaternionic case.
In that scalar quaternionic case
the infinitesimal character of the representation $\pi_{\lambda_k}$
is given by 
$$
\lambda_k = -k\nu^* + \rho_\gg, \, \rho_\gk = \rho_\gl + c_2 \nu^* 
\,\,\text{ and } \,\,\rho_{\gg/\gk} = c_1 \nu^*
$$
for positive constants $c_1$ and $c_2$ depending only on the root system.
Then for the standard root order $\Delta^+_{\lambda_k}$ we get
$$
\rho'_\gk = \rho_\gl - c_2 \nu^* \text{ and } \rho'_n = -c_1 \nu^*
$$
so that the lowest $K_0$--type has highest weight 
$$ 
\Lambda = -k\nu^* + \rho_\gg + \rho'_\gg - 2 \rho'_\gk 
	= -k\nu^* + 2\rho_{\gk/\gl}
$$
where $\rho_{\gk/\gl} = c_2 \nu^*$ is exactly the shift coming from
the square root of the canonical bundle $\K \to Y$.  Thus this corresponds 
to the lowest $K_0$--type above, viz. $ W = H^s(Y, \cO(\L_{-k}))$.              
\hfill $\diamondsuit$
\medskip

It is an interesting problem to study the structure, including
unitarity, of $\pi_{\lambda_k}$ for smaller values of $k$, and to
relate this to the projective varieties defined by the
relative invariants - this will be taken up in a sequel to
this paper. In particular the ring of
regular functions on $L_0$--orbits will be important,
as in the paper by Gross and Wallach for the case
of the quaternionic discrete series. 
For now we remark as an application of the
admissibility above, that branching problems
will be manageable in a way
similar to the case of holomorphic discrete series.
This will require that the embedding of the smaller group
respects the relevant structure, i.e. that the orderings are
compatible. For example if we want to branch to a symmetric
subgroup $H_0$ of $G_0$, then the embedding will be compatible provided
the symmetry fixes $K_1$.  
Namely, we simply use the admissibility of the action of
$K_1$, so that admissibility for the branching law to
$H_0$ will follow for the Borel -- de Siebenthal discrete series.  
In this case a Borel -- de Siebenthal discrete series will branch
as a direct sum of Borel -- de Siebenthal discrete series
representations.
\medskip

\noindent{\bf Remark.}
As is evident from the case of indefinite orthogonal groups
as in \cite{Kn3}, the question of continuation of the
discrete series modules is closely connected with the
geometry of the relative invariants for the holomorphic
normal (to the maximal compact
subvariety) $V$. Already the case where $G$ is of type     
$E_8$ and $K$ of type $D_8$ is an interesting example; here
$V$ is of dimension $64$ and admits a relative invariant
of degree $8$, and the maximal compact subvariety is
the Grassmannian of $2$--planes in $16$--space. 
Let us here be a little more explicit about this example:

Let 
$$\gamma_0 = n_2\xi_2 + ... + n_8\xi_8$$ 
and we find 
       $$\rho_\gk = 14\psi_1 + 28\psi_2 + 35\psi_3 + 55\psi_4
                  + 46\psi_5 + 36\psi_6 + 25\psi_7 + 13\psi_8$$
and
        $$\rho_\gg =  46\psi_1 + 68\psi_2 + 91\psi_3 + 135\psi_4
                  + 110\psi_5 + 84\psi_6 + 57\psi_7 +  29\psi_8$$
so the "sufficient $G_0$--negativity" condition of Theorems 2.12
and 5.18 is, using 
$$ \mu = 2\psi_1 + 3\psi_2 + 4\psi_3 + 6\psi_4
         + 5\psi_5 + 4\psi_6 + 3\psi_7 + 2\psi_8$$
that both
$t < -\frac{1}{2}(3n_2 + 4n_3 + 6n_4 + 5n_5 + 4n_6 + 3n_7
+ 2n_8) - \frac{29}{2}$ 
and
$t < -\langle \gamma_0 + \rho_\gg, w_\gl^0(\nu)\rangle$. 
Now the range from $G_0$--sufficiently negative to $K_0$--sufficiently negative
is indicated by the condition
(as in Corollary 5.16) $t < -\frac{1}{2}(3n_2 + 4n_3 + 6n_4 + 5n_5 + 4n_6 + 3n_7
+ 2n_8) + \frac{1}{2}$.

Hence we see that there is an interval, where the $K_0$--types still
exist as cohomology groups, even though the large cohomology group
carrying the $G_0$--representation ceases to exist.       
We shall study in more detail what happens here in a sequel to the
present paper.

\enddocument
\end